\documentclass{article}

\usepackage{arxiv}
\usepackage[utf8]{inputenc}
\usepackage[T1]{fontenc}    

\usepackage{graphicx}
\usepackage{subfig}       
\graphicspath{{figures/}} 

\usepackage{psfrag}       
\usepackage{amsmath,bm,amsfonts,mathrsfs}   
\usepackage{shortcuts}    
\usepackage{color}

\usepackage{hyperref}       
\usepackage{url}            
\usepackage{enumitem}
\usepackage{soul} 
\usepackage[normalem]{ulem} 

\usepackage{algorithm,algorithmic}

\usepackage{scalerel}
\DeclareMathOperator*{\foo}{\scalerel*{A}{\sum}}

\usepackage{xr}
\usepackage[color=green!20]{todonotes}

\title{Label-free learning of elliptic partial differential equation solvers with generalizability across boundary value problems}

\author{
  Xiaoxuan Zhang$^1$, ~Krishna Garikipati$^{1,2,3}$ \thanks{Corresponding author. E-mail address: krishna@umich.edu \hfill \today} \\[3mm]
  $^1$Department of Mechanical Engineering, University of Michigan, United States \\
  $^2$Department of Mathematics, University of Michigan, United States \\
  $^3$Michigan Institute for Computational Discovery \& Engineering, University of Michigan, United States \\
}

\begin{document}

\maketitle

\begin{abstract}

Traditional, numerical discretization-based solvers of partial differential equations (PDEs) are fundamentally agnostic to domains, boundary conditions and coefficients. In contrast, machine learnt solvers have  a limited generalizability across these elements of boundary value problems. This is strongly true in the case of surrogate models that are typically trained on direct numerical simulations of PDEs applied to one specific boundary value problem. In a departure from this direct approach, the label-free machine learning of solvers is centered on a loss function that incorporates the PDE and boundary conditions in residual form. However, their generalization across boundary conditions is limited and they remain strongly domain-dependent. Here, we present a framework that generalizes across domains, boundary conditions and coefficients simultaneously with learning the PDE in weak form. Our work explores the ability of simple, convolutional neural network (CNN)-based encoder-decoder architectures to learn to solve a PDE in greater generality than its restriction to a particular boundary value problem. In this first communication, we consider the elliptic PDEs of Fickien diffusion, linear and nonlinear elasticity. Importantly, the learning happens independently of any labelled field data from either experiments or direct numreical solutions. We develop probabilistic CNNs in the Bayesian setting using variational inference. Extensive results for these problem classes demonstrate the framework's ability to learn PDE solvers that generalize across hundreds of thousands of domains, boundary conditions and coefficients, including extrapolation beyond the learning regime. Once trained, the machine learning solvers are orders of magnitude faster than discretization-based solvers. They therefore could have relevance to high-throughput solutions of PDEs on varying domains, boundary conditions and coefficients, such as for inverse modelling, optimization, design and decision-making. We place our work in the context of recent continuous operator learning frameworks, and note extensions to transfer learning, active learning and reinforcement learning.


\end{abstract}



\section*{Introduction}

Partial differential equation (PDE) solvers play a central role in computational science and engineering. They bridge between the mathematical physics of field theories to applications in engineering science.
Popular, discretization-based numerical methods to solve PDEs include, but are not limited to, the finite element method (FEM), finite difference method, finite volume method and their variants, each with its own advantages and limitations.
The FEM, in its many variant forms, is notable for the natural treatment of complex domain geometries and boundary conditions. 
However, when a large number of boundary value problems need to be solved, such as for inverse modelling, optimization, design and decision-making or these discretization-based numerical solvers can prove very expensive. Scientific machine learning (ML) techniques have proved to be natural candidates.

ML approaches to solving mathematical descriptions of physical systems can be categorized as surrogate models and PDE solvers. 
The first category typically requires a vast amount of training data, either from measurement or direct numerical simulations (DNSs), whose acquisition can pose challenges of availability and expense, (see \cite{Zhu2018Zabaras-UQ-Bayesian-ED-CNN,Winovich2019Lin-ConvPDE-UQ,Bhatnagar2019CNN-encoder-decoder,Li2020Reaction-diffusion-prediction-CNN}, and many others).
For example, in Ref \cite{Zhu2018Zabaras-UQ-Bayesian-ED-CNN} a Bayesian uncertainty quantification (UQ) approach to convolutional neural networks (CNNs) was proposed for flows in heterogeneous media. 
CNNs are also used to predict the velocity and pressure fields in aerodynamics \cite{Bhatnagar2019CNN-encoder-decoder} and the concentration field for single-species reaction-diffusion systems \cite{Li2020Reaction-diffusion-prediction-CNN}.
The second category requires little or no pre-labeled data to solve PDEs \cite{Lagaris1998NN-PDEs,Han2018Solve-high-dimension-PDE-deep-learning,Sirignano2018DGM-solve-PDEs,Raissi2019physics-informed-forward-inverse-jcp,Zhu2019Perdikaris-Physics-PDE-CNN,Geneva2020Zabaras-JCP-auto-regressive-NN-PDE,Yang2021BPINNs-PDE,Berg2018unified-deep-ann-PDE-complex-geometries,Sun2020Surrogate-modeling-flow-physics-constrained,Samaniego2020energy-approach-PDE-ML}.
For example, high-dimensional, free-boundary PDEs have been solved by fully connected NNs \cite{Han2018Solve-high-dimension-PDE-deep-learning}, and by the Deep Galerkin Method  \cite{Sirignano2018DGM-solve-PDEs} using NNs, which satisfy the differential operators, initial conditions (ICs), and boundary conditions (BCs).
The Physics-informed Neural Networks (PINNs) approach has been proposed to solve steady and transient systems \cite{Raissi2019physics-informed-forward-inverse-jcp}.
In PINNs, the strong form of PDEs, the ICs and BCs are incorporated in the loss \cite{Raissi2019physics-informed-forward-inverse-jcp}.
PINNs have been extended to solve numerous systems \cite{Jin2020NSFnets-PINN,Pang2019fPINNs,Meng2020PPINN,Geneva2020Zabaras-JCP-auto-regressive-NN-PDE,Yang2021BPINNs-PDE,Wang2020DL-free-boundary,Jagtap2020cPINN-discrete-domain-conservation-law,meng2022learning}.
NN-based PDE solvers constructed from the weak/variational formulation  also have been studied \cite{Zang2020Weak-adversarial-network-PDE,Chen2020Friedrichs-learning-weak-PDE,Khodayi-mehr2019VarNet,Li2020D3M-domain-decomposition,Kharazmi2021hp-VPINNs}.

To become viable alternates to discretization-based PDE solvers, ML frameworks have to extend to solutions of the same PDE system, but with different ICs, BCs, and on different problem domains.
However, this is difficult to achieve with NN-based approaches that typically enforce only one specific set of BCs \cite{Raissi2019physics-informed-forward-inverse-jcp}, parameterized BCs and domains \cite{Gao2020PhyGeoNet-PDE-on-Irregular-domain} via the loss function \cite{Raissi2019physics-informed-forward-inverse-jcp} or the NN architecture \cite{Sun2020Surrogate-modeling-flow-physics-constrained,Gao2020PhyGeoNet-PDE-on-Irregular-domain}.
Such NN-based solvers must be retrained for each new BC and domain. In one recent approach, the notion of a genome has been introduced to learn solutions on subdomains and use them to construct solutions on larger, and to some extent varying shape domains \cite{wang2022mosaic}. However, labelled training data are needed. A related domain-decomposition approach with NNs to impose desired regularity at subdomain boundaries has also been used \cite{dong2021local}.

In this work, we address these challenges by a new class of physics-constrained NN solvers where the BCs are specified as inputs.
We draw from the FEM, where the weak formulation incorporates the governing PDE, the natural and essential BCs, and the solution corresponds to the vanishing the discretized residual. Central to our framework is transformation of the NN predicted solution to a discretized PDE residual that defines the loss function used to train the NNs, through efficient, convolutional operator-based, and vectorized calculations.
We introduce weak PDE loss layers via kernels whose parameters are not trainable, and are independent of the NN that learns the PDE solver.
Such features offer us great flexibility to choose the NN architecture. 

In our framework, the trainable NN learns, from a number of boundary value problems, a representation that recognizes domains, boundary conditions and coefficients. To the extent that this learning is imperfect, there will be uncertainties in its predicted solutions. 
One common categorization of uncertainties in modelling frameworks is as epistemic and aleatoric. The former represents model-form error that can be reduced by learning from more data or using a better model--and is natural for any finite-capacity NN. The latter typically represents measurement errors, is less prone to reduction \cite{Kiureghian2009Ditlevsen-aleatory-or-epistemic}--and is not applicable to the proposed framework, in which boundary value problems are exact statements. 
Probabilistic machine learning models have been developed for uncertainty quantification with NN-based PDE solvers \cite{Karniadakis2019Zhang-Quantify-UQ-dropout,Yang2019Perdikaris-Adversarial-UQ,Geneva2020Zabaras-JCP-auto-regressive-NN-PDE,Yang2021BPINNs-PDE}. Techniques, such as dropout \cite{Gal2016Dropout}, adversarial inference \cite{Yang2019Perdikaris-Adversarial-UQ}, Bayesian methods \cite{Yang2021BPINNs-PDE}, Stochastic Weight Averaging Gaussian \cite{Maddox2019SWAG} have been used, among many others. A recent work \cite{kadeethum2021framework} uses conditional generative adversarial networks to map between images for forward and inverse solution of PDEs; the authors treatment of domains as images bears similarity to our approach.
In this communication, we present both deterministic and probabilistic ML-PDE solvers, with Bayesian NNs (BNNs) for the latter. More specifically, we adopt variational inference  in the Bayesian setting, motivated by its efficiency over Monte Carlo approaches.
We focus on an encoder-decoder architecture, which has been investigated for other physical systems \cite{Zhu2018Zabaras-UQ-Bayesian-ED-CNN,Winovich2019Lin-ConvPDE-UQ,Bhatnagar2019CNN-encoder-decoder}.
The encoder-decoder structure can be easily adapted to BNNs with the modularized probabilistic layers available in current ML platforms, of which we use the TensorFlow Probability (TFP) library.
In our framework, deterministic/probabilistic convolutional NN layers learn representations for problem domains and the applied BCs (both Dirichlet and Neumann) through carefully designed input data structures. 
Thus, a single NN can be used to simultaneously solve different boundary value problems that are governed by the same PDEs but on different domains with different BCs.
Having learnt the PDE, the ML solvers can make predictions for interpolated and, to a certain extent, extrapolated domains and BCs that they were not exposed to during training. Similar to other NN-based PDE solvers, our learning approach is free of labelled data on field solutions.

We note the recent development of operator learning approaches for nonlinear mappings between continuous functional spaces as inputs and solution fields as outputs. Of particular interest in this direction are DeepONets \cite{lu2021learning,goswami2022deep} as well as graph kernel networks \cite{li2020neural} and Fourier neural operators \cite{li2020fourier} as settings for PDE solvers. Our framework differs from these approaches in its focus upon learning solvers for specific PDEs, but  with generalizability across boundary value problems spanning domains, BCs and coefficients with the added feature of UQ. Our proposed framework is generalizable and applicable to both steady-state and transient problems.
Here, we present it in detail and focus on its application to elliptic PDEs of steady-state diffusion, linear and nonlinear elasticity.
We defer the investigation of transient problems to a subsequent work.
We demonstrate that, with the proposed framework, a single NN can learn a solver that can be applied to tens to hundreds of thousands of boundary value problems. For brevity we use NN-PDE-S for the deterministic neural network-based PDE solver, and BNN-PDE-S for the Bayesian version.

\section*{Results}

\subsection*{(Bayesian) NN-based PDE solver}

\begin{figure}[t!]
    \centering
    \includegraphics[width=1.0\linewidth]{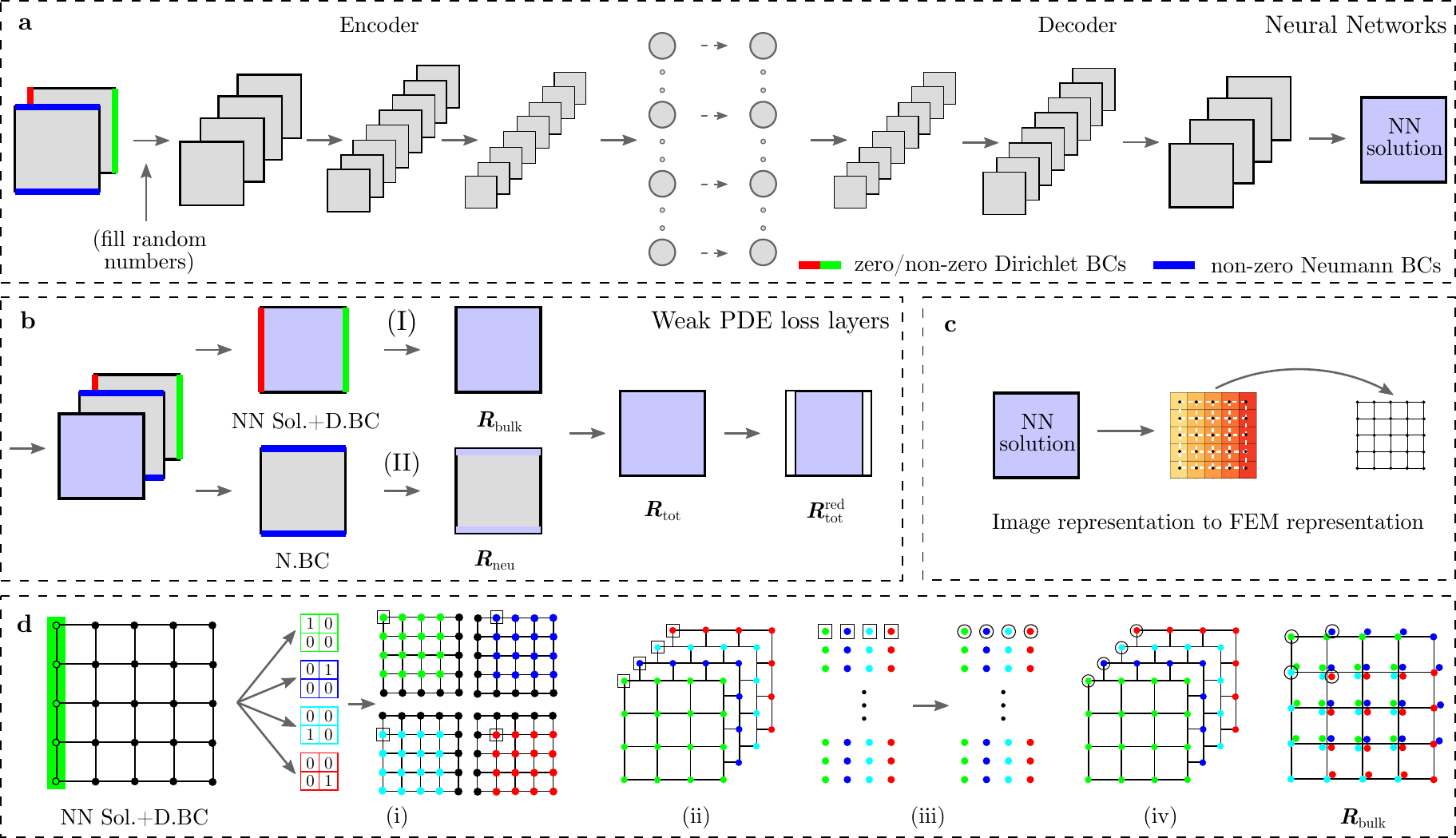}
    \caption{Architecture of the NN-PDE-S. 
      (a) An encoder-decoder NN to store the nonlinear mapping between NN inputs, which consist of the geometry of problem domains and the applied BCs, and the solution of PDEs. Random numbers are given to the pixels in the interior domain of the channel with Dirichlet BCs when working with the augmented dataset with only a few unique boundary value problems. BCs are color coded with red for zero Dirichlet BCs; green for non-zero Dirichlet BCs; blue for non-zero Neumann BCs.  
      (b) Both NN inputs and outputs are passed to the Weak PDE loss layers to calculate the discretized residual, where the bulk residual (I) is computed from NN solutions with imposed Dirichlet BCs and the residual contribution from the Neumann BCs (II) are computed based on NN inputs. The reduced discretized residual by excluding contributions from Dirichlet boundary locations is used to form the loss of both deterministic and probabilistic NNs. 
      (c) Illustration of the construction of finite element meshes from pixelated image representations.
      (d) Details of the bulk residual calculation (I). Four filters are applied to the NN solutions with imposed Dirichlet BCs to select different nodes (i), resulting in a multi-channel data representation (ii), which has the structure of the local nodes of one finite element. The multi-channel data is reshaped into a two-dimensional matrix to perform residual calculation (iii). The two-dimensional residual matrix is reshaped to a multi-channel data structure (iv) and x reduced sum operations are performed to get the bulk residual,  avoiding the time-consuming assembly process in the traditional FEM.
}
    \label{fig:NN}
\end{figure}

In the proposed PDE solver (Fig. \ref{fig:NN}), NNs are used to represent the nonlinear mapping between BCs and the resulting solutions of PDEs. 
The discretized residuals of PDEs are used to construct the losses and therefore regularize the NN's solutions of PDEs.
We studied both deterministic NNs and BNNs, where the uncertainty of the latter is represented by computing the statistical moments of their outputs via the predictive expectation and the predictive variance.

Deterministic NNs have fixed model parameter  values, and their losses  are mean squared Euclidean norms of the discretized reduced residual vectors.
For BNNs, the model parameters are drawn from a posterior distribution that is computed from Bayes' theorem.
The loss of BNNs is formulated using variational inference  \cite{Blei2017Variational-inference-review}, which consists of a data-independent contribution and a data-dependent contribution. The latter is the log likelihood function, which has the form of a joint distribution of the discretized residual with an added Gaussian noise.
Both NNs are trained via a mini-batch optimization process with standard stochastic optimizers.

In the proposed framework, the PDE loss layers (Fig. \ref{fig:NN}b) are independent of the NNs, which offers  flexibility in choosing the NN architecture. 
An encoder-decoder NN architecture is explored (Fig. \ref{fig:NN}a).
The NNs, which store the nonlinear mappings between BCs and the PDE solutions, accept image-type inputs that contain physically meaningful boundary values and markers for different regions to allow the convolutional NN layers to learn the BCs and problem domain. 
This allows a well-trained NN to make predictions for new problem domains with new BCs when these are provided as inputs.
As the image-type NN outputs can be treated as  FEM meshes (Fig. \ref{fig:NN}c), we evaluate the discretized residuals of PDEs based on the imposed BCs and NN predictions by following the FEM and using  standard numerical integration schemes. This is achieved through an efficient, discrete, convolution operator-based, and vectorized implementation.
Calculation of the bulk residual is illustrated in Fig. \ref{fig:NN}d with detailed implementations and procedures provided in the SI.

In this work, we define one unique BVP as imposing a certain PDE on a specific domain with specific boundary values at specific boundary conditions.
Changes in any of these elements defines a new boundary value problem.
We found that training a single NN to solve multiple boudnary value prooblems with both Dirichlet and Neumann BCs can be very challenging. 
We introduced a zero-initialization step to address the slower convergence for problems under Neumann BCs (see SI for detailed discussions).
When training BNNs to solve multiple boundary value problems simultaneously, their parameters can stagnate around some local minima, leading to poor performance. 
To address this issue, we use a warm start approach by initializing the mean of a BNN with the optimized parameters from a deterministic NN with the identical architecture.
Our method works for both small and large datasets of boundary value problems. If the number of unique boundary value problems is small, we replicate them to obtain an augmented dataset for training.

\subsection*{Steady-state diffusion problem with small dataset}

\begin{figure}
    \centering
    \includegraphics[width=1.0\linewidth]{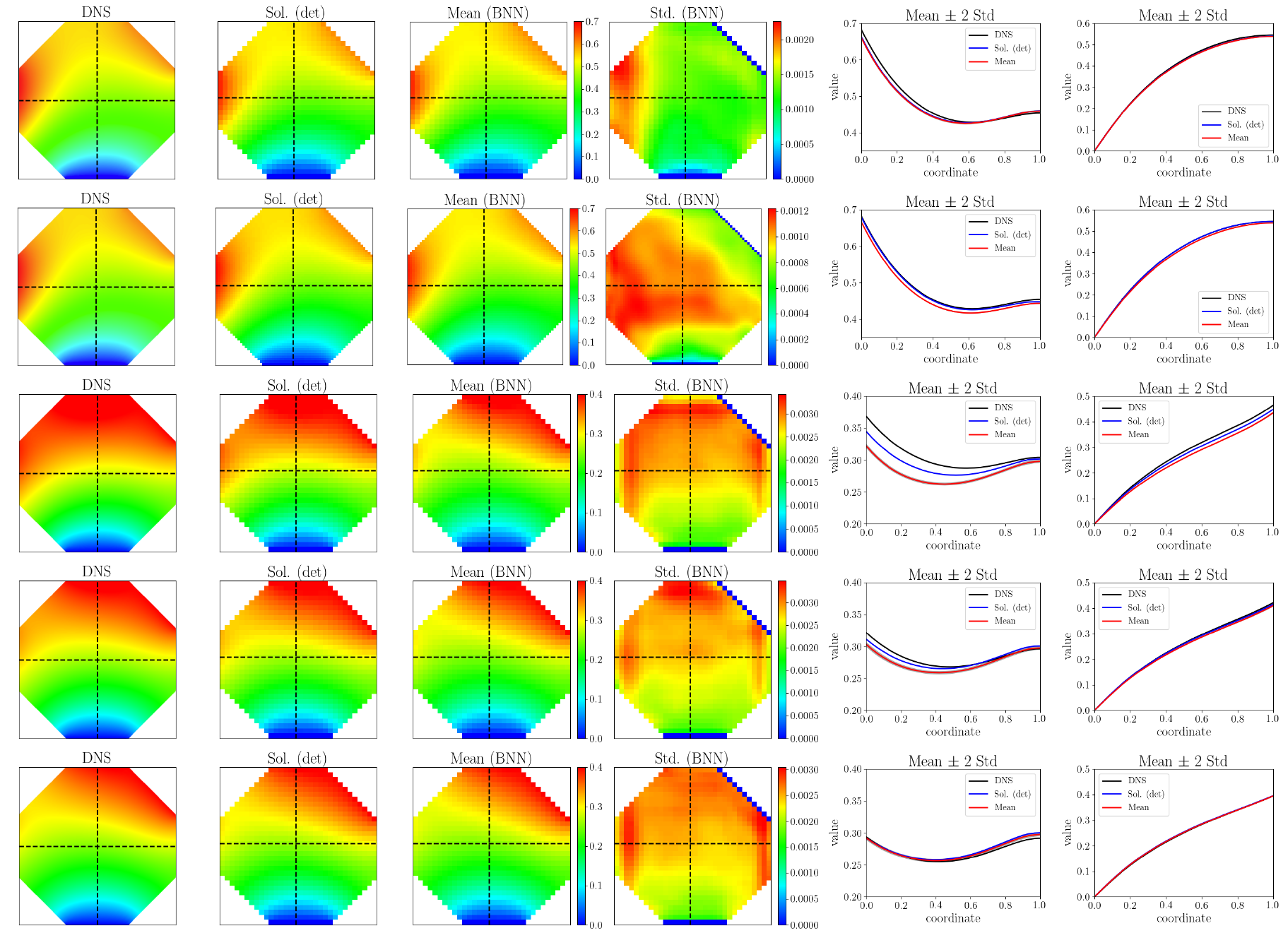}
    \caption{(B)NN-PDE-S for steady-state diffusion. Each row contains the FEM solution (labeled as DNS),  solutions from NN-PDE-S, the mean and std of BNN-PDE-S solutions (warm started from the  NN-PDE-S), and a quantitative comparison between the FEM solution, the NN-PDE-S solution, and the mean $\pm$ 2 std of BNN-PDE-S results along the two dashed lines. Rows 1 and 2: results for the same boundary value problem with mixed BCs at different mesh resolutions with 32 $\times$ 32 for Row 1 and 64 $\times$ 64 for Row 2. The problem setup for Rows 1 and 2 is shown in the SI. Rows 3-5 results for boundary value problems with identical BCs but different material parameters simultaneously solved by a single (B)NN-PDE-S at a mesh resolution of 32 $\times$ 32.  
  }
    \label{fig:diffusion}
\end{figure}

In this first example (Fig. \ref{fig:diffusion}a-b), we use the (B)NN-PDE-S on a single steady-state diffusion boundary value problem on an octagonal domain with imposed mixed BCs (zero/non-zero Dirichlet and non-zero Neumann BCs) at two different mesh resolutions. 
The DNS solution from FEM, NN solution from a deterministic NN, and the mean and std. of BNN results from 50 Monte Carlo samplings for a mesh resolution of 32 $\times$ 32 are shown in Fig. \ref{fig:diffusion}a.
The optimized parameters from the deterministic NN are used for the warm start of the BNNs. 
The deterministic NN results, the mean $\pm$ 2 std of BNN results, and the FEM solution, which is considered the grand truth, are quantitatively commpared along the two dashed lines in Fig. \ref{fig:diffusion}a. These quantitative comparisons confirm the accuracy of the NN results.
The results in Fig. \ref{fig:diffusion}b for a mesh resolution of 64 $\times$ 64 show the same accuracy of the NN results. 
In the second example (Fig. \ref{fig:diffusion}c-e), we use the NN-PDE-S to simultaneously solve three boundary value problems with identical BCs but different material parameters at a mesh resolution of 32 $\times$ 32. The quantitative comparisons in Fig. \ref{fig:diffusion}c-e demonstrate the ability of a single, trained NN to simultaneously solve multiple booundary value prooblems.



\subsection*{Linear elasticity problem with small dataset}

\begin{figure}
    \centering
    \includegraphics[width=1.0\linewidth]{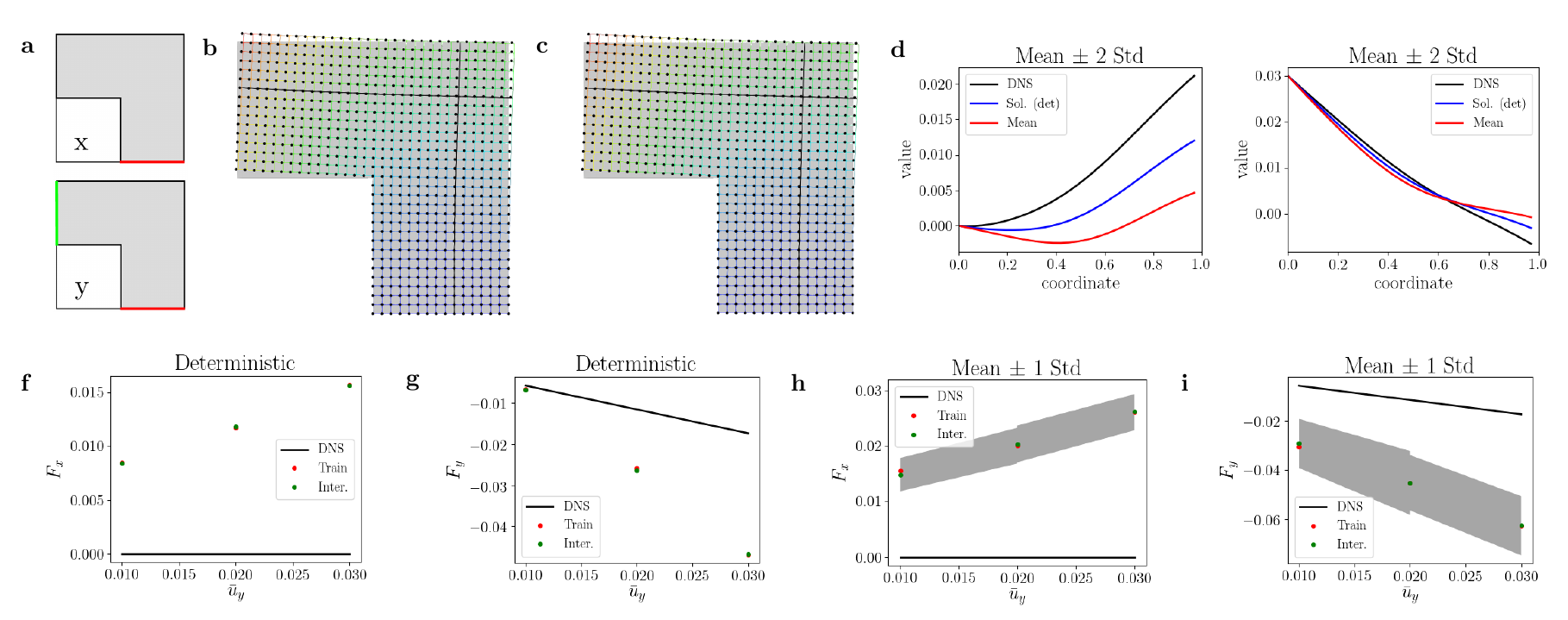}
    \caption{(B)NN-PDE-S for linear elasticity. 
    (a) setup of BCs for the L-shape problem. 
  (b) comparison of solutions between the FEM and the  NN-PDE-S. Both solutions are plotted in the deformed configuration with the FEM solution being illustrated by the mesh grid and the deterministic NN solution being illustrated by the solid black dot. 
  (c) similar solution comparison as in (b), but between the FEM solution and the mean of the BNN-PDE-S solutions among 50 Monte Carlo samplings. 
  (d) comparison of displacements in the X-direction along the vertical line and in the Y-direction along the horizontal line among the FEM solution, the  NN-PDE-S solution, whose parameters are used to warm start the BNN-PDE-S, and the mean and std of solutions from a BNN-PDE-S for a solution resolution of $32\times32$.
  (f, g) comparison of reaction forces in both X- and Y-direction between the FEM and the NN-PDE-S.
  (h, i) comparison of reaction forces in both $X$ and $Y$-direction between the FEM and the BNN-PDE-S.}
    \label{fig:linear-elasticity}
\end{figure}

It is challenging to solve linear elasticity mainly because the governing PDE is written in terms of the infinitesimal strain, which is the gradient of the displacement field and has a very small magnitude $\sim 10^{-4}$. 
Here, we use the (B)NN-PDE-S  for the displacement on an L-shape domain, which is fixed in both directions on the bottom edge and has a vertical displacement applied on the left edge.
A mesh resolution of 32 $\times$ 32 is considered. 
The problem is defined in Fig. \ref{fig:linear-elasticity}a. 
We consider three incremental loading levels and treat each loading level as a different boundary value problem. 
The BNN used to solve these boundary value problems is trained with a warm start.
The deformed geometries from both the FEM results and the deterministic NN results are compared in Fig. \ref{fig:linear-elasticity}b, where the FEM solution is illustrated by the mesh and the NN solution by the solid black dots. The corresponding comparison between the FEM and the mean of the BNN solution is shown in \ref{fig:linear-elasticity}c. 
The quantitative comparisons for $\bu_x$ along the vertical lines and $\bu_y$ along the horizontal lines in Fig. \ref{fig:linear-elasticity}(b,c) between the FEM results, the deterministic NN solution, and the mean $\pm$ 2 std of BNN results are shown in Fig. \ref{fig:linear-elasticity}(d). 
Those results confirm the accuracy of the NN-based solver.
Additionally, we also show the comparison of reaction forces for different loading levels in both the $x-$ and $y-$ directions between the FEM solution and the deterministic NNs (Fig. \ref{fig:linear-elasticity}f,g) and between the FEM solution and the mean $\pm$ 2 std of BNN results (Fig. \ref{fig:linear-elasticity}h,i).
However, because the magnitude of the solution at the bottom of the L-shape is low, the NNs have difficulty computing the total reaction forces.
Another example using a single deterministic NN and BNN to solve 30 linear elastic boundary value prooblems for five different domains with six sets of BCs applied to each appears in the SI.

\subsection*{Nonlinear elasticity problem with small dataset}

\begin{figure}
    \centering
    \includegraphics[width=1.0\linewidth]{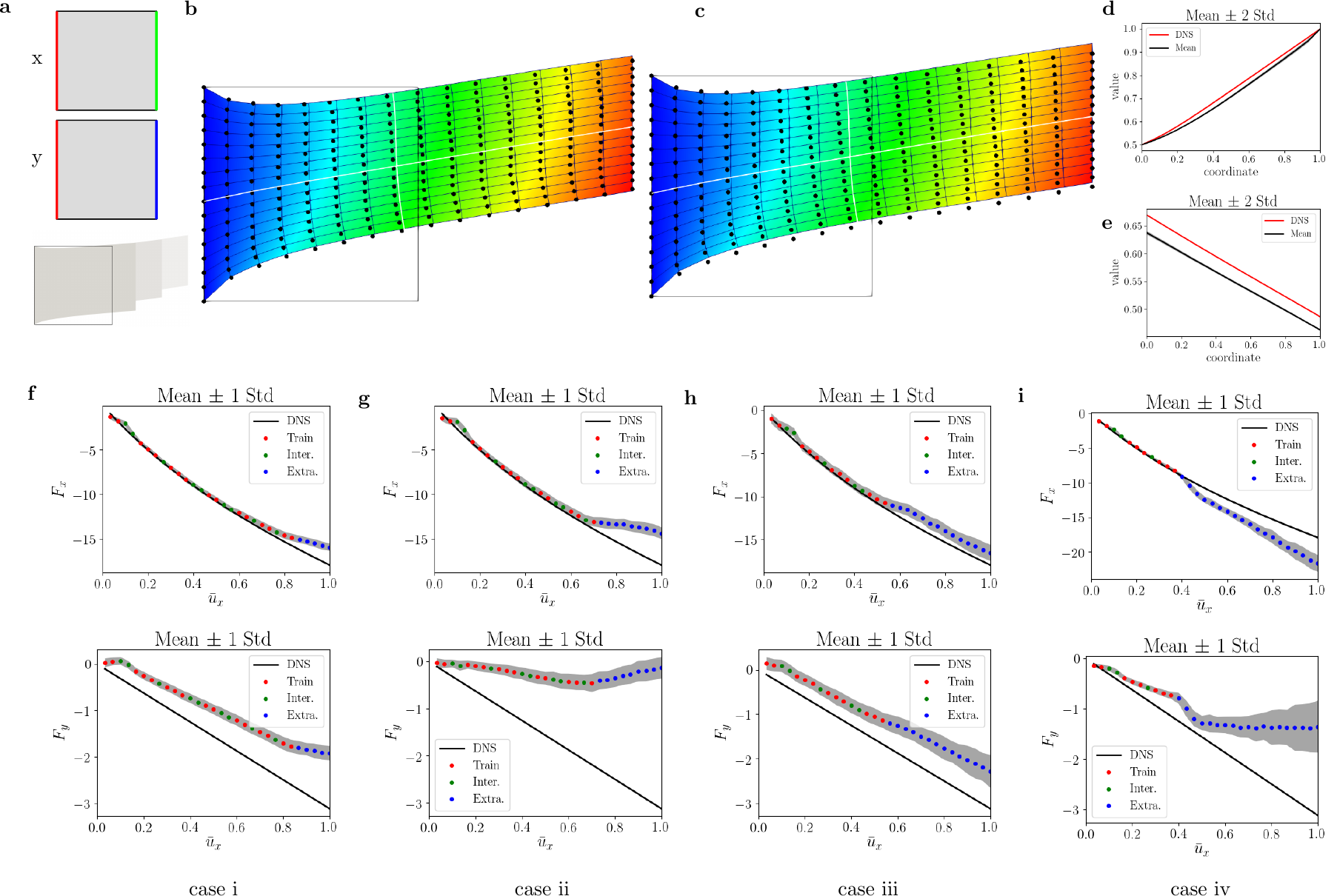}
    \caption{(B)NN-PDE-S for nonlinear elasticity. 
    (a) definition of the boundary value problem for the interpolated and extrapolated loading example. Three loading levels and the deformed shapes are illustrated. 
    (b) Comparison between the FEM solution and the  NN-PDE-S solution with a resolution of $16\times16$ for the last extrapolated BCs for case i. The deformed shapes are plotted with the FEM solution  illustrated by the mesh and the NN solution by the solid black dot. 
    (c) The corresponding comparison between the FEM solution and the mean of the BNN-PDE-S solution over 50 MC samplings. 
    (d) Comparison of displacement in the $X$-direction along the horizontal white lines in (b, c) between the solutions of the FEM, the  NN-PDE-S whose parameters are used to warm start the BNN-PDE-S, and the mean and std of solutions from the BNN-PDE-S.
    (e) Similar comparison as in (d), but for displacement in the $Y$-direction along the vertical white lines in (b, c). 
    (f-i) reaction forces in both $X$- and $Y$-directions for four different cases with each containing a different number of training and testing data.
  }
    \label{fig:nonlinear-elasticity}
\end{figure}

\begin{figure}
    \centering
    \includegraphics[width=1.0\linewidth]{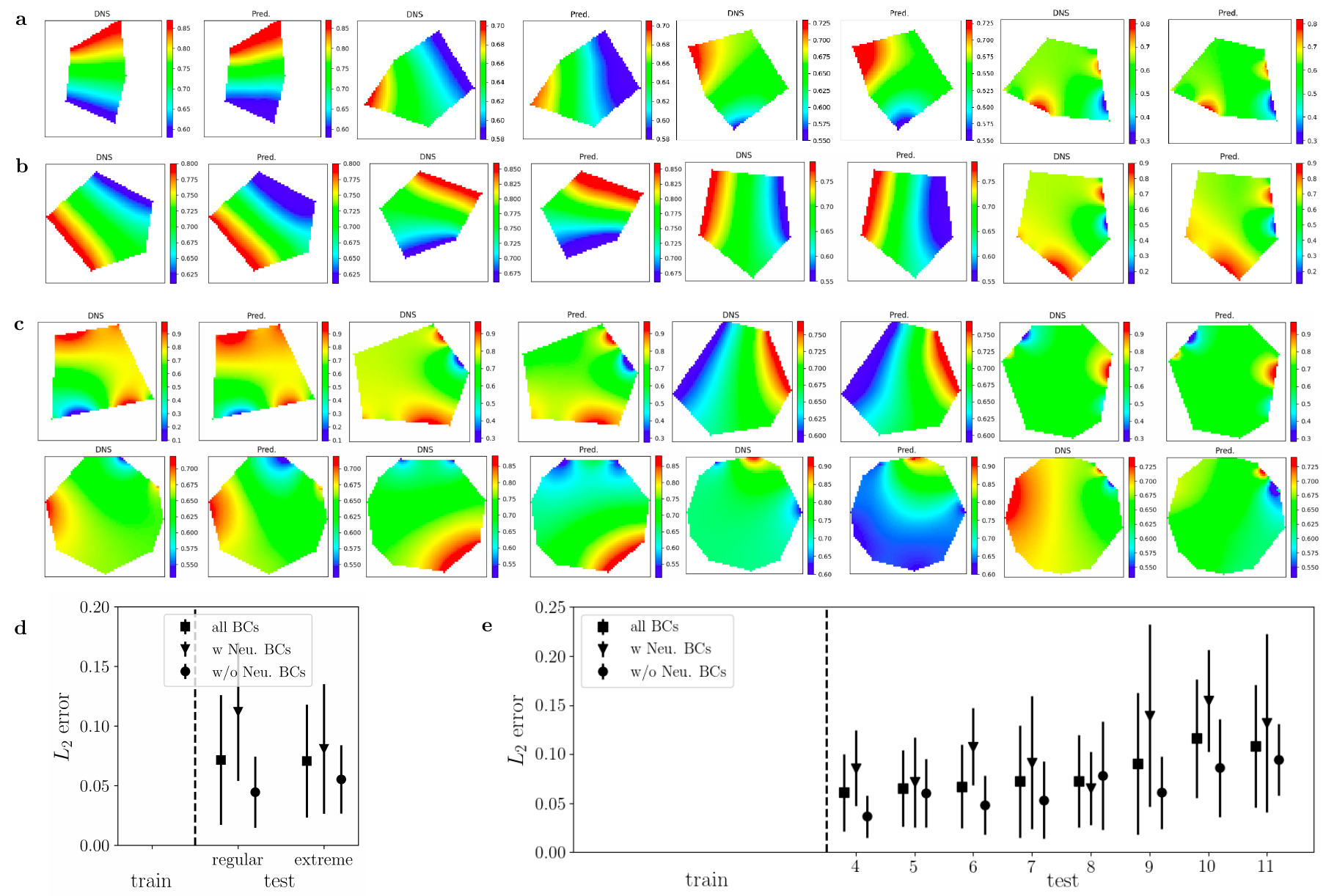}
    \caption{NN-PDE-S for steady-state diffusion with large dataset. 
    (a) Selected good NN-PDE-S predictions for new, randomly generated extreme pentagons. 
    (b) Selected good NN-PDE-S predictions for new, randomly generated regular pentagons. 
    (c) Selected good NN-PDE-S predictions for new, randomly generated polygons with different total numbers of edges. 
    (d) $L_2$ error of NN-PDE-S results for the pentagon example study. $L_2$ error from all boundary value problems, those with Neumann BCs, and others without Neumann BCs are plotted for both training and testing datasets. 
  (e) $L_2$ error of NN-PDE-S results for the polygon example study. $L_2$ error from all boundary value problems, those with Neumann BCs, and others without Neumann BCs are plotted for both training and testing datasets. 
}
    \label{fig:large-diffusion}
\end{figure}

We use the (B)NN-PDE-S solver framework on a nonlinear elasticity problem on a square domain fixed in both directions on the left edge, and with a horizontal displacement loading and vertical traction loading applied on the right edge. 
We considered 30 incremental loading levels and treated each as a different boundary value problem, as shown in Fig. \ref{fig:nonlinear-elasticity}c. 
The BNNs are trained with a warm start for four cases.
The training dataset for each case has different distributions of loading levels. The testing dataset  contains both interpolated and extrapolated loading levels, as indicated by different colors in Fig. \ref{fig:nonlinear-elasticity}f-i. 
The results, including the FEM solution, the deterministic NN, and the mean and std from the BNN, for the last extrapolated loading level for case (i) in both $X$- and $Y$-directions are shown in Fig. \ref{fig:nonlinear-elasticity}b and \ref{fig:nonlinear-elasticity}c.
Additionally, we show the quantitative comparison between FEM results and the mean $\pm$ 2 std of BNN results along the dashed lines (Fig. \ref{fig:nonlinear-elasticity}d and \ref{fig:nonlinear-elasticity}e).
These results confirm the accuracy of the NN-based solver and demonstrate its predictivity for unseen extrapolated BCs.
The reaction forces computed from the NN full field solution are plotted in Fig. \ref{fig:nonlinear-elasticity}f-i, which shows improved accuracy compared to the results for linear elasticity  (Fig. \ref{fig:linear-elasticity}f-i).
For each case, we observe that the results from interpolated BCs are generally accurate. For extrapolated BCs, NN predictions are accurate to a certain degree, particularly for the reaction force in the $X$-direction.
Details of the treatment for computing the nonlinear kinematics, constitutive relations, and discretized residuals with NNs, NN architecture and training scheme, and an additional example on solving 30 boundary value problems with the proposed method are provided in the SI.

\subsection*{Steady-state diffusion problem with large dataset}

Lastly, we use the BNN-PDE-S on the steady-state diffusion problem for two very large datasets with a $64\times64$ resolution. 
In these datasets, the geometries of problem domains and the values of BCs are randomly generated. The locations of BCs are randomly selected from two non-adjacent edges. The boundary values along one edge could be constant, a linear distribution, quadratic distribution, or sinusoidal distribution. The details of the data generation scheme and the statistics of our datasets are discussed in the SI.

In the first case, the training dataset contains 168,000 unique boundary value problems with pentagons, where none of the inner angles exceed $160^\circ$. 
The testing dataset is newly generated with 80 boundary value problems on pentagons with all angles $<160^\circ$ and 80 on ``extreme" pentagons (those with at least one angle $\ge 160^\circ$). Some of the selected results for different geometries are shown in Fig. \ref{fig:large-diffusion}a, which shows that the trained NNs could predict the solution for unseen geometries, including near-degenerate polygons, and unseen BCs. 
The accuracy of the NN results is evaluated by the volume averaged $L_2$ error across all boundary value problems in the corresponding (training/testing) dataset, which is calculated as
    \begin{equation}
      \Vert e\Vert_2 = \frac{1}{N_\mathrm{BVP}} \sum_{l=1}^{N_\mathrm{BVP}} \left(  \sqrt{ \frac{1}{K_\mathrm{px}} \sum_{k=1}^{K_\mathrm{px}} \left(  y_{l,k}^\mathrm{DNS} - y_{l,k}^\mathrm{NN} \right)^2}\right).
    \end{equation}
    The $L_2$ errors of NN results for both training and testing datasets are plotted in Fig. \ref{fig:large-diffusion}(d). We observe that boundary value problems with only Dirichlet BCs generally perform better than those with Neumann BCs. 

In the second case, the training dataset contains 192,000 unique boundary value problems with randomly generated geometries of quadrilaterals, pentagons, and hexagons, whereas the testing dataset is newly generated with the number of edges of spanning from four to eleven with 16 boundary value problems for each type of polygon. Some of the selected results for different polygons appear in Fig. \ref{fig:large-diffusion}c,  again demonstrating that the trained NNs predict  solutions on unseen geometries and unseen BCs. 
The $L_2$ errors of NN results for both training and testing datasets are plotted in Fig. \ref{fig:large-diffusion}(e). We observe some degradation of NN results as the number of polygon edges in the testing set significantly exceed those in the training dataset. However, for training up to hexagons, this degradation in accuracy sets in only for nonagons, decagons and hendecagons.

\begin{table}[t!]
    \centering
    \caption{Summary of the accuracy and speed of the (B)NN-PDE-S for the example presented in Fig. \ref{fig:nonlinear-elasticity}. Note: There may be speedup of the FEM simulation with further code optimization. The wall times of (B)NN-PDE-S results do not include the time to load the solver.}

    \begin{tabular}{c|c|c|c|c}
        \hline
        Solver & Hardware & Software & Wall-time & Averaged $L_2$ error\\ \hline
        FEM & Intel i7-8750, 2.2GHz (use single core) & mechanoChemFEM & 110ms & -\\
        deterministic NN & GeForce GTX 1050 Ti, 4GB memory & Tensorflow & 0.22ms & 2.45e-3 \\ 
        BNN & GeForce GTX 1050 Ti, 4GB memory & Tensorflow & 0.29ms & 3.07e-3 \\ \hline
    \end{tabular}
    \label{tab:table}
\end{table}

\section*{Discussion}

Our framework successfully learns PDE-specific solvers, whether restricted to a single or multiple boundary problems. It has been our focus to develop solvers that generalize across boundary value problems with different domains, boundary conditions and coefficients, and are orders of magnitude faster than the traditional FEM, as shown in Table \ref{tab:table}.
Such features distinguish our approach from certain other NN-based PDE solvers \cite{Raissi2019physics-informed-forward-inverse-jcp,Gao2020PhyGeoNet-PDE-on-Irregular-domain}.

We have used the framework to solve the diffusion problem over two fairly complex, and large datasets of boundary value problems to demonstrate its performance. 
From the corresponding examples, we observe that the (B)NN-PDE-S makes more accurate predictions for interpolated BCs. 
To improve its performance, one can manually introduce new BCs to the training dataset to improve the poorly trained region or targeted prediction region; for instance, expanding the dataset with many geometries if prediction across domains is the goal. If the same geometries are considered for training and testing, the dataset could be augmented by filling the interior domain with random numbers. Naturally, good sampling of BCs and boundary locations is important for learning.
The performance could further be improved via careful hyper-parameter tuning. The learning or cross-validation error computed over this dataset of $\sim 10^5$ boundary value problems with randomly generated polygonal domains (order and shape), boundary conditions and coefficients can be used to drive an active learning workflow that detects regimes (domains, boundary conditions, boundary value functions and coefficients) in which additional boundary value problems are needed for improved learning. 

For BNNs, we applied a constant additive noise to the NN solutions, which is propagated through the PDE residual calculation to compute the training loss function. While the additive noise often is used to account for aleatoric uncertainty in Bayesian inference, here it accounts for model form error of the BNN, and thus corresponds to epistemic uncertainty. The applied noise $\Sigma_1$ functions as the convergence threshold used in the traditional numerical methods for PDEs. 
  The loss of BNNs consists of data-independent (K-L divergence) and  data-dependent (negative log-likelihood) contributions (Methods, Eqs (\ref{eq:bayes}-\ref{eq:bnn-loss-approx-batch})). 
  During training, the data-dependent log-likelihood contribution to the loss is, in general much larger than the data-independent K-L divergence, and decreases. As it does, the data-independent K-L divergence weighs more, leading the mean of the BNN solutions to drift away from the ground truth. With training, both contributions decrease and the mean of the BNN solutions gradually converges to the DNS solution. The K-L divergence term tends to introduce more uncertainty to the model parameters especially with poorly informed priors. The log-likelihood, which depends on $\Sigma_2$ tends to reduce the uncertainty. It controls convergence of the NN solution to the DNS solution as $\Sigma_2$ itself converges.
  We report the evolution of $\Sigma_2$ in Figs S21, S25 and SS29.

  We use a ``warm start'' by initializing the means of the BNN parameters to the optimized parameters from deterministic NNs with an identical architecture. An alternative is to initialize the means of the prior to the optimized parameters from deterministic NNs.

  BVPs with very small variations in the solutions, such as the L-shape linear elasticity example,  pose challenges since the NNs have difficulty capturing small differences. Formally, of course, these problems have low information content. Additionally, the linear elasticity residual is determined by the displacement gradient (infinitesimal strain) field, which is invariant to data normalization. We found that though the residual is computed from the physical, un-normalized, solution, learning is more effective with data normalization. It ensures that the variations of NN outputs (scaled solution) is large $\in [-1,1]$ unlike the infinitesimal strain $\sim 10^{-4}$. Large variations in the NN outputs drive NN parameter variations and favor training. This also applies to diffusion when the solution range is very small. In the residual-based loss, however,  the NN output is scaled back to the physical range, to prevent violation of physics. 

  Hyper-parameter searches are essential for optimal (B)NN-PDE-S. The PDEs differ in their optimal hyperparameters. NNs targeted at solving a wider range of steady state diffusion boundary value problems required wider layers. The vector elasticity problems, even with isotropic properties, have greater information content in their solution field and in general demanded wider layers. The numbers of layers were more closely aligned, and optimal kernel sizes were the same across the PDEs and targeted boundary value problem ranges. See Tables S1, S3, S6, S8, S10, S12.

The NN solvers presented here were trained on a single GeForce GTX 1050 Ti GPU with 4GB memory. Training could take hours for networks designed to solve $\sim 20$ boundary value problems,  and days for those to solve $\sim 10^5$ of boundary value problems. In addition  to more high-powered GPUs, as well as multi-GPU training, optimization of  the training workflow remains unexplored. 
The training time would be further reduced if transfer learning or multi-fidelity learning is used by continued training of previously trained networks. 
Training the (B)NN-PDE-S takes more time than training regular NNs because of the PDE constrained loss layer. The prediction time, however, is  unaffected by these loss layers, as they are not activated during prediction.

In this work, the problem domains are represented via pixels on a square background grid for simplicity. 
Thus, domain boundaries are not smooth curves, but have a pixel-level resolution.
This treatment is of importance as it applies  directly to solving PDEs on pixel-based,  experimental images as domain data--a target future application for our work.
For smooth boundary representations, one can leverage recent work for approaches to map complex and irregular domains onto a regular mesh \cite{Bhatnagar2019CNN-encoder-decoder,Li2020Reaction-diffusion-prediction-CNN,Gao2020PhyGeoNet-PDE-on-Irregular-domain}. 
Such geometric transformations can be taken into account in the proposed PDE loss layers, for instance via the mapping of the physical domain from parent hyper-cubes, as is commonly done in the FEM. While we have considered polygonal domains for their approximation of other geometries, the above mapping could be exploited to remove this restriction, also.

Our approach is formally different from recent operator networks which are focused on learning nonlinear mappings between input function spaces and output spaces, and therefore are mesh independent. In this regard we note that a NN solver that has learnt a PDE on a given discretization (pixel resolution) can serve as the source network for a target finer or coarser mesh within a transfer learning context. The most important difference between the presented (B)NN-PDE-S and DeepONets \cite{lu2021learning,goswami2022deep}, graph kernel networks \cite{li2020neural} and Fourier operator networks \cite{li2020fourier} is that these operator network approaches need labelled field data for training--typically from DNS using the underlying PDE. By not presenting and labelled field solution, but only the domain, BCs and coefficients   to the network, our approach in addition to being label free allows room for our claim that the network is forced to learn the PDE, which by definition holds across boundary value problems. We note that the recent TL-DeepONet \cite{goswami2022deep} uses a source Banach space from a DeepONet trained on labelled data for specific boundary conditions. Further training the final layer of the TL-DeepONet allows transfer learning to some new boundary value problems. We are not aware of the extensiveness to which this transfer learning across boundary value problems has been studied by the authors.
Our approach is very appealing for high-throughput solution of PDEs ranging from inverse and other optimization problems through design and decision-making. Its generalizability across domains and boundary conditions also presents opportunities in topology optimization problems.
Ongoing developments will extend it beyond elliptic PDEs.

\section*{Methods}

\textbf{General elliptic PDEs}
In this work, we develop (B)NN-PDE-S for steady-state diffusion, linear elasticity, and nonlinear elasticity.
These three physical systems are described by a general elliptic PDE on a continuum domain $\Omega \subset \mathbb{R}^n$ with  Dirichlet BCs on $\Gamma^{\Bvarphi}$ and the Neumann BCs on $\Gamma^{\Bk}$:
\begin{equation}
  \begin{aligned}
    \nabla \cdot \BA(\Bvarphi) = \Bzero \quad & \text{on} \quad \Omega, \\
    \Bvarphi (\BX) = \bar{\Bvarphi} (\BX)  \quad & \text{on} \quad  \Gamma^{\Bvarphi}, \\
    \Bk (\BX) = \bar{\Bk} (\BX)  \quad & \text{on} \quad  \Gamma^\Bk,
  \end{aligned}
  \label{eq:general-pde}
\end{equation}
where $\Bvarphi(\BX)$ represents the spatially-dependent unknown field and $\BX \in \mathbb{R}^n$ is the position vector.
The  boundary of the continuum domain satisfies $\overline{\Gamma} = \overline{\Gamma^{\Bvarphi} \bigcup \Gamma^{\Bk}}$ and $\Gamma^{\Bvarphi} \bigcap \Gamma^{\Bk}=\emptyset$.
We use bold typeface for $\Bvarphi$, $\BA$, and $\Bk$ in \eref{eq:general-pde}, depending on the physical system, they could represent either scalar, vector, or tensor fields.
For example, in the diffusion problem, $\Bvarphi$, $\BA$, and $\Bk$ represent the compositional order parameter (scalar), the diffusive flux (vector), and the outflux (scalar), respectively, whereas in elasticity problems, $\Bvarphi$, $\BA$, and $\Bk$ represent the deformation field (vector), the stress field (second-order tensor), and the surface traction (vector), respectively. 
The details of each system appear below.

The weak form of \eref{eq:general-pde} states: For variations $\Bomega$ satisfying $\forall \Bomega \in \scrV$ with $\scrV = \left\{ \Bomega | \Bomega=\Bzero~\text{on}~\Gamma^{\Bvarphi} \right\}$, seek trial solutions $\Bvarphi \in \scrS $ with $ \scrS = \left\{ \Bvarphi|\Bvarphi=\bar{\Bvarphi}~\text{on}~\Gamma^{\Bvarphi} \right\}$ such that
\begin{equation}
  \int_\Omega  \nabla \Bomega \cdot \BA(\Bvarphi) ~ dV - \int_{\Gamma^\Bk} \bar{\Bk} \cdot \Bomega ~dS =0.
  \label{eq:general-pde-weak}
\end{equation}
Eq. \eref{eq:general-pde-weak} is obtained by multiplying (\ref{eq:general-pde}$_1$) with $\Bomega$, integrating by parts, and then incorporating the Neumann BC in (\ref{eq:general-pde}$_3$).
For the diffusion problem, $\Bomega$ is a scalar field. For elasticity problems, $\Bomega$ is a vector field.

Approximate, numerical solutions of \eref{eq:general-pde-weak} can be obtained using its finite-dimensional form. Finite-dimensional approximations of $\Bomega$ and $\Bvarphi$, denoted by $\Bomega^h$ and $\Bvarphi^h$, are constructed with $\Bomega^h \in \scrV^h = \left\{ \Bomega^h | \Bomega^h=\Bzero~\text{on}~\Gamma^{\Bvarphi} \right\}$ and $\Bvarphi^h \in \scrS^h = \left\{ \Bvarphi^h|\Bvarphi^h=\bar{\Bvarphi}~\text{on}~\Gamma^{\Bvarphi} \right\}$. 
The finite-dimensional fields $\Bomega^h$, $\nabla\Bomega^h$, and $\Bvarphi^h$ are computed as
\begin{equation}
  \Bomega^h = \BN \Bd_{\Bomega},
  \quad
  \nabla\Bomega^h = \BB \Bd_{\Bomega},
  \quad
  \text{and}
  \quad
  \Bvarphi^h = \BN \Bd_{\Bvarphi} 
  \label{eq:discretized-form-u}
\end{equation}
in terms of the nodal solutions $\Bd_{\Bomega}$ and $\Bd_{\Bvarphi}$, the basis functions $\BN$, and the basis function gradient operator $\BB = \nabla\BN$.
Inserting \eref{eq:discretized-form-u} into \eref{eq:general-pde-weak} we obtain the discretized residual as an assembly over subdomains $\Omega^e$ and their associated boundary $\Gamma^{e}$ as
\begin{equation}
  \BR = \foo_{e=1}^{n_\text{elem}} \left\{ \int_{\Omega^e} \BB^T \BA(\Bvarphi^h) dV - \int_{\Gamma^{e,\Bk}} \BN^T \bar{\Bk}~dS \right\}
  \label{eq:discretized-residual}
\end{equation}
where $\foo$ is the assembly operator and $n_\text{elem}$ represents the total number of subdomains.
The volume and surface integrations in \eref{eq:discretized-residual} are evaluated numerically via Gaussian quadrature rules.
In this work, the problem domain $\Omega$ is treated as an image. 
Single component, connected graphs whose vertices are pixels form the subdomains $\Omega^e$. Pixel connectivity is preserved as graph edges in the image data.
Field values at each pixel of the image are treated as nodal values. 
Additional discussion on constructing the subdomains based on image pixels is provided in the SI. Of interest, but tangential, in this context is a recent work in which NN layers map between FE meshes of different resolutions \cite{uriarte2022finite}.

\textbf{Steady-state diffusion}
This problem is described by a linear elliptic PDE in the scalar composition field following \eref{eq:general-pde}
\begin{equation}
  \begin{aligned}
    \nabla \cdot \BH = \Bzero \quad & \text{on} \quad \Omega, \\
    C (\BX) = \bar{C} (\BX)  \quad & \text{on} \quad  \Gamma^{C}, \\
    H = \bar{H} (\BX)  \quad & \text{on} \quad  \Gamma^H.
  \end{aligned}
  \label{eq:general-pde-diffusion}
\end{equation}
 
In \eref{eq:general-pde-diffusion}, $C$ represents the composition, $\BH$ is the diffusive flux defined as
\begin{equation}
  \BH = - D \nabla C,
  \label{eq:diffusion-flux}
\end{equation}
with $D$ as the diffusivity, and $H$ is the outward surface flux in the normal direction.
The discretized residual function \eref{eq:discretized-residual} for steady-state diffusion is written as
\begin{equation}
    \begin{aligned}
  \BR  & = \foo_{e=1}^{n_\text{elem}} \left\{ \int_{\Omega^e} \BB^T \BH dV - \int_{\Gamma^{e,H}} \BN^T \bar{H}~dS \right\}. \\
    \end{aligned}
  \label{eq:discretized-residual-diffusion}
\end{equation}
Diffusivity $D \in [1.0, 100.0]$ has been used.

\textbf{Linear elasticity}
This problem also is posed as a linear elliptic PDE, but in terms of a vector field, $\Bu\in \mathbb{R}^n$. Following \eref{eq:general-pde} we have:
\begin{equation}
  \begin{aligned}
    \nabla \cdot \Bsigma = \Bzero \quad & \text{on} \quad \Omega, \\
    \Bu (\BX) = \bar{\Bu} (\BX)  \quad & \text{on} \quad  \Gamma^{\Bu}, \\
    \BT = \bar{\BT} (\BX)  \quad & \text{on} \quad  \Gamma^\BT.
  \end{aligned}
  \label{eq:general-pde-linear-elasticity}
\end{equation}

In \eref{eq:general-pde-linear-elasticity}, $\Bu$ represents the displacement field, $\Bsigma$ is the stress tensor, and $\BT$ is the surface traction.
Here, $\Bsigma$ is related to the infinitesimal strain $\Bvarepsilon=\frac{1}{2}\left( \nabla \Bu + (\nabla \Bu)^T\right)$ via the following constitutive relationship
\begin{equation}
  \Bsigma = \lambda\tr (\Bvarepsilon) \Bone + 2\mu\Bvarepsilon
  \label{eq:constitutive-linear}
\end{equation}
where $\lambda$ and $\mu$ are the Lam\'e constants, and $\Bone$ is the second-order identity tensor.
The discretized residual function \eref{eq:discretized-residual} for the linear elasticity problem is written as
\begin{equation}
  \BR = \foo_{e=1}^{n_\text{elem}} \left\{ \int_{\Omega^e} \BB^T \Bsigma dV - \int_{\Gamma^{e,T}} \BN^T \bar{\BT}~dS \right\}.
  \label{eq:discretized-residual-linear}
\end{equation}
We used $\lambda=14.4231$ and $\mu=9.61538$ in both DNSs with FEM for comparison with the (B)NN-PDE-S and in the PDE loss layers.

\textbf{Nonlinear elasticity}
With the displacement $\Bu \in \mathbb{R}^n$ as the vector field unknown, we write following \eref{eq:general-pde}:
\begin{equation}
  \begin{aligned}
    \nabla \cdot \BP = \Bzero \quad & \text{on} \quad \Omega_0, \\
    \Bu (\BX) = \bar{\Bu} (\BX)  \quad & \text{on} \quad  \Gamma^{\Bu}_0, \\
    \BT = \bar{\BT} (\BX)  \quad & \text{on} \quad  \Gamma^\BT_0,
  \end{aligned}
  \label{eq:general-pde-nonlinear-elasticity}
\end{equation}
for the nonlinear elasticity problem, with the subscript $0$ indicating the reference configuration. 
In \eref{eq:general-pde-nonlinear-elasticity}, $\BP$ is the first Piola-Kirchhoff stress tensor, and $\BT$ is the surface traction.
In the nonlinear elasticity problem, the deformation gradient is defined as $\BF=\Bone + \partial\Bu/\partial\BX$ with $\Bone$ being the second-order identity tensor.
The right Cauchy-Green deformation tensor is written as $\BC=\BF^{T}\BF$.
The following compressible Neo-hookean hyperelastic free energy function is considered 
\begin{equation}
  W=\frac{1}{2} \mu(\text{tr}(\BC)3-3-2 \ln(J))+\lambda \frac{1}{2}(J-1)^2,
\end{equation}
with $\mu$ and $\lambda$ as the Lam\'e constants and $J=\det(F)$.
The Piola stress tensor $\BP$ is computed as
\begin{equation}
  \BP = \frac{\partial W}{\partial \BF} 
  =  \lambda(J^2-J)\BF^{-T} +\mu(\BF-\BF^{-T}).
\end{equation}
The discretized residual function \eref{eq:discretized-residual} for the nonlinear elasticity problem is written as
\begin{equation}
  \BR = \foo_{e=1}^{n_\text{elem}} \left\{ \int_{\Omega^e} \BB^T \BP dV - \int_{\Gamma^{e,T}} \BN^T \bar{\BT}~dS \right\}.
  \label{eq:discretized-residual-nonlinear}
\end{equation}
The same Lam\'{e} constants were used as in linear elasticity.

\textbf{Deterministic loss}
When using mini-batch optimization to train the NN-PDE-S over a dataset $\calD$, where each data point is a boundary value problem with information on problem domain and BCs,
the batch loss $\calL_i$ is written in terms of the reduced total residual $\BR_\text{tot}^\text{red}$, as illustrated in Fig. \ref{fig:NN}, as
\begin{equation}
  \calL_i = \frac{1}{N} \sum_{n=1}^{N} \left(\BR_\text{tot}^\text{red}( \calD_i, \BTheta) \right)^2,
  \label{eq:loss-deterministic}
\end{equation}
for each mini-batch $i=1,2,\cdots,M$ with $N$ indicating the number of data points (boundary value problems) in each mini-batch.
The detailed training of NN-PDE-S is discussed in the SI.

\textbf{Probabilistic loss}
In BNN-PDE-S, each model parameter is sampled from a posterior distribution. 
We solve for the posterior distribution of model parameters with variational inference instead of Markov Chain Monte Carlo (MCMC) sampling, as the latter involves expensive iterative inference steps and is not suitable for systems with a large number of parameters \cite{Blei2017Variational-inference-review,Kingma2014Welling+Variational-bayes}.
In our work, the likelihood function is constructed based on the discretized PDE residuals.
An additive noise is often applied to the NN predicted solution to represent the aleatoric uncertainty \cite{Luo2020Bayesian-deep-learning,Wang2020Garikipati-Xun-perspective-system-id-bayesian,Geneva2020Zabaras-JCP-auto-regressive-NN-PDE,Zhu2018Zabaras-UQ-Bayesian-ED-CNN,Zhu2019Perdikaris-Physics-PDE-CNN}.
Here, we also applied an additive noise to the solution to represent epistemic uncertainty stemming from model form error between BNN-PDE-S, whose perturbed solutions are still constrained by the PDEs, and the FEM solver that yields the DNS solution.

The BNN model parameters $\BTheta$ are stochastic and  sampled from a posterior distribution $P(\BTheta | \calD)$ instead of being represented by single values as in deterministic NNs. 
The posterior distribution $P(\BTheta | \calD)$ is computed based on Bayes' theorem
\begin{equation}
  P(\BTheta | \calD) = \frac{P(\calD|\BTheta) P(\BTheta)}{P(\calD)},
  \label{eq:bayes}
\end{equation}
where $\calD$ denotes the i.i.d.~observations (training data) and $P$ represents the probability density function.
In \eref{eq:bayes}, $P(\calD|\BTheta)$ is the likelihood, $P(\BTheta)$ is the prior probability, and $P(\calD)$ is the evidence, respectively.
The likelihood is the probability of $\calD$ given $\BTheta$, which describes the probability of the observed data for given parameters $\BTheta$.
A larger value of $P(\calD | \BTheta )$ implies that $\BTheta$ is more likely to yield $\calD$.
The prior must be specified to begin Bayesian inference  \cite{gelman2013bayesian}.

To compute the posterior distributions of $\BTheta$, one can use  popular sampling-based methods, such as MCMC. However, MCMC involves expensive iterative inference steps and would be difficult to use when datasets are large or models are very complex \cite{Kingma2014Welling+Variational-bayes,Liu2016Wang-Stein-variational-gradient-descent,Blei2017Variational-inference-review}.
An alternative is variational inference, which approximates the exact posterior distribution $P(\BTheta | \calD)$ with a more tractable surrogate distribution $Q(\BTheta)$ by minimizing the Kullback-Leibler (KL) divergence \cite{Liu2016Wang-Stein-variational-gradient-descent,Blei2017Variational-inference-review,Graves2011Varational-inference}
\begin{equation}
  Q^* = \text{arg~min}~D_\text{KL}(Q(\BTheta)||P(\BTheta | \calD)).
  \label{eq:kl-divergence}
\end{equation}
Variational inference is faster than MCMC and easier to scale to large datasets. 
We therefore explore it in this work, even though it is less rigorously studied than MCMC \cite{Blei2017Variational-inference-review}. 
The KL divergence is computed as 
\begin{equation}
  D_\text{KL}(Q(\BTheta)||P(\BTheta | \calD)) = \mathbb{E}_Q[\log Q(\BTheta)] - \mathbb{E}_Q[\log P(\BTheta, \calD)]  + \log P(\calD),
  \label{eq:kl-divergence-calculation}
\end{equation}
which requires computing the logarithm of the evidence, $\text{log}P(\calD)$ in \eref{eq:bayes} \cite{Blei2017Variational-inference-review}. 
However, computation of $P(\calD)$ would require marginalization over all realizations of $\Theta$--an intractable task. It is also difficult to estimate $P(\calD)$. Consequently, it is challenging to directly evaluate  the objective function in \eref{eq:kl-divergence}.
Alternatively, we can optimize the so-called evidence lower bound (ELBO) which is equivalent to the KL-divergence up to an added constant. Therefore, an optimal distribution determined using the ELBO is also optimal for the KL-divergence.
\begin{equation}
  \begin{aligned}
    \text{ELBO}(Q) &=  -D_\text{KL}(Q(\BTheta)||P(\BTheta | \calD)) + \log P(\calD)  \\
  \text{ELBO}(Q) &= \mathbb{E}_Q[\log P(\BTheta, \calD)] - \mathbb{E}_Q[\log Q(\BTheta)] \\
  & = \mathbb{E}_{Q}[\log P(\calD | \BTheta )] - \left( \mathbb{E}_{Q}[\log Q(\BTheta)] - \mathbb{E}_{Q}[\log P(\BTheta)] \right) \\
   & = \mathbb{E}_{Q}[\log P(\calD | \BTheta )] - D_\text{KL}\left( Q(\BTheta) ||P(\BTheta)  \right).
  \end{aligned}
  \label{eq:elbo}
\end{equation}

So, the loss function for the BNN is written as
\begin{equation}
  \calL = D_\text{KL}\left( Q(\BTheta) ||P(\BTheta)  \right) - \mathbb{E}_{Q}[\log P(\calD | \BTheta )],
  \label{eq:bnn-loss-exact}
\end{equation}
which consists of a prior-dependent but data-independent part and a data-dependent part.
The former is the KL-divergence of the surrogate posterior distribution $Q(\BTheta)$ and the prior $P(\BTheta)$, and the latter is the negative log-likelihood.
For mini-batch optimization, the batch loss is written as 
\begin{equation}
  \calL_i = \frac{1}{M} D_\text{KL}\left( Q(\BTheta) ||P(\BTheta)  \right) - \mathbb{E}_{Q}[\log P(\calD_i | \BTheta^{(i)} )],
  \label{eq:bnn-loss-batch}
\end{equation}
for each mini-batch $i=1,2,\cdots,M$ \cite{Blundell2015Weight-Uncertainty-NN}.
With \eref{eq:bnn-loss-batch}, we have $\calL = \sum_i \calL_i$.
Following Ref \cite{Blundell2015Weight-Uncertainty-NN}, Monte Carlo (MC) sampling is used to approximate the expectation in \eref{eq:bnn-loss-batch} as 
\begin{equation}
  \calL_i \approx \frac{1}{M} D_\text{KL}\left( Q(\BTheta) ||P(\BTheta)  \right) - \frac{1}{N} \sum_{n=1}^{N} \log P(\calD_i^n | \BTheta^{(i)} ),
  \label{eq:bnn-loss-approx-batch}
\end{equation}
where $N$ is the size of each mini-batch dataset, and $\BTheta^{(i)}$ denotes the $i$th batch sample drawn from the posterior distribution $Q(\BTheta)$. 
Even though only one set of parameters $\BTheta^{(i)}$ is drawn from $Q(\BTheta)$ for each mini-batch, the perturbation approach proposed by Flipout (see SI) ensures that parameters are de-correlated for each individual example $\calD_i^n$ in calculating the log-likelihood.
Probabilistic dense layers and convolutional layers with the Flipout weight perturbation technique have been implemented in the TFP Library \footnote{\href{https://www.tensorflow.org/probability/api_docs/python/tfp/layers}{www.tensorflow.org/probability/api\_docs/python/tfp/layers}} and are used to construct the BNNs in this work.


\textbf{Neural network structure and loss function}
Using modularized implementation of probabilistic layers in the TFP library it is easy to construct the BNN-PDE-S to have the encoder-decoder architecture shown in Fig. \ref{fig:NN}, which is similar to the NN-PDE-S but with all weights being drawn from probability distributions. 
The loss of the BNNs is given in \eref{eq:bnn-loss-exact}. 
The probabilistic layers in the TFP library automatically calculate the prior-dependent KL-divergence and add it to the total loss.

The data-dependent loss is accounted for by the likelihood. Assuming Gaussian noise $\Bepsilon \sim \calN(\Bzero, \Sigma_1\BI)$ with a zero-mean and a pre-specified constant covariance $\Sigma_1$, the NN representation $\Bf(\Bx, \BTheta)$ is augmented by noise to yield the output $\By$:
\begin{equation}
    \By = \Bf(\Bx, \BTheta) + \Bepsilon.
    \label{eq:f-surrogate}
\end{equation}

Forward solutions are obtained by seeking to drive the residual to zero, and correspond to satisfaction of the weak form. For BNN-PDE-S, the likelihood function is constructed from the residual value, rather than the NN predicted solutions, thus ensuring that the framework remains label free.
With the noise $\Bepsilon$ in \eref{eq:f-surrogate} propagating through the residual calculation, the likelihood function is written as
\begin{equation}
  P(\BR_\text{tot}^\text{red}(\calD_i,\BTheta^{(i)})|\Bzero, \Sigma\BI) = \prod_{k=1}^{K} \calN \left( R_\text{tot}^\text{red,$k$}|0, \Sigma_2 \right)
  \label{eq:likelihood-function}
\end{equation}
where index $k$ indicates the pixel number with $K$ total pixels and $\BR_\text{tot}^{\text{red},k}$ is the component of $\BR_\text{tot}^\text{red}$ at pixel $k$. For systems where nonlinear operations are involved in the residual calculation, the residual noise distribution is in general non Gaussian even if the noise in the BNN outputs is assumed to be Gaussian. Under the conditions that $\Sigma_1$ is small and the nonlinear operations are smooth, there exists a neighborhood of any point in which the linear approximation is valid (higher-order terms being negligible), we assume that $\Sigma_2$, the noise distribution of the residual, is approximately Gaussian.
As it is challenging to directly calculate $\Sigma_2$ via error propagation based on $\Sigma_1$, we treat $\Sigma_2$ as a learnable parameter to be optimized based on the NN loss.
In \eref{eq:likelihood-function}, $\Sigma_2$ essentially serves as a convergence threshold for the residual. 
The batch-wise loss of the residual constrained BNNs has the following form
\begin{equation}
  \calL_i \approx \frac{1}{M} D_\text{KL}\left( Q(\BTheta) ||P(\BTheta)  \right) - \frac{1}{N} \sum_{n=1}^{N} \sum_{k=1}^{K} \log \left( \calN \left( R_\text{tot}^\text{red,$k$}(\calD_i^n,\BTheta^{(i)}) | 0, \Sigma_2 \right) \right).
  \label{eq:loss-probabilistic}
\end{equation}
The detailed training scheme for BNN-PDE-S is discussed in the SI.

\textbf{Uncertainty quantification}
The BNNs allow us to quantify the epistemic uncertainty from model parameters.
With the discretized residual constrained BNNs, the posterior predictive distribution $P(\By^*|\Bx^*, \calD)$ of the BNN-predicted full field solution $\By^*$ for a specific testing data point $\Bx^*$, i.e. one boundary value problem with information on problem domain and BCs is \cite{Zhu2018Zabaras-UQ-Bayesian-ED-CNN, Luo2020Bayesian-deep-learning} 
\begin{equation}
  \begin{aligned}
  P(\By^*|\Bx^*, \calD) & = \int P(\By^*|\Bx^*, \BTheta) P(\BTheta | \calD) d \BTheta \\
  & \approx \int P(\By^*|\Bx^*, \BTheta) Q(\BTheta) d \BTheta, \\
  \end{aligned}
\end{equation}
which can be numerically evaluated via MC sampling as
\begin{equation}
  \begin{aligned}
    P(\By^*|\Bx^*, \calD) & \approx \frac{1}{S} \sum_{s=1}^{S}  P(\By^*|\Bx^*, \BTheta^s) 
    \quad \text{where} \quad \BTheta^s \sim Q(\BTheta), \\
  \end{aligned}
\end{equation}
with $s$ indicating each sampling.
To represent the uncertainty, we compute the statistical moments of $\By^*$ via the predictive expectation 
\begin{equation}
  \mathbb{E}[\By^*|\Bx^*, \calD] \approx \frac{1}{S} \sum_{s=1}^{S} \Bf(\Bx^*, \BTheta^s)
\end{equation}
and the predictive variance 
\begin{equation}
  \begin{aligned}
    \text{Var}[\By^*|\Bx^*, \calD] 
    & = \mathbb{E}[(\By^* + \Bepsilon)^2] - (\mathbb{E}[\By^* + \Bepsilon])^2 \\
  & \approx \frac{1}{S} \sum_{s=1}^{S} \left(  \Bf(\Bx^*, \BTheta^s)\Bf^T(\Bx^*, \BTheta^s)  + \Sigma_1\BI \right)  
  - \left( \frac{1}{S} \sum_{s=1}^{S} \Bf(\Bx^*, \BTheta^s)  \right)\left(\frac{1}{S} \sum_{s=1}^{S} \Bf(\Bx^*, \BTheta^s)  \right)^T. \\
  \end{aligned}
\end{equation}


\section*{Code Availability}
Our modularized code implementation is publicly available\footnote{\href{https://github.com/mechanoChem/mechanoChemML}{github.com/mechanoChem/mechanoChemML }}, which will assist the extension to other PDE systems. 

\section*{Acknowledgements}
We gratefully acknowledge the support of Toyota Research Institute, Award \#849910: ``Computational framework for data-driven, predictive, multi-scale and multi-physics modeling of battery materials''.  
Computing resources were provided in part by the National Science Foundation, United States via grant 1531752 MRI: Acquisition of Conflux, A Novel Platform for Data-Driven Computational Physics (Tech. Monitor: Ed Walker). 
This work also used the Extreme Science and Engineering Discovery Environment (XSEDE) Comet at the San Diego Supercomputer Center and Stampede2 at The University of Texas at Austin's Texas Advanced Computing Center through allocation TG-MSS160003 and TG-DMR180072.

\section*{Author Contributions}

\section*{Competing Interests statement}

\bibliographystyle{unsrt} 

\bibliography{lib.bib}
\end{document}


\maketitle

\section{Training NNs with stochastic weights - Flipout}

Different methods are available for training neural networks (NNs) with stochastic weights, such as weight perturbation \cite{Graves2011Varational-inference,Blundell2015Weight-Uncertainty-NN,Wen2018GrosseFlipout}, activation perturbation \cite{Ioffe2015Szegedy-batch-normalization}, reparameterization \cite{Kingma2014Welling+Variational-bayes}, and many others.
In this work, we follow a specific weight perturbation method, the so-called Flipout, proposed in Ref \cite{Wen2018GrosseFlipout}. 
Compared with other weight perturbation algorithms that suffer from high variance of the gradient estimates because the same perturbation is shared in a mini-batch for all training examples, Flipout is an efficient method, which decorrelates the gradients in a mini-batch by implicitly sampling pseudo-independent weight perturbation for each example, and thus reduces the variance of NNs with stochastic weights \cite{Wen2018GrosseFlipout}.
This method can be efficiently implemented in a vectorized manner with unbiased stochastic gradients.

A brief description of Flipout is summarized here. Readers are directed to Ref. \cite{Wen2018GrosseFlipout} for details.
Flipout assumes that the perturbations of different weights are independent, and the perturbation distribution is symmetric around zero.
Under such assumptions, the perturbation distribution is invariant to element-wise multiplication by a random sign matrix.
To minimize the loss $\calL$, the distribution of $Q(\BTheta)$ can be described in terms of perturbations with $W=\lbar{W}+\Delta W$, where $\lbar{W}$ and $\Delta W$ are the mean and a stochastic perturbation for NN parameters $\BTheta$, respectively.
Flipout uses a base perturbation $\widehat{\Delta W}$ shared by all samples (training data points) in a mini-batch, and arrives at the perturbation for an individual sample by multiplying $\widehat{\Delta W}$ with a different rank-one sign matrix
\begin{equation}
  \Delta W_n =  \widehat{\Delta W} \circ r_n s_n^t,
\end{equation}
where the subscript $n$ indicates an individual example in a mini-batch, the superscript $t$ denotes the transpose operation, and $r_n$ and $s_n$ are entries of random vectors uniformly sampled from $\pm 1$.
Using different perturbations for each sample in a mini-batch rather than an identical perturbation for all the example in a mini-batch ensures the reduction of the variance of the stochastic gradients in Flipout during training.
For BNNs, the $\lbar{W}$ and $\widehat{\Delta W}$ are the mean and standard deviation of the posterior distribution $Q(\BTheta)$, which are obtained via backpropagation with stochastic optimization algorithms.

\section{Efficient implementation of the residual calculation} \label{sec:NN-FEM} 

\begin{figure}[h!]
    \centering
    \includegraphics[width=0.95\linewidth]{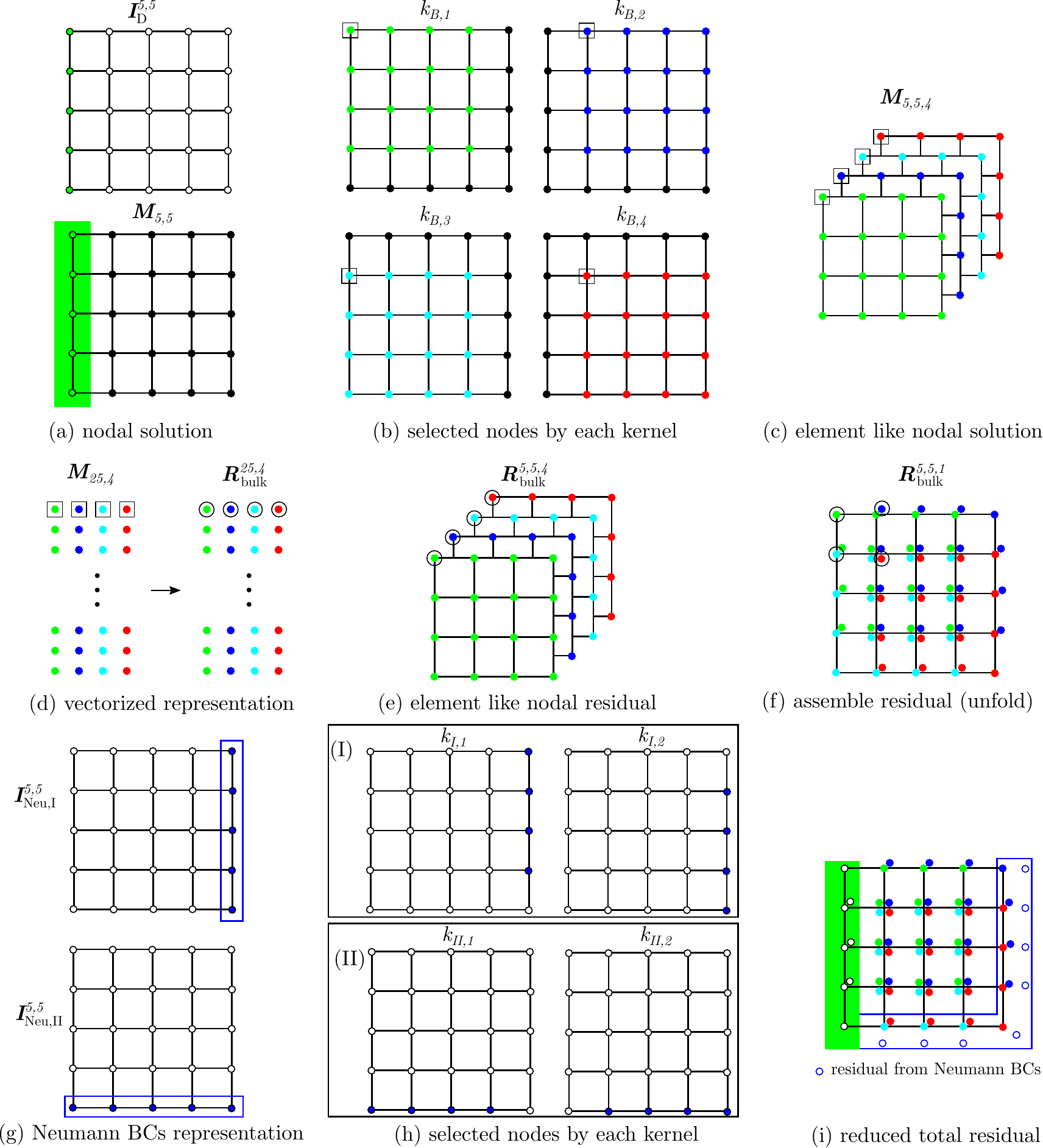}
    \caption{Illustration of the implementation steps of computing the residual in the weak PDE loss layers. Paddings of zeros are not shown in (c,d,e). }
    \label{fig:schematic-weak}
\end{figure}

In this section, we describe the implementation details of the weak PDE loss layers.
We heavily utilize the convolutional operation, and the vector/matrix/tensor operations to achieve numerical efficiency.
Readers are directed to our source code for additional details\footnote{\href{https://github.com/mechanoChem/mechanoChemML}{github.com/mechanoChem/mechanoChemML}}. 
As shown in Fig. \ref{fig:NN}b, the weak PDE loss layers take both NN inputs (BCs information) and outputs (NN predicted solution) as their inputs.
The data structure to represent the BCs is discussed in detail in Section \ref{sec:data-representation}.
Schematics of the major implementation steps of the bulk residual calculation and the residual calculation of Neumann BCs in the weak PDE loss layers are shown in Fig. \ref{fig:schematic-weak}, respectively. 

We choose a steady state diffusion problem with a scalar unknown at each node for the purpose of illustration, with Dirichlet BCs being applied on the left boundary, non-zero Neumann BCs being applied on the bottom and right boundaries, and zero Neumann BCs on the top.
To fix ideas, we consider an example such that the output of the NN shown in Fig. \ref{fig:schematic-weak} is a $5\times5$ matrix (an image with $5\times5$ pixels), denoted as $\BM_{5,5}^\text{NN}$ with the value of each entry being the actual concentration, $\BM_{5,5}^\text{NN}$ is equivalent to the nodal solution on a domain, which is discretized by a $4\times 4$ block of elements, as shown in Fig. \ref{fig:schematic-weak}(a).\footnote{In the source code, $\BM_{5,5}^\text{NN}$ is stored as $\BM_{5,5,1}^\text{NN}$ with a third dimension of 1, which indicates the DOF per node. For elasticity problems, the third dimension has a size of 2. Here, we drop the ``1'' to simplify the notations.}
The implementation procedure is summarized in the Algorithm Box \ref{algo:compute-residual}.

\subsection{Dirichlet BCs}
The channel of NN inputs with Dirichlet BCs information is denoted as $\BI_\text{D}^{5,5}$, as shown in Fig. \ref{fig:schematic-weak}(a). 
To enforce the Dirichlet BCs, we replace the nodal values of $\BM_{5,5}^\text{NN}$ at the location of Dirichlet boundary with the actual values of $\BI_\text{D}^{5,5}$ to obtain a new matrix, denoted as $\BM_{5,5}$, as indicated by the green color in Fig. \ref{fig:schematic-weak}(a).
The Dirichlet BCs are then automatically incorporated into the residual vector during the bulk residual calculation discussed in the next Section.

\subsection{Bulk residual}

The matrix representation of the nodal solution automatically contains the element connectivity information of the mesh.
To compute the bulk residual, we first apply convolutional operations to $\BM_{5,5}$ with the following kernels
\begin{equation}
  k_{B,1} =
\begin{bmatrix}
1 & 0\\
0 & 0\\
\end{bmatrix}, 
\quad
k_{B,2} =
\begin{bmatrix}
0 & 1\\
0 & 0\\
\end{bmatrix}, 
\quad
k_{B,3} =
\begin{bmatrix}
0 & 0\\
1 & 0\\
\end{bmatrix}, 
\quad
k_{B,4} =
\begin{bmatrix}
0 & 0\\
0 & 1\\
\end{bmatrix}.
  \label{eq:kernel-node-to-bulk-elem}
\end{equation}
Each convolutional operation results in a matrix with a size of $5\times5$,\footnote{The resulting matrix size is $4\times4$. Zero paddings are used to ensure the resulted matrix with a dimension of $5\times5$. Keeping the matrix size unchanged during the convolutional operations is not necessary and might require a small amount of extra floating-point operations, but it is less prone to errors if we handle matrices with a fixed size.} which corresponds to the selected nodes, as highlighted with colors in Fig. \ref{fig:schematic-weak}(b). 
With these four convolutional operations, we now have a matrix with a size of $5\times5\times4$ ($\BM_{5,5,4}$), as shown in Fig. \ref{fig:schematic-weak}(c).
We then reshape the matrix to an array $25\times4$ ($\BM_{25,4}$), as shown in Fig. \ref{fig:schematic-weak}(d).
Each row of $\BM_{25,4}$ corresponds to the local nodal solution vector inside one finite element, the subdomain $\Omega^e$ in \eref{eq:discretized-residual}, which can then be used to efficiently evaluate the residual via the matrix-vector operation. But it leaves one blank element at each of rows 5, 10,\dots 25.

To evaluate the residual of the steady-state diffusion problem with $2\times2$ Gauss quadrature points, the B-operator matrix in \eref{eq:discretized-residual-diffusion} has a size of $4\times2\times4$ (\# of Gauss quadrature points $\times$ spatial dimensions $\times$ \# of nodes), denoted as $\BB_{4,2,4}$. Its transpose with respect to its last two slots is denoted as $\BB_{4,4,2}^T$.
The bulk residual at each Gauss quadrature point $i$ is evaluated as
\begin{equation}
(\BR_\text{bulk}^{25,4})^i = \omega_i D \BM_{25,4} \BB_{i,4,2} \BB_{i,2,4} 
  \label{eq:bulk-residual-matrix-vector-evaluation}
\end{equation}
with $\omega_i$ denoting the weights. 
The total bulk residual is computed as
\begin{equation}
\BR_\text{bulk}^{25,4} = \sum_{i=1}^{n_\text{quad}} \BR_\text{bulk}^i,
\end{equation}
as shown in Fig. \ref{fig:schematic-weak}(d). $\BR_\text{bulk}^{25,4}$ is then reshaped to $\BR_\text{bulk}^{5,5,4}$, and stored in the element-like form, as shown in Fig. \ref{fig:schematic-weak}(e).

Next, we use the \verb tf.roll ~function to shift the element-like residual to the correct nodal position, as shown Fig. \ref{fig:schematic-weak}(f), with 
\begin{equation}
  \begin{aligned}
    \BR_\text{bulk}^{5,5,0:1} & = \BR_\text{bulk}^{5,5,0:1} \\
    \BR_\text{bulk}^{5,5,1:2} & = \texttt{tf.roll} (\BR_\text{bulk}^{5,5,1:2}, [1], [2])\\
    \BR_\text{bulk}^{5,5,2:3} & = \texttt{tf.roll} (\BR_\text{bulk}^{5,5,2:3}, [1], [1])\\
    \BR_\text{bulk}^{5,5,3:4} & = \texttt{tf.roll} (\BR_\text{bulk}^{5,5,3:4}, [1,1], [1,2])\\
  \end{aligned}
\end{equation}
where the ``='' sign represents an assignment operation, and the ``:'' sign represents the slicing operation in Python. Readers are directed to TensorFlow documentation for the usage of \texttt{tf.roll}.
The assemble operation in \eref{eq:discretized-residual-diffusion} for the bulk integration is now achieved by the \verb tf.reduce_sum ($\BR_\text{bulk}^{5,5,4}$) function without looping over all the elements to get $\BR_\text{bulk}^{5,5,1}$ as done traditionally in the FEM.
Readers are directed to our source code for the implementation of the linear/nonlinear elasticity problems.

\begin{algorithm}[t!]
  \caption{Residual calculation for the steady-state diffusion example. \label{algo:compute-residual}}
  \textbf{Bulk residual with applied Dirichlet BCs:} $\BR_\text{tot}$
  \begin{algorithmic}[1]
    \STATE Apply Dirichlet BCs to NN predicted solutions $\BM^\text{NN}_{5,5}$ by replacing the nodal values at the corresponding locations to obtain $\BM_{5,5}$ (Fig. \ref{fig:schematic-weak}a).
    \STATE Convert $\BM_{5,5}$ from nodal value representation to a four-node element representation $\BM_{5,5,4}$ by convolutional operations with kernels $k_{B,1}, k_{B,2}, k_{B,3}$, and $k_{B,4}$ (Fig. \ref{fig:schematic-weak}b, \ref{fig:schematic-weak}c). For two-dimensional elasticity problems, NN predicted solutions have two channels to represent both the components of the displacement vector $\Bu = \Bvarphi$. The same four kernels are applied to both channels, resulting in the element representation $\BM$ with a third dimension of 8 instead of 4.
    \STATE Get the vectorized representation $\BM_{25,4}$ with each row being the local nodal solutions for one element (Fig. \ref{fig:schematic-weak}d).
    \STATE Compute bulk residual $\BR_{25,4}$ for each element (Fig. \ref{fig:schematic-weak}d). Readers are directed to our source code for details of the bulk residual calculation of linear/nonlinear elasticity.  
    \STATE Switch back to matrix representation of element-like nodal residual $\BR_\text{bulk}^{5,5,4}$ (Fig. \ref{fig:schematic-weak}e).
    \STATE Assemble bulk residual $\BR_\text{bulk}^{5,5,1}$ (Fig. \ref{fig:schematic-weak}f).
  \end{algorithmic}
  \textbf{Residual contributions from Neumann BCs:} $\BR_\text{Neu}$
  \begin{algorithmic}[1]
    \STATE Use kernels $k_{I,1}, k_{I,2}$ and $k_{II,1}, k_{II,2}$ to construct two groups of two-node surface elements $\BI^{5,5,2}_\text{Neu,I}$ and $\BI^{5,5,2}_\text{Neu,II}$ corresponding to the two edges with Neumann BCs. 
    For two-dimensional elasticity problems, we have four groups of surface elements with two each for the traction vector $\BT = \Bk$.
    \STATE Get the vectorized representation of surface elements $\BI^{25,2}_\text{Neu,I}$ and $\BI^{25,2}_\text{Neu,II}$.
    \STATE Compute residuals $\BR_\text{Neu,I}^{25,2}$ and $\BR_\text{Neu,II}^{25,2}$ from Neumann BCs.
    \STATE Switch back to matrix representation of element-like nodal residual $\BR_\text{Neu,I}^{5,5,2}$ and $\BR_\text{Neu,II}^{5,5,2}$.
    \STATE Assemble residual at Neumann BCs $\BR_\text{Neu}^{5,5,1}$.
  \end{algorithmic}
  \textbf{Reduced total residual:} $\BR_\text{tot}^\text{red}$
  \begin{algorithmic}[1]
    \STATE Create a mask matrix $\BM_\text{bulk}^{5,5}$ based on $\BI_\text{D}^{5,5}$ to represent the pixel locations with valid bulk residual values. The entries of $\BM_\text{bulk}^{5,5}$ are zero for the components of $\BI_\text{D}^{5,5}$ with a value of $-1$, which indicates the margins between the true problem domain and the background grid (for more details see Section \ref{sec:data-representation}). For the steady-state diffusion examples, all entries of  $\BM_\text{bulk}^{5,5}$ are one, since the problem domain matches the background grid.
    \STATE Create a reverse mask matrix $\BM_\text{D,rev}^{5,5}$ based on $\BI_\text{D}^{5,5}$ to represent the pixel locations that are not at the Dirichlet boundary. The entries of $\BM_\text{D,rev}^{5,5}$ are zero corresponding to the elements of $\BI_\text{D}^{5,5}$ with a value larger than zero. 
    \STATE Compute total residual $\BR_\text{tot}$ based on \eref{eq:compute-total-residual-example}.
    \STATE Multiply (element-wise) $\BR_\text{tot}$ with $\BM_\text{D,rev}^{5,5}$ and $\BM_\text{bulk}^{5,5}$ to get $\BR_\text{tot}^\text{red}$.
  \end{algorithmic}
\end{algorithm}

\subsection{Neumann BCs}
One channel of the inputs that contains purely Neumann BCs, denoted as $\BI^{5,5}_\text{Neu}$, is shown in Fig. \ref{fig:schematic-weak}(g), where the matrix contains only non-zero entries at the non-zero Neumann boundary locations. 
The Neumann residual needs to be evaluated within surface elements. 
Similar to computing the bulk residual, we apply convolutional operations to $\BI^{5,5}_\text{Neu}$ to construct surface elements. 
Two sets of kernels are used to construct two groups of surface elements, with group I for computing the residual contributions on edges with a surface normal in the X-direction (zero padding is required), and group II for edges with a surface normal in the Y-direction (zero padding is required).
We use the following two kernels 
\begin{equation}
  k_{I,1} =
\begin{bmatrix}
1 & 0\\
0 & 0\\
\end{bmatrix}, 
\quad
k_{I,2} =
\begin{bmatrix}
0 & 0\\
1 & 0\\
\end{bmatrix}.
  \label{eq:kernel-node-to-y-surf-elem}
\end{equation}
to construct surface elements $\BI^{5,5,2}_\text{Neu,I}$ for the first group, with the selected nodal information being shown in Fig. \ref{fig:schematic-weak}(h-I), and the following kernels 
\begin{equation}
  k_{II,1} =
\begin{bmatrix}
1 & 0\\
0 & 0\\
\end{bmatrix}, 
\quad
k_{II,2} =
\begin{bmatrix}
0 & 1\\
0 & 0\\
\end{bmatrix}, 
  \label{eq:kernel-node-to-x-surf-elem}
\end{equation}
to construct surface elements $\BI^{5,5,2}_\text{Neu,II}$ for the second group, with the selected nodal information being shown in Fig. \ref{fig:schematic-weak}(h-II). 

Similar to the bulk residual calculation, we form two matrices, $\BI^{25,2}_\text{Neu,I}$ and $\BI^{25,2}_\text{Neu,II}$, to compute the Neumann residual.
We use two Gauss quadrature points for surface integration.
The shape function $\BN$ in \eref{eq:discretized-residual-diffusion} has a size of $2\times2$ (\# of Gauss quadrature points $\times$ \# of nodes), and is denoted by $\BN_{2,2}$.
We evaluate the Neumann residual at each Gauss quadrature point $i$ via
\begin{equation}
  (\BR_\text{Neu,I}^{25,2})^i = \omega_i \BI^{25,2}_\text{Neu,I} \BN_{i,2} \BN_{i,2}
  \quad \text{and} \quad
  (\BR_\text{Neu,II}^{25,2})^i = \omega_i \BI^{25,2}_\text{Neu,II} \BN_{i,2} \BN_{i,2} 
\end{equation}
with $\omega_i$ denoting the weights. The total Neumann residual is computed as
\begin{equation}
  \BR_\text{Neu,I}^{25,2} = \sum_{i=1}^{n_\text{sq}} \BR_\text{Neu,I}^i
  \quad \text{and} \quad
  \BR_\text{Neu,II}^{25,2} = \sum_{i=1}^{n_\text{sq}} \BR_\text{Neu,II}^i.
\end{equation}
where $n_\text{sq}$ is the number of surface quadrature points. Again, we use the \verb tf.roll ~function to unfold the element-like residual to the correct nodal position, similar to those shown Fig. \ref{fig:schematic-weak}(f), for group I 
\begin{equation}
  \begin{aligned}
    \BR_\text{Neu,I}^{5,5,0:1} & = \BR_\text{Neu,I}^{5,5,0:1} \\
    \BR_\text{Neu,I}^{5,5,1:2} & = \texttt{tf.roll} (\BR_\text{Neu,I}^{5,5,1:2}, [1], [1])\\
  \end{aligned}
\end{equation}
and for group II
\begin{equation}
  \begin{aligned}
    \BR_\text{Neu,II}^{5,5,0:1} & = \BR_\text{Neu,II}^{5,5,0:1} \\
    \BR_\text{Neu,II}^{5,5,1:2} & = \texttt{tf.roll} (\BR_\text{Neu,II}^{5,5,1:2}, [1], [2]).\\
  \end{aligned}
\end{equation}
The assemble operation in \eref{eq:discretized-residual-diffusion} for the surface integration is now achieved by the \verb tf.reduce_sum ($\BR_\text{Neu,I}^{5,5,2}$) and \verb tf.reduce_sum ($\BR_\text{Neu,II}^{5,5,2}$) without looping over elements.
We obtain the final residual contributions from the Neumann BCs:
\begin{equation}
  \BR_\text{Neu}^{5,5,1} = \BR_\text{Neu,I} + \BR_\text{Neu,II}.
\end{equation}
The total residual $\BR_\text{tot}$ in \eref{eq:discretized-residual-diffusion}, as shown in Fig. \ref{fig:NN}, is computed as
\begin{equation}
  \BR_\text{tot}^{5,5,1} = \BR_\text{bulk}^{5,5,1} - \BR_\text{Neu}^{5,5,1} 
  \label{eq:compute-total-residual-example}
\end{equation}
by applying the Neumann residual to the bulk residual.
To construct the deterministic loss in \eref{eq:loss-deterministic} and the likelihood function in \eref{eq:likelihood-function}, the reduced residual $\BR_\text{tot}^\text{red}$ obtained by excluding the residual at the Dirichlet boundary location from $\BR_\text{tot}$ is used, as shown Fig. \ref{fig:schematic-weak}(i).
It is worth mentioning that additional auxiliary matrix/vector/tensor operations have been introduced, which are not included in the description, to complete this efficient residual evaluation. Readers are invited to refer to our code for the detailed implementation.

\section{Data representation and numerical aspects} \label{sec:data}
To demonstrate the performance of our methods, we investigate different definitions of BVPs for the three considered physical systems. 
We prepared both small datasets, which contain one or multiple BVPs, and large datasets, which could contain hundreds of thousands BVPs.
The NN inputs corresponding to BVPs in these datasets are synthetically generated to train the discretized residual-constrained NNs. 
To compare the solution accuracy between NNs and DNSs, we also solve these BVPs with \verb mechanoChemFEM, \footnote{Code available at \href{https://github.com/mechanoChem/mechanoChemFEM}{github.com/mechanoChem/mechanoChemFEM}.} which is a publicly available multiphysics code developed by us based on the deal.II library \cite{dealII90}.
In general, for small datasets, the number of unique sets of BCs is much smaller than the number of parameters of the NNs that represent the PDE solution. 
We therefore augment the unique sets of BVPs by duplicating them multiple times to form an augmented dataset during training.
In the remaining part of this section, we present  details on the data structure of NN inputs, domain/boundary detection, and the NN training procedure.

\subsection{Data structure of NN inputs} \label{sec:data-representation}
\begin{figure}[h!]
  \centering
  \includegraphics[width=0.75\linewidth]{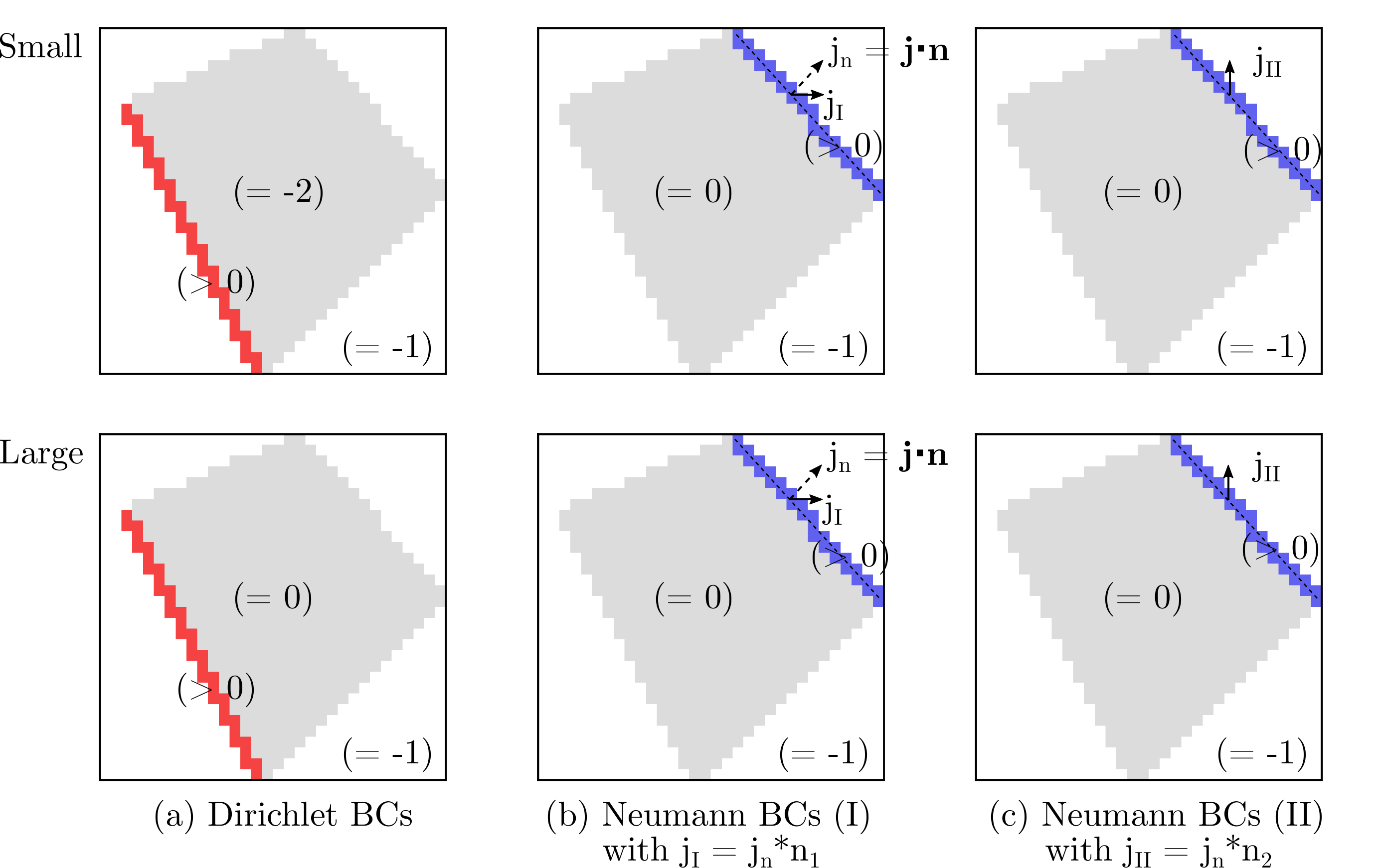}
  \caption{Illustration of the data structure of NN inputs for a steady-state diffusion problem. The NN inputs contain three channels. (a) the Dirichlet BCs (red),  (b,c) the Neumann BCs (blue) on edges with a surface normal in the X- and Y-direction, respectively.\protect\footnotemark~
    Only the boundary locations have values that are greater than 0, which is physically meaningful.
    Top row: input structure for small datasets, which contain one or multiple BVPs.
    Bottom row: input structure for large datasets, which could contain hundreds of thousands BVPs.
    In the first channel, the problem domain (gray color) is filled with a value of $-2$ (top) or $0$ (bottom). If the value is $-2$, the region is filled with random numbers during the training process. 
    The margin (white region) is filled with a value of $-1$. 
    The residual contribution from this region is excluded when computing $\BR_\text{tot}^\text{red}$.
    In the second and third channel, the problem domain is filled with a value of $0$. 
    Similarly, the margin is filled with a value of $-1$. 
    We use $\bj=[j_\text{I},~j_\text{II}]$ and $\bn=[n_1,~n_2]$ to denote the surface flux and surface normal, with $j_\text{I}$ and $j_\text{II}$ being the projected values in X- and Y-direction, respectively.
  }
  \label{fig:data-diffusion-illustration}
\end{figure}
\footnotetext{For elasticity problems, the inputs contain four channels, with the first two representing Dirichlet BCs for $\bar{u}_x$ and $\bar{u}_y$ and the last two presenting Neumann BCs for $\bar{T}_x$ and $\bar{T}_y$.}

Since the discretized residual constrained NNs do not require labels for training, the NN inputs are synthetically generated  with only information on problem domains and the applied BCs.
We consider a fixed square background grid of $[0,~1]\times[0,~1]$, with \verb nx ~ and \verb ny ~ total pixels in each dimension.
For both diffusion and elasticity problems, each input data point is a three-dimensional matrix $\BI_\text{nx,ny,3$\times$DOF}$ to represent a set of BCs. The first two indices of $\BI$ indicate the pixels locations in X- and Y- directions. 
For steady-state diffusion problem with one scalar DOF per node, there are three channels in the third dimension, which contain information of Dirichlet BCs, Neumann BCs on the edge with a surface normal in the X-direction, and Neumann BCs on the edge with a surface normal in the Y-direction, respectively.
For elasticity problems, there are six channels in the third dimension with the first two channels containing Dirichlet BCs in X- and Y- directions, the third and fourth channels containing Neumann BCs in X- and Y- directions on the edge with a surface normal in the X-direction, and the fifth and sixth channels containing Neumann BCs in X- and Y- directions on the edge with a surface normal in the Y-direction, respectively.
Data normalization between $[-1,~1]$ is used to ensure that all the physically meaningful data in our study has a value greater than 0.

The structure of the input data is illustrated in Fig. \ref{fig:data-diffusion-illustration} with the diffusion problem as an example.
In our study, the problem domain does not necessarily occupy the whole background grid, which results in the margin region as shown in Figs \ref{fig:data-diffusion-illustration} and \ref{fig:bvp-problem-setup}.
In small datasets, for the channel(s) containing Dirichlet BCs, the problem domain is filled with $-2$ except the Dirichlet boundary values, which is greater than 0.
The auxiliary number $-2$ exists in augmented datasets and serves as an indicator to be filled with random numbers during the training process. 
For the margin region, which represents the space between the background grid and the problem domain, if there is any, is filled with $-1$. 
The auxiliary number\footnote{The auxiliary numbers $-1$ and $-2$ are arbitrary choices with no physical meaning. Users can choose different values to assign to the margin and the problem domain for the inputs.} $-1$ serves as an indicator to evaluate $\BR_\text{tot}^\text{red}$ with the residual in this region being excluded.
In large datasets, for the channel(s) containing Dirichlet BCs, the problem domain is filled with $0$ except the Dirichlet boundary values, which is greater than 0.
For the channel(s) containing Neumann BCs, the problem domain is filled with a value of $0$ except the Neumann boundary values. 
When convolutional kernels operate on the problem domain, only the boundary makes a non-zero contribution. 
Similarly, the margin  is filled with a value of $-1$ for assisting the calculation of $\BR_\text{tot}^\text{red}$.
Examples of the actual inputs for steady state diffusion are shown in Fig. \ref{fig:diffusion-20bvp-results}(a,b).

\subsection{Domain/Boundary detection}\label{sec:data-domain-bc-detection}
As discussed in Section \ref{sec:data-representation}, a fixed value of $-1$ is assigned to the margins.
When calculating the residual, a mask matrix is created for domain detection.
This mask matrix is created based on the information on Dirichlet BCs from the inputs and ensures that only the residual inside the actual problem domain is evaluated.
Readers are directed to Refs \cite{Bhatnagar2019CNN-encoder-decoder,Li2020Reaction-diffusion-prediction-CNN,Gao2020PhyGeoNet-PDE-on-Irregular-domain} and many others for approaches to map complex and irregular domains onto a regular mesh. 
Such geometric transformations can be easily taken into account in the proposed PDE loss layers with the parent to physical domain mapping commonly used in the FEM via basis functions. 
The proposed approach, using a mask matrix for domain detection, should still be applicable to other parametric domain representations, though it is not the focus of this work.

In our study, each input data point represents a unique BVP for a specific problem domain and set of applied BCs.
To detect the Dirichlet BCs in small datasets, during the NN training, the input augmented data is first passed to a customized Keras layer, called \verb LayerFillRandomNumber, which fills the pixel locations with values of  $-2$ in the Dirichlet BCs channel with uniformly generated random numbers in the range of $[0,~1]$ to ensure that all the data points are independent from each other.
As the problem domain is filled with random numbers, the convolutional kernels iteratively learn the actual Dirichlet boundary values.
The data structure in the Neumann BCs channel ensures that the kernels learn to recognize the information on the Neumann BCs, as the problem domain is filled with zeros.
For large datasets, the interior domain of the Dirichlet BCs channel is filled with zero. The convolutional kernels will learn the actual problem domain, boundary locations, and boundary values.

\subsection{Neural network training} \label{sec:NN-training}

\begin{figure}[t!]
  \centering
  \includegraphics[width=1.0\linewidth]{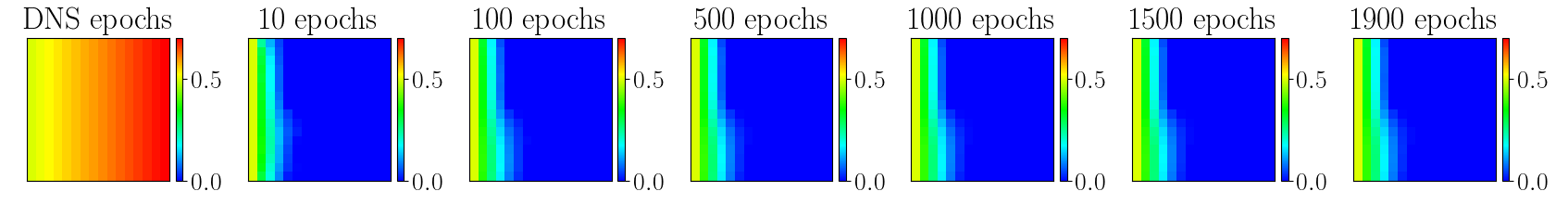}
  \includegraphics[width=1.0\linewidth]{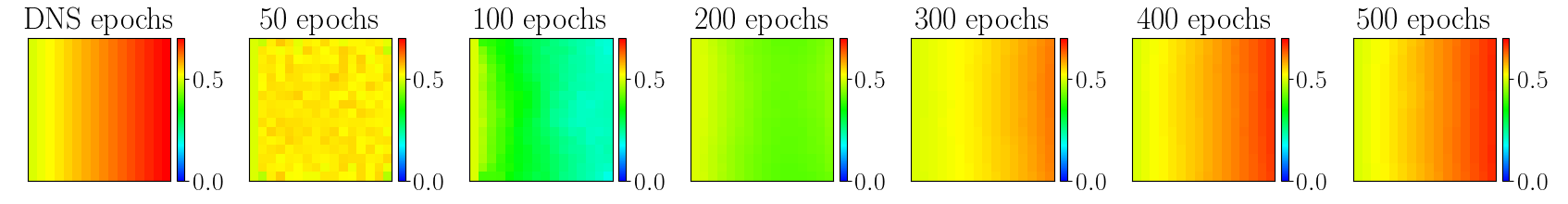}
  \caption{Illustration of the deterministic NNs predicted solution at different epochs for a diffusion BVP setup with domain id 5 and BCs id 2 (flux loading from the right edge), as shown in Fig. \ref{fig:bvp-problem-setup}(a). Top: without zero initialization. Bottom: with zero initialization for the first 100 epochs. } 
  \label{fig:zero-initialization}
\end{figure}

\begin{algorithm}[t!]
  \caption{Training procedure for deterministic NNs. \label{algo:cnn-training}}
  \begin{algorithmic}[1]
    \STATE Load NN inputs with each data point being a unique set of BCs.
    \STATE Use a large dataset $\calD$ or an augmented dataset $\calD$ by by duplicating the small dataset multiple times.
    \STATE Split $\calD$ into training, validation, and testing datasets. 
    \STATE Setup the encoder-decoder deterministic NN structure, with the first layer being a customized layer to fill the locations that have values of $-2$ in $\calD$ with uniform random numbers between $[0,1]$ to ensure $\calD$ is i.i.d.
    \FOR{epoch < total\_epochs}
    \STATE Batch train the NNs
    \IF{use zero initialization \AND epoch < total\_zero\_initialization\_epochs}
    \STATE Use dummy labels with values of 0.5, which is equivalent to an actual zero before data normalization, to form the MSE loss to train the NN.
    \ELSE
    \STATE Use $\BR_\text{tot}^\text{red}$ to form the deterministic loss to train the NN.
    \ENDIF
    \ENDFOR
    \STATE Make prediction.
  \end{algorithmic}
\end{algorithm}

For deterministic NNs, a fixed learning rate is used to batch optimize the loss function \eref{eq:loss-deterministic} and solve the PDE systems.
In our study, we found that pure Dirichlet problems are learned (converge) faster than Dirichlet-Neumann problems. 
In the latter case, the solution the NNs fail to learn the solution in some instances. 
This observation holds for all three PDEs considered here. 
This is mainly because, for Dirichlet-Neumann problems, it is the gradient of the unknown field(s) that drives the loss instead of the field(s) itself (themselves) as in the case of the Dirichlet problem.
We demonstrate this by showing the NN predicted solution at different epochs in Fig. \ref{fig:zero-initialization} for a diffusion BVP  with domain id 5 and BCs id 2 (see Fig. \ref{fig:bvp-problem-setup}a) with zero concentration on the left edge and non-zero flux on the right edge.
The top row of Fig. \ref{fig:zero-initialization} shows that even though the NN predicted concentration changes, it does so very slowly via a front progressing from the left edge (zero Dirichlet BCs) to the right edge (flux BCs). This indicates that the solution has not yet been learned with an accuracy that is comparable to the DNS results.

This difficulty arises mainly because the parameters of NNs are randomly initialized. As a result, the
NN predicted solutions at early training stages are random numbers close to zero. Since data normalization is used, the NN solution of zero corresponds to an actual value of $-1$. 
Such random outputs in early stages can violate the governing equations, potentially in catastrophic manner e.g. resulting in a deformation gradient with negative determinant in nonlinear elasticity. 
In order to rectify this inconsistency, we draw from the conventional initialization of the solution vector to zero in the FEM, and adopt the same approach for the NNs. 
For the first few epochs, we train NNs with dummy labels with values of 0.5 (equivalent to an actual value of 0) without enforcing the PDE constraint. 
We call this as the zero initialization process.
This process helps to improve the initialization of NN parameters. 
After the zero initialization procedure is completed, the PDE constraints are enabled to train the NNs to solve the PDE systems. 
We found this remedy to drastically improve the learned solution as well as speed up the training, as shown in the bottom row of Fig. \ref{fig:zero-initialization}, where the NN predicted solutions approach the DNS results at 500 epochs, much faster than the case without the zero initialization process.
The training process for deterministic NNs is summarized in the Algorithm Box \ref{algo:cnn-training}. 

For probabilistic NNs, we can use the proposed approach successfully solve a single BVP. 
However, when we try to solve multiple BVPs, we notice that the BNNs converge faster to purely Dirichlet problems (boundary id 1, 3 for the diffusion problem and boundary id 1, 4 for elasticity problems) than those with non-zero Neumann BCs.
Once the BNN parameters stagnate at sub-optimal solutions, it is very difficult to optimize them further for other BVPs with Neumann BCs.
To overcome this challenge, we first train deterministic NNs with identical architectures as the BNNs. 
Once the deterministic NNs have converged to a desired tolerance, we then initialize the mean of the posterior distribution of parameters in the BNNs with the optimized parameters from the deterministic model.
We refer to this as the optimal parameter initialization process.
During the subsequent training of the BNNs we use a small learning rate to explore the local parameter space around these optimized parameters. A similar approach has been adopted in Ref \cite{Geneva2020Zabaras-JCP-auto-regressive-NN-PDE}, 
The training process for BNNs is summarized in the Algorithm Box \ref{algo:bnn-training}. 

Here an epoch corresponds to a single application of an augmented dataset or a large dataset $\calD$ as inputs to predict a NN solution by training parameters.
For augmented dataset, when splitting $\calD$ into training, validation, and testing groups, each group could potentially contain all the unique BCs. 
The evaluation of the NN results based on such validation and testing dataset only indicates how well the NNs solve the exposed sets of BCs.
To test the predictability of the trained NNs in Sections \ref{sec:linear-lshape} and \ref{sec:nonlinear-rectangle-interpolation}, the unseen testing sets of BCs were set aside before data augmentation. 

\begin{algorithm}[t!]
  \caption{Training procedure for BNNs. \label{algo:bnn-training}}
  \begin{algorithmic}[1]
    \STATE Load NN inputs with each data point being a unique set of BCs.
    \STATE Use a large dataset $\calD$ or an augmented dataset $\calD$ by by duplicating the small dataset multiple times.
    \STATE Split $\calD$ into training, validation, and testing datasets. 
    \STATE Setup the encoder-decoder probabilistic NN structure, with the first layer being a customized layer to fill the locations that have values of $-2$ in $\calD$ with uniform random numbers between $[0,1]$.
    \IF{use optimal parameter initialization}
      \STATE Load the optimized parameters from the deterministic NNs to initialize the mean of the posterior distribution of BNN parameters.
    \ELSE
      \STATE Use random initialization for the posterior distribution of BNN parameters.
    \ENDIF
    \FOR{epoch < total\_epochs}
    \STATE Batch train the NNs
    \IF{use zero initialization \AND epoch < total\_zero\_initialization\_epochs \AND (\NOT use optimal parameter initialization)}
    \STATE Use dummy labels with values of 0.5, which is equivalent to an actual zero before data normalization, to form the MSE loss to train the NN.
    \ELSE
    \STATE Use $\BR_\text{tot}^\text{red}$ to form the likelihood loss to train the NN.
    \ENDIF
    \ENDFOR
    \STATE MC sampling for UQ.
  \end{algorithmic}
\end{algorithm}

\subsection{Data generation for large dataset}
We use the \verb scikit-geometry $~$package to generate the Polygons that were used for the large dataset study, where the vertices of polygons are randomly sampled approximately along a circle with applied perturbation. The detailed code implementation is available on GitHub. To choose the edges to impose the BCs, we first generate all the possible combinations of two non-adjacent edges. We then apply the \verb numpy.random.shuffle() $~$to the combinations and select the first few of them from the combinations. For the boundary values, we use the \verb numpy.random.uniform() function to sample control points from the range $[0,~1]$. If the boundary value has a constant distribution, a single value is sampled. If the boundary value has a linear distribution, two extremes are sampled. And interpolated values are imposed to each pixel on the edge. One can refer to our code on Github for details to impose BCs with a quadratic or sinusoidal distribution. 

\section{Numerical results}
In this section, we provide additional information for the examples presented in the main text along with extra examples to support the performance of the proposed method. 

\subsection{Steady state diffusion - small dataset}
\subsubsection{Single octagon domain with Dirichlet BCs - cold start versus warm start}\label{sec:results-diffusion-octagon-cold-warm}

\begin{figure}[h!]
  \centering
  \subfloat[]{\includegraphics[width=0.15\linewidth]{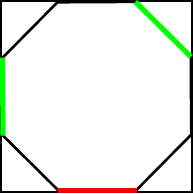}}
  \hspace{10mm}
  \subfloat[]{\includegraphics[width=0.15\linewidth]{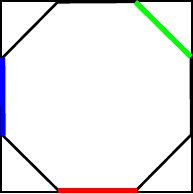}}
  \caption{Illustration of the setup of BVPs with (a) purely Dirichlet BCs and (b) mixed BCs}
  \label{fig:oct-bc-setup}
\end{figure}

\begin{figure}[h!]
  \centering
  \includegraphics[width=1.0\linewidth]{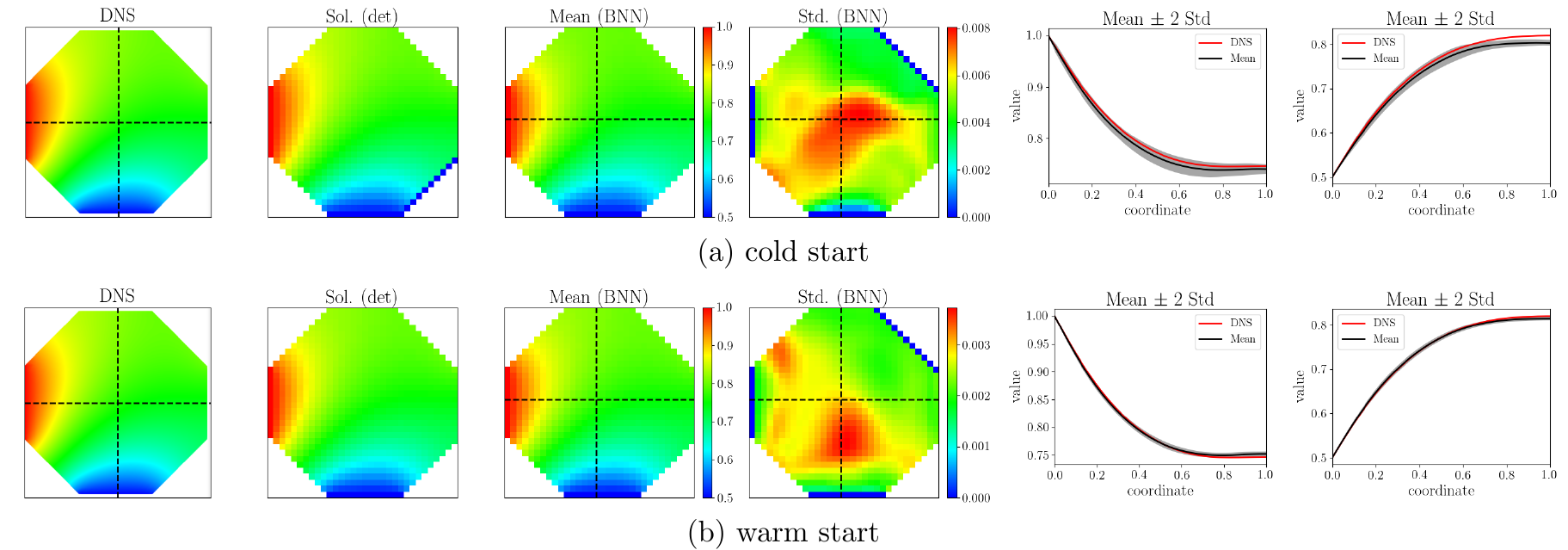}
  \caption{Comparison of BNN results between (a) cold start and (b) warm start with purely Dirichlet BCs. The results are very similar, except the uncertainty level in (a) is higher than (b).}
  \label{fig:cold-warm-results-comparison}
\end{figure}

In this example, we use the proposed method to solve the steady-state diffusion problem on an octagon domain with purely Dirichlet BCs, as shown in Fig. \ref{fig:oct-bc-setup}(a). The BNN is trained in two cases with (i) cold start and (ii) warm start. The results are shown in Fig. \ref{fig:cold-warm-results-comparison}. One can see that the BNNs with cold start can also find the correct solution for this specific BVP, comparable with the results from the warm start. This confirms that the loss of the BNN is correct. However, in our study, we found that it is very challenging, if not impossible, for BNNs to directly (cold start) solve multiple BVPs. All the other BNNs results in this work are solved based on the warm start approach. 

\subsubsection{Single octagon domain with mixed BCs}\label{sec:results-diffusion-octagon}

In Fig. \ref{fig:diffusion}(a-g), we use the proposed PDE constrained NNs to solve steady-state diffusion on an octagonal domain with mixed BCs, resulting in a solution with spatially varying gradients along both X- and Y- directions.
NN structure information for the octagonal domain simulation are summarized in Table \ref{tab:diffusion-octagon-NNs} and \ref{tab:diffusion-octagon-NNs-others}.
The NN hyperparameters are manually tuned to achieve a desired performance.

\begin{figure}[t!]
  \centering
  \subfloat[deterministic loss]{\includegraphics[width=0.33\linewidth]{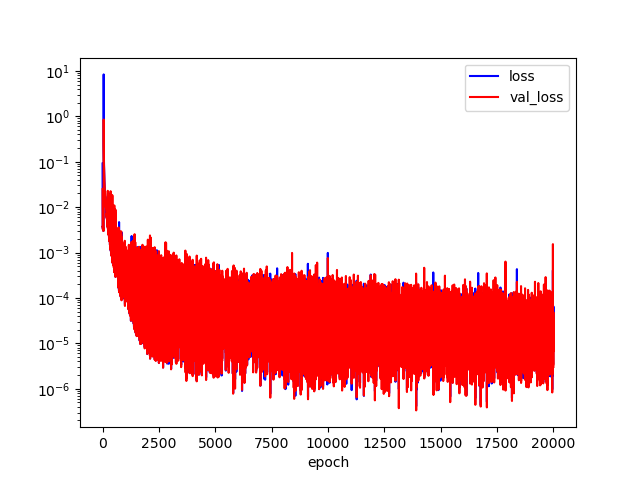}}
  \subfloat[probabilistic loss]{\includegraphics[width=0.33\linewidth]{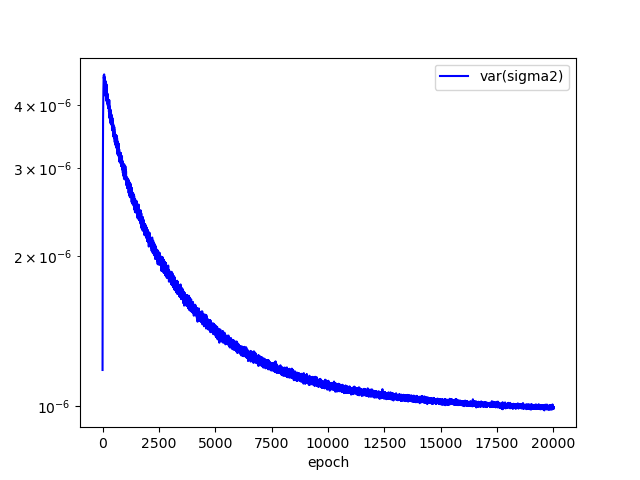}}
  \caption{NN training information for octagon problem with mixed BCs. (a) Loss from the deterministic NN. (b) Loss from the BNN. (c) Evolution of $\Sigma_2$. 
}
  \label{fig:diffusion-oct-32-mixed-loss}
\end{figure}

The training losses for both deterministic and probabilistic NNs are given in Fig. \ref{fig:diffusion-oct-32-mixed-loss}(a, b).
Depends on the choice of the initial value of $\Sigma_2$, Fig. \ref{fig:diffusion-oct-32-mixed-loss}(b) could have a negative values. The negative value is reasonable because the total loss of BNNs in \eref{eq:loss-probabilistic} consists two terms.
The first term in \eref{eq:loss-probabilistic} is non-negative, whereas the second term could be either positive or negative depending on values of both $\BR_\text{tot}^\text{red}$ and $\Sigma_2$. 
The evolution of $\Sigma_2$ from the BNN is shown in Fig. \ref{fig:diffusion-oct-32-mixed-loss}(c), which converges to a sharp value during training. 
The evolution of $\Sigma_2$ is correlated to the sign change of the BNN loss.
The logarithm of the probability density function converges to the maximum, or the negative of it converges the minimum, as shown in Fig. \ref{fig:diffusion-oct-32-mixed-loss}(b), when $\calN \left(0, \Sigma_2 \right)$ best represents $ R_\text{tot}^\text{red}$.
Such behavior is expected as the BNN is initialized with optimal parameters from the deterministic NNs and is trained with a very small learning rate to only explore the local parameter space around the optimized parameters.

\begin{table}
  \centering
  \begin{tabular}{l | l | l | l}
    \hline
    Deterministic         & Probabilistic         & Size         & Layer arguments \\ \hline
    Input                 & Input                 & -            & - \\
    LayerFillRandomNumber & LayerFillRandomNumber & -            & - \\
    Conv2D                & Convolution2DFlipout  & filters = 8  & kernel (5,5), padding: same, ReLU \\
    MaxPooling2D          & MaxPooling2D          & -            & kernel (2,2), padding: same\\
    Conv2D                & Convolution2DFlipout  & filters = 8  & kernel (5,5), padding: same, ReLU \\
    MaxPooling2D          & MaxPooling2D          & -            & kernel (2,2), padding: same\\
    Conv2D                & Convolution2DFlipout  & filters = 8  & kernel (5,5), padding: same, ReLU \\
    MaxPooling2D          & MaxPooling2D          & -            & kernel (2,2), padding: same\\
    Flatten               & Flatten               & -            & - \\
    Dense                 & DenseFlipout          & units = 32   & ReLU \\
    Dense                 & DenseFlipout          & units = 32   & ReLU \\
    Reshape               & Reshape               & -            & $[4,4,4]$ \\
    Conv2D                & Convolution2DFlipout  & filters = 8  & kernel (5,5), padding: same, ReLU \\
    UpSampling2D          & UpSampling2D          & -            & size (2,2) \\
    Conv2D                & Convolution2DFlipout  & filters = 8  & kernel (5,5), padding: same, ReLU \\
    UpSampling2D          & UpSampling2D          & -            & size (2,2) \\
    Conv2D                & Convolution2DFlipout  & filters = 8  & kernel (5,5), padding: same, ReLU \\
    UpSampling2D          & UpSampling2D          & -            & size (2,2) \\
    Conv2D                & Convolution2DFlipout  & filters = 16 & kernel (5,5), padding: same, ReLU \\
    Conv2D                & Convolution2DFlipout  & filters = 1  & kernel (5,5), padding: same, ReLU \\
    \hline
  \end{tabular}
  \caption{Details of both deterministic and probabilistic NNs for solving diffusion BVPs on the octagon domain with an output resolution of $32\times 32$.}
  \label{tab:diffusion-octagon-NNs}
\end{table}

\begin{table}
  \centering
  \begin{tabular}{l | l | l }
    \hline
    Description                     & Deterministic                 & Probabilistic       \\ \hline
    Total parameters                & 16,049                        & 31,970                 \\
    Size of $\calD$                 & 1 $\times$ Aug: $2^{12}$      & 1 $\times$ Aug: $2^{11}$      \\
    Epochs                          & 20,000                        & 100                  \\
    Zero initialization epochs      & 100                           & -                     \\
    Optimizer                       & Nadam                         & Nadam                 \\
    Learning Rate                   & 2.5e-4                        & 1e-8                  \\
    Batch Size                      & 256                           & 64                    \\
    $\Sigma_1$                      & -                             & 1e-8                  \\
    Initial value of $\Sigma_2$     & -                             & 1e-8                  \\
    \hline
  \end{tabular}
  \caption{Training related parameters for solving steady-state diffusion on the octagonal domain. Aug: data augmentation. }
  \label{tab:diffusion-octagon-NNs-others}
\end{table}

\subsubsection{Single octagon domain with mixed BCs and different material parameters}\label{sec:results-diffusion-octagon-different-parameters}

\begin{figure}[h!]
  \centering
  \includegraphics[width=1.0\linewidth]{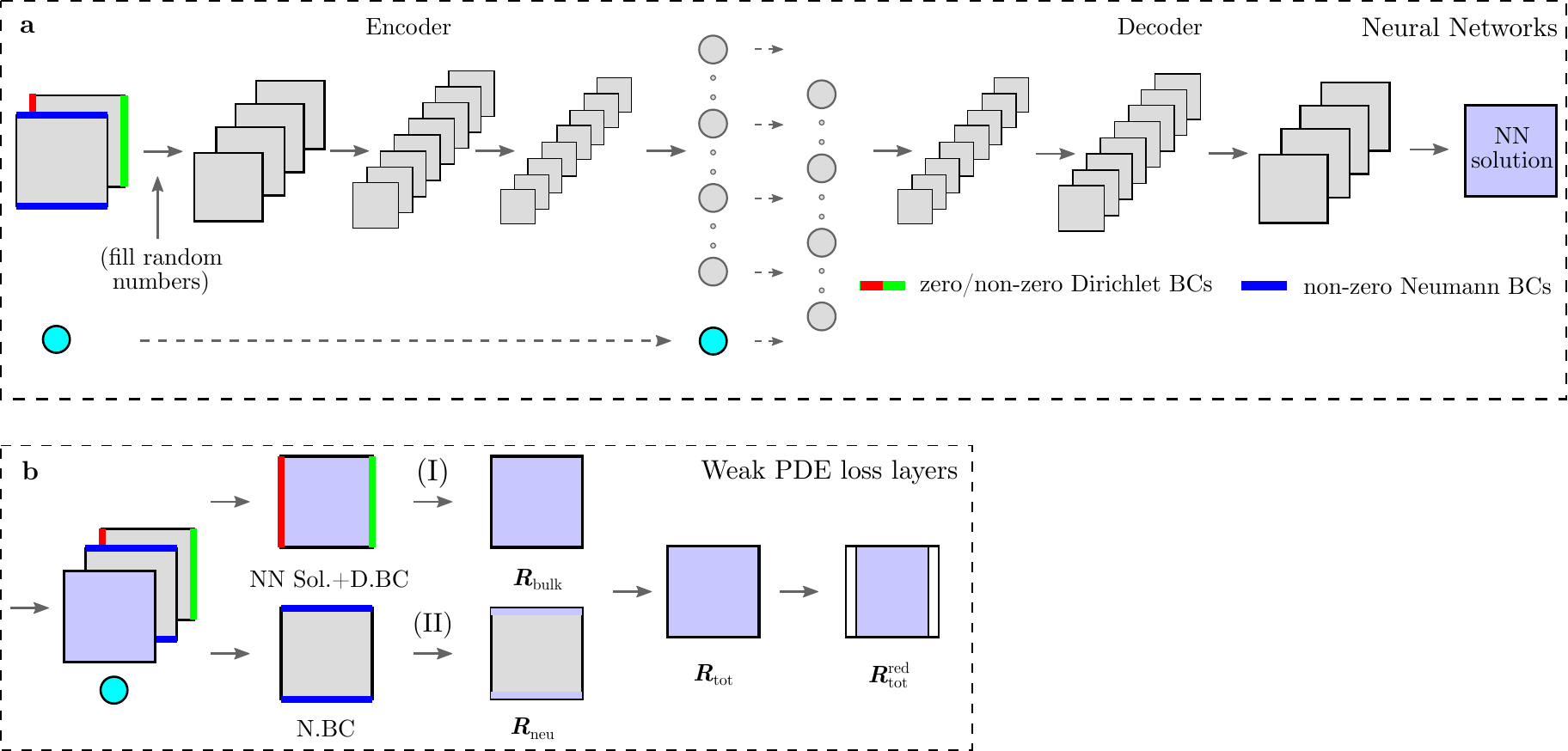}
  \caption{Illustration of the NN architectures with heterogeneous data inputs to account for different material parameters. The image-type input contains the information of problem domains and BCs. The scalar(s) represent the material parameters.}
  \label{fig:NN-heter-ill}
\end{figure}

In this example, we use the proposed method to solve steady-state diffusion problem on an octagon domain with mixed BCs and varying diffusivity.
The NN architecture is shown in Fig. \ref{fig:NN-heter-ill} with heterogeneous inputs, which consists of images and scalars. The image-type input contains the information of problem domains and BCs. The scalar(s) represent the material parameters. 

\subsubsection{Multiple rectangular domains with different BCs}\label{sec:results-diffusion-rectangle}
\begin{table}
  \centering
  \begin{tabular}{l | l | l | l}
    \hline
    Deterministic         & Probabilistic         & Size         & Layer arguments \\ \hline
    Input                 & Input                 & -            & - \\
    LayerFillRandomNumber & LayerFillRandomNumber & -            & - \\
    Conv2D                & Convolution2DFlipout  & filters = 8  & kernel (5,5), padding: same, ReLU \\
    MaxPooling2D          & MaxPooling2D          & -            & kernel (2,2), padding: same\\
    Conv2D                & Convolution2DFlipout  & filters = 16 & kernel (5,5), padding: same, ReLU \\
    MaxPooling2D          & MaxPooling2D          & -            & kernel (2,2), padding: same\\
    Conv2D                & Convolution2DFlipout  & filters = 16 & kernel (5,5), padding: same, ReLU \\
    MaxPooling2D          & MaxPooling2D          & -            & kernel (2,2), padding: same\\
    Flatten               & Flatten               & -            & - \\
    Dense                 & DenseFlipout          & units = 64   & ReLU \\
    Dense                 & DenseFlipout          & units = 64   & ReLU \\
    Reshape               & Reshape               & -            & $[4,4,4]$ \\
    Conv2D                & Convolution2DFlipout  & filters = 16 & kernel (5,5), padding: same, ReLU \\
    UpSampling2D          & UpSampling2D          & -            & size (2,2) \\
    Conv2D                & Convolution2DFlipout  & filters = 16 & kernel (5,5), padding: same, ReLU \\
    UpSampling2D          & UpSampling2D          & -            & size (2,2) \\
    Conv2D                & Convolution2DFlipout  & filters = 16 & kernel (5,5), padding: same, ReLU \\
    Conv2D                & Convolution2DFlipout  & filters = 1  & kernel (5,5), padding: same, ReLU \\
    \hline
  \end{tabular}
  \caption{Details of both deterministic and probabilistic NNs for solving 20 steady-state diffusion BVPs. Readers are directed to TensorFlow documentation for detailed description of the functionality of each layer.}
  \label{tab:diffusion-20bvp-NNs}
\end{table}

\begin{table}
  \centering
  \begin{tabular}{l | l | l }
    \hline
    Description                     & Deterministic                 & Probabilistic         \\ \hline
    Total parameters                & 33,209                        & 66,202                 \\
    Size of $\calD$                 & 20 $\times$ Aug: $2^{10}$     & 20 $\times$ Aug: $2^{9}$      \\
    Epochs                          & 20,000                        & 5,000                 \\
    Zero initialization epochs      & 100                           & -                     \\
    Optimizer                       & Nadam                         & Nadam                 \\
    Learning Rate                   & 2.5e-4                        & 1e-8                  \\
    Batch Size                      & 256                           & 64                    \\
    $\Sigma_1$                      & -                             & 1e-8                  \\
    Initial value of $\Sigma_2$     & -                             & 1e-8                  \\
    \hline
  \end{tabular}
  \caption{Training related parameters for solving 20 steady-state diffusion BVPs. Aug: data augmentation. }
  \label{tab:diffusion-20bvp-NNs-others}
\end{table}

\begin{figure}[h!]
  \centering
  {\includegraphics[height=0.08\linewidth]{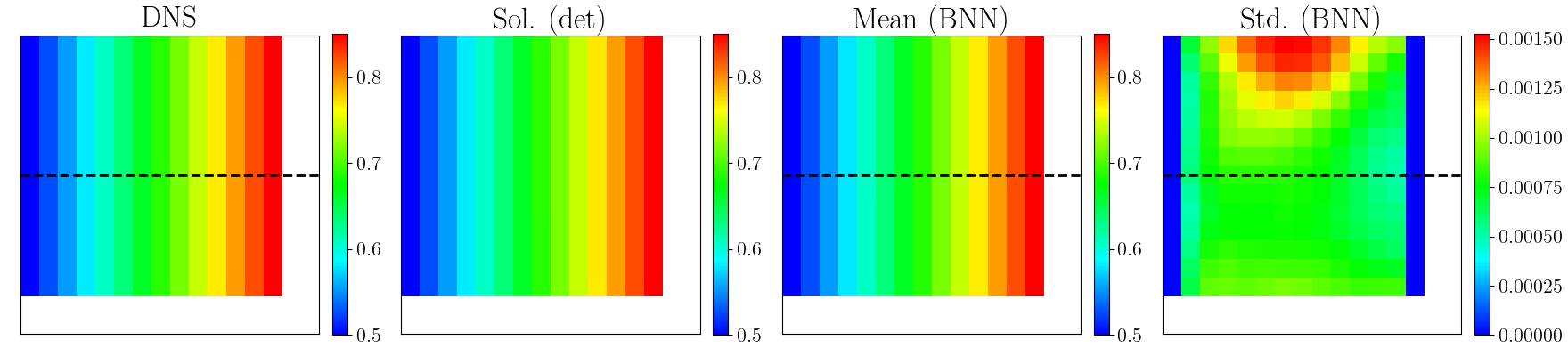}} 
  {\includegraphics[height=0.08\linewidth]{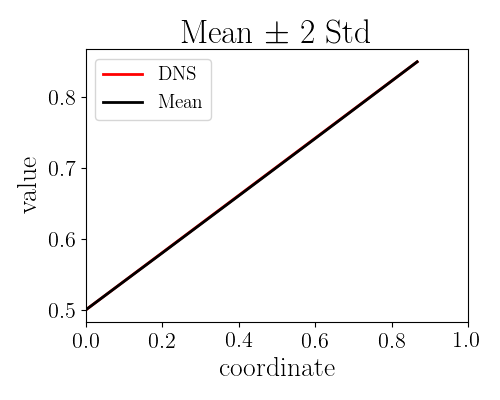}} 
  {\includegraphics[height=0.08\linewidth]{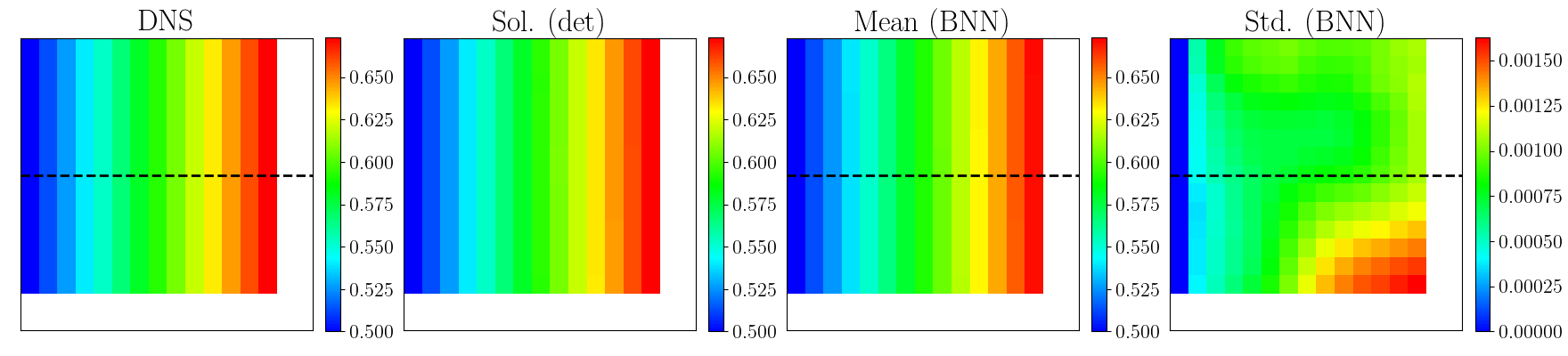}} 
  {\includegraphics[height=0.08\linewidth]{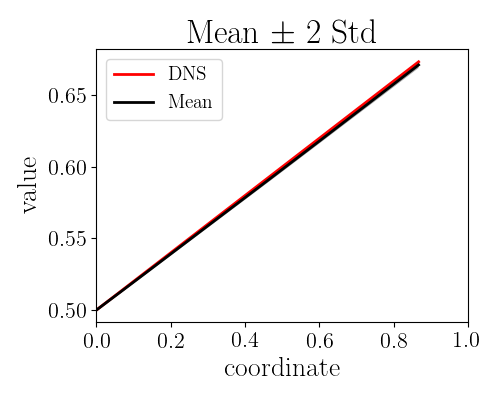}} \\\
  {\includegraphics[height=0.08\linewidth]{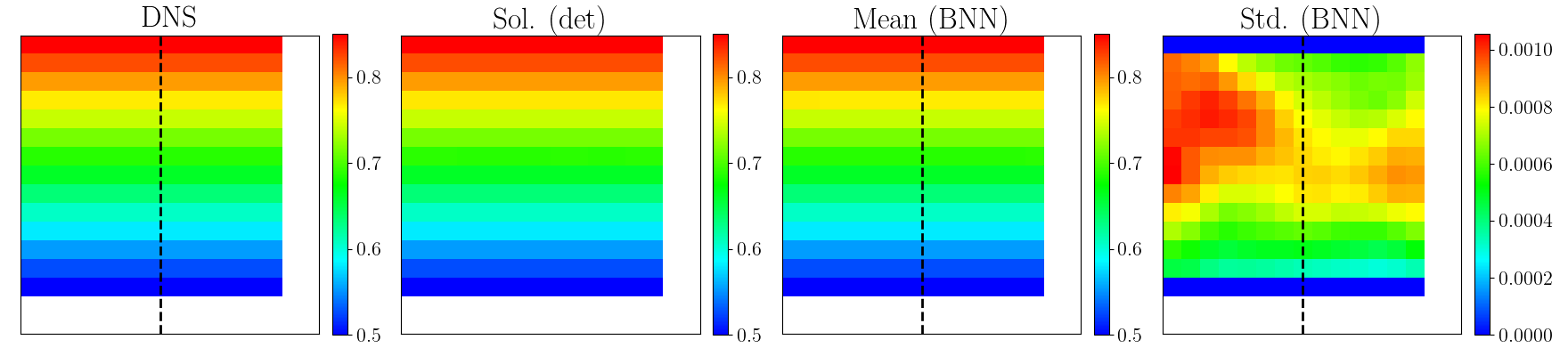}} 
  {\includegraphics[height=0.08\linewidth]{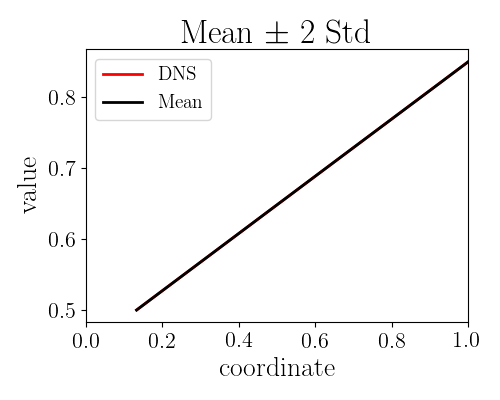}} 
  {\includegraphics[height=0.08\linewidth]{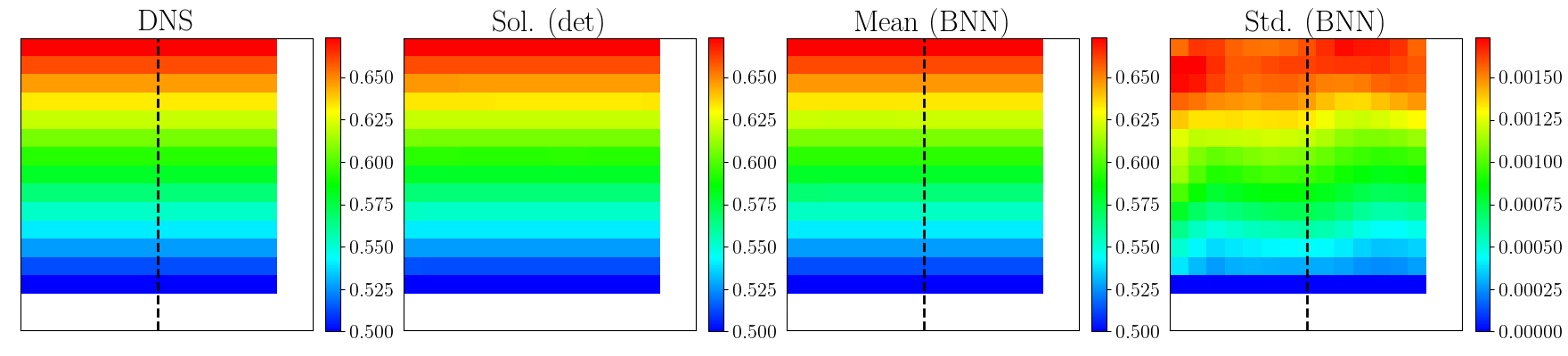}} 
  {\includegraphics[height=0.08\linewidth]{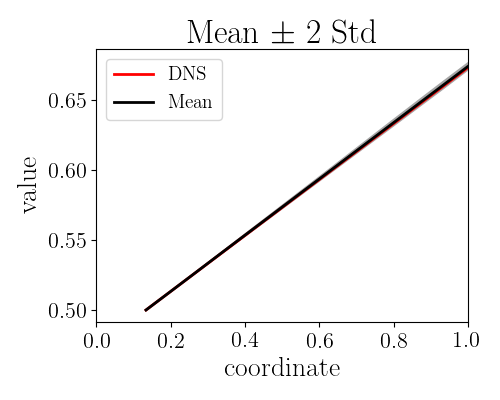}} \\ 
  {\includegraphics[height=0.08\linewidth]{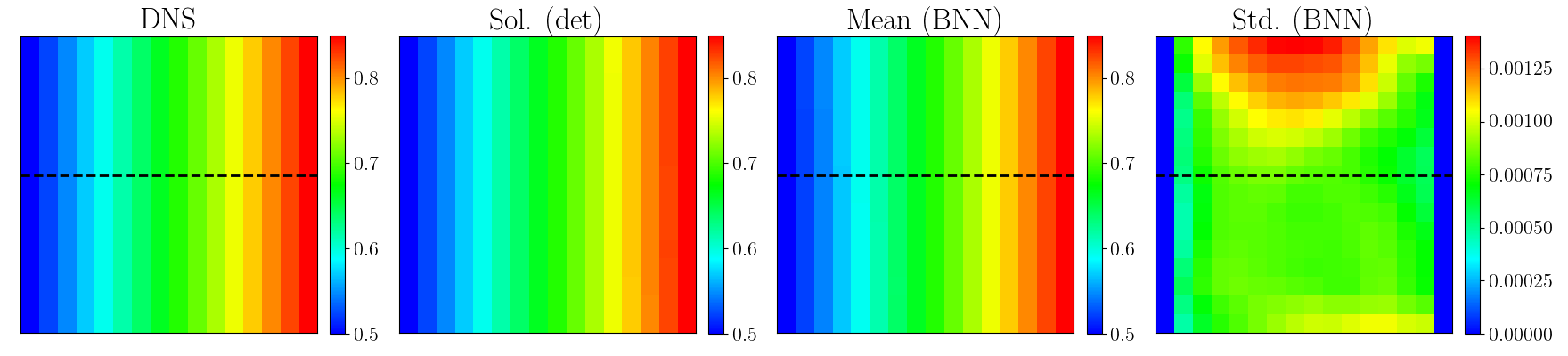}} 
  {\includegraphics[height=0.08\linewidth]{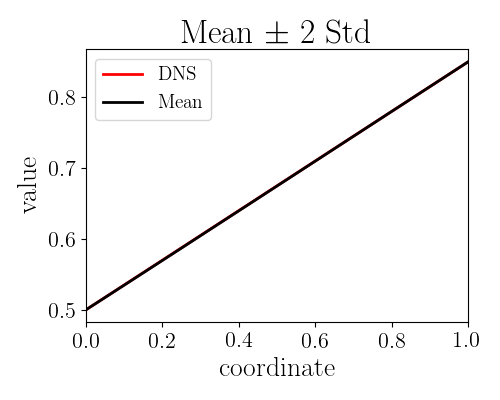}} 
  {\includegraphics[height=0.08\linewidth]{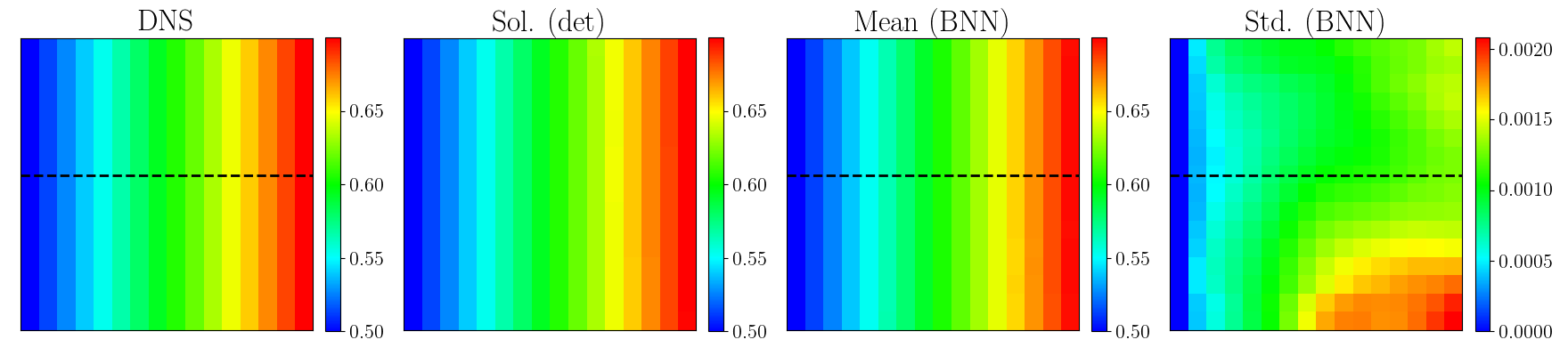}} 
  {\includegraphics[height=0.08\linewidth]{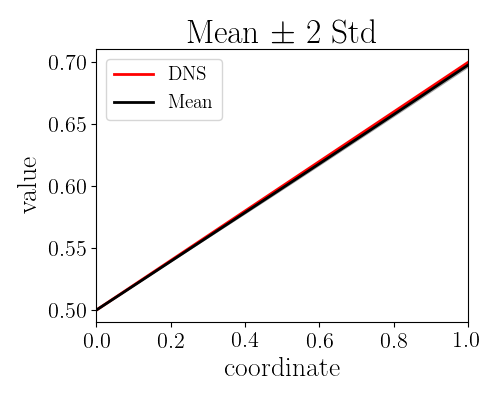}} \\
  {\includegraphics[height=0.08\linewidth]{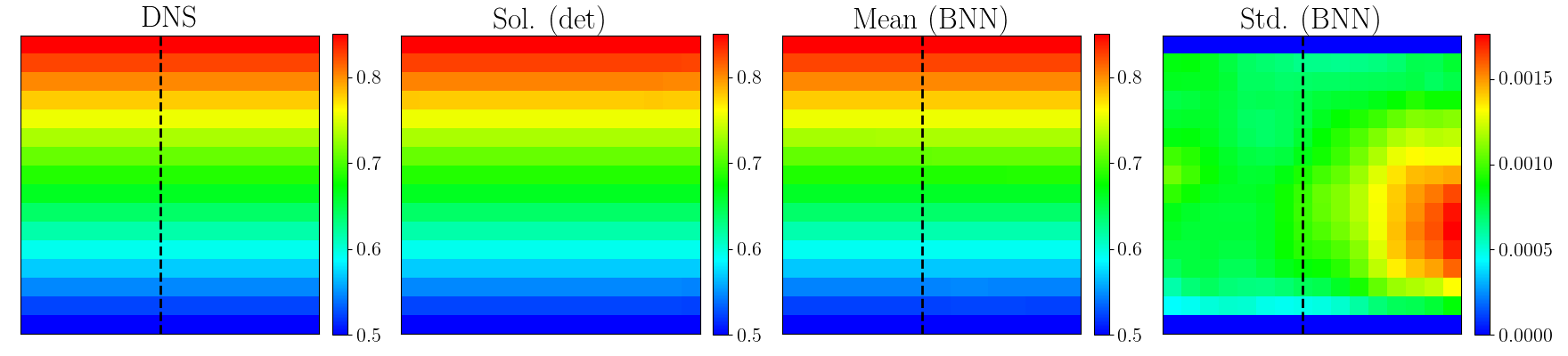}} 
  {\includegraphics[height=0.08\linewidth]{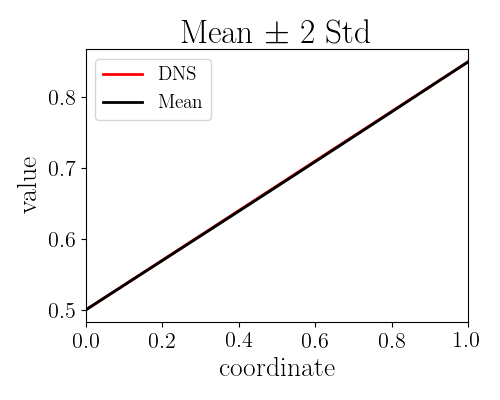}} 
  {\includegraphics[height=0.08\linewidth]{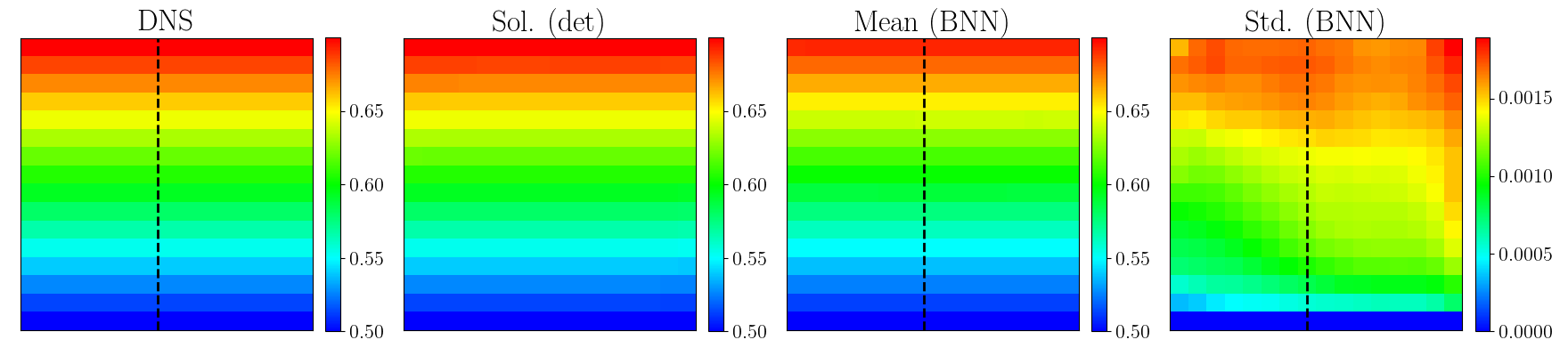}} 
  {\includegraphics[height=0.08\linewidth]{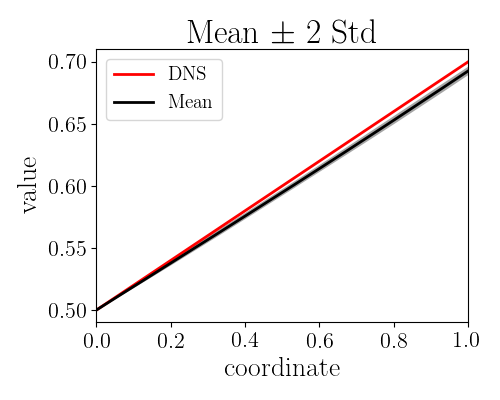}} \\
  {\includegraphics[height=0.08\linewidth]{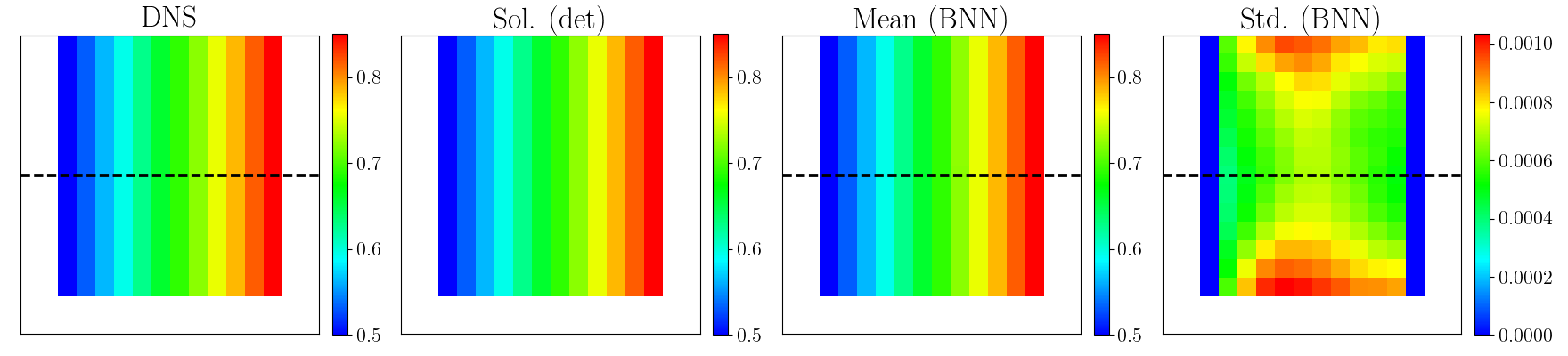}} 
  {\includegraphics[height=0.08\linewidth]{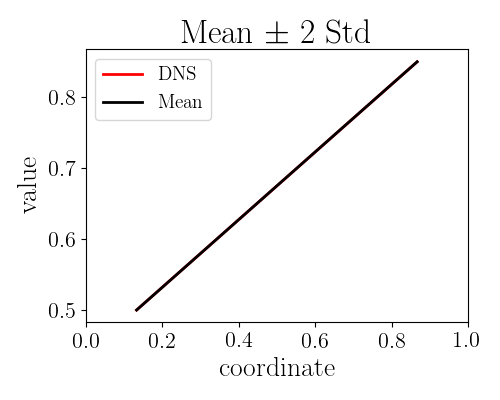}} 
  {\includegraphics[height=0.08\linewidth]{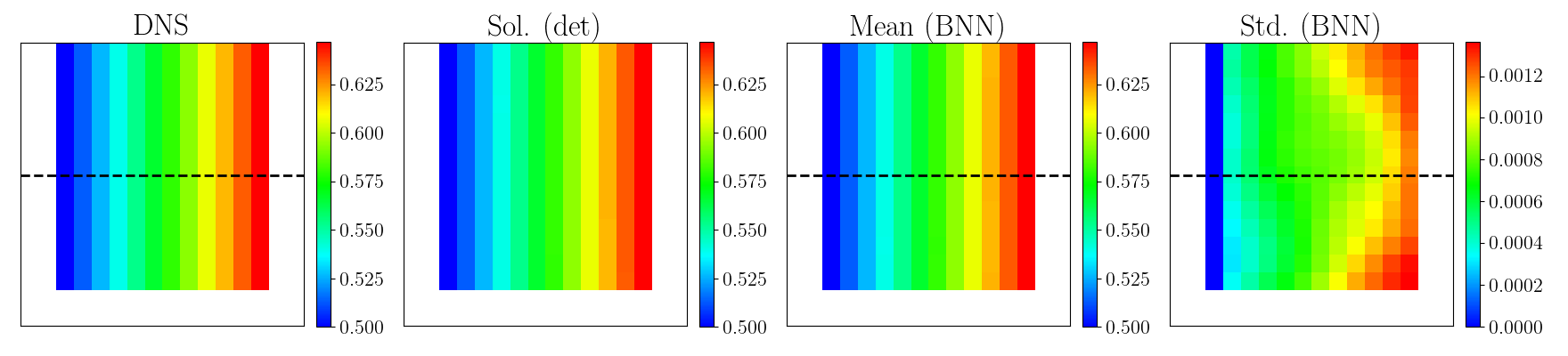}} 
  {\includegraphics[height=0.08\linewidth]{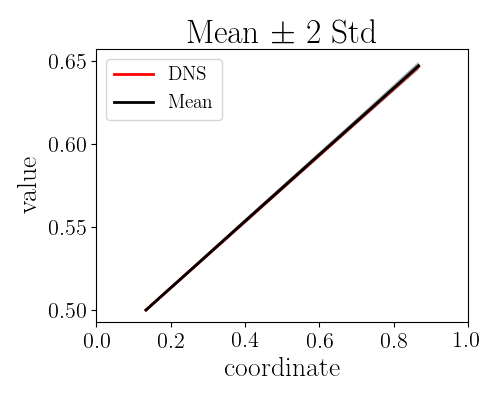}} \\
  {\includegraphics[height=0.08\linewidth]{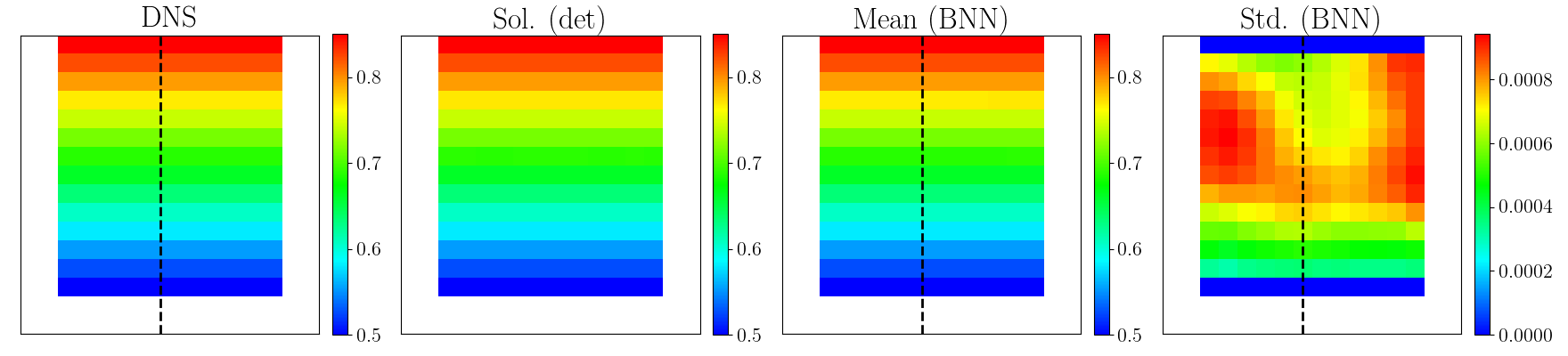}} 
  {\includegraphics[height=0.08\linewidth]{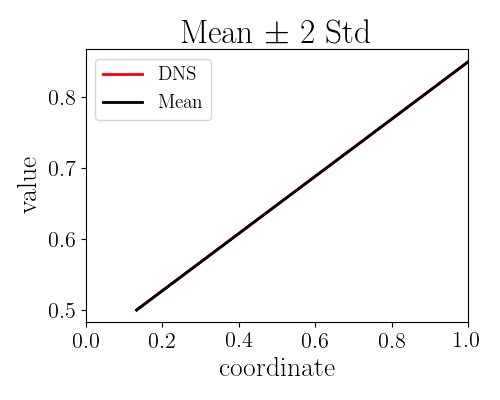}} 
  {\includegraphics[height=0.08\linewidth]{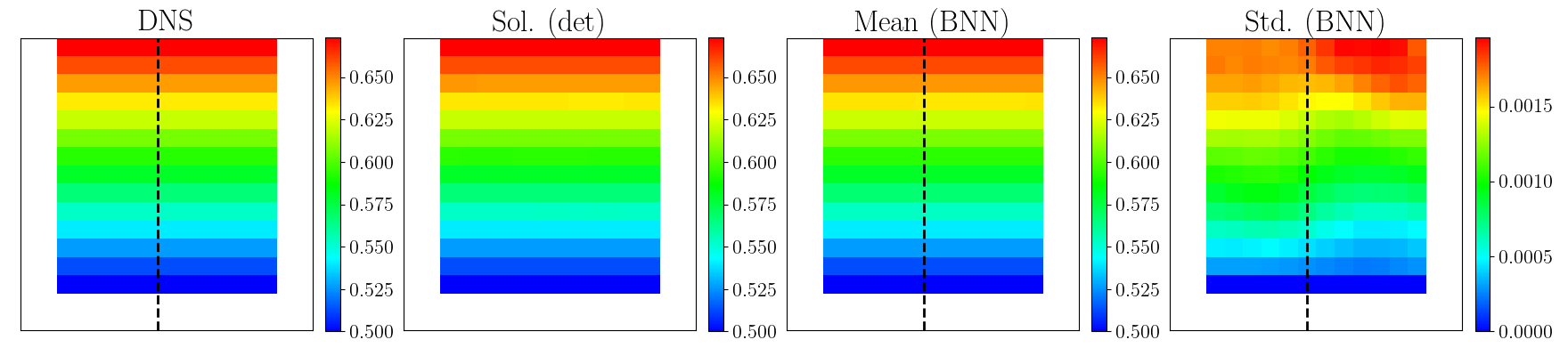}} 
  {\includegraphics[height=0.08\linewidth]{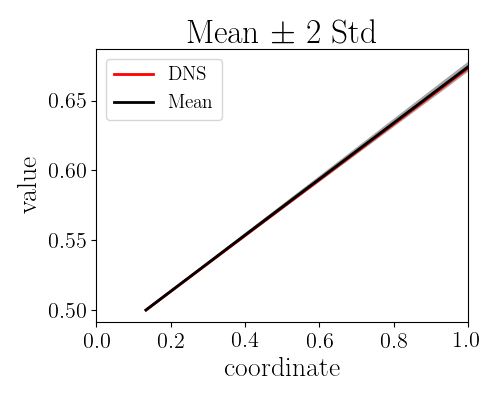}} \\
  {\includegraphics[height=0.08\linewidth]{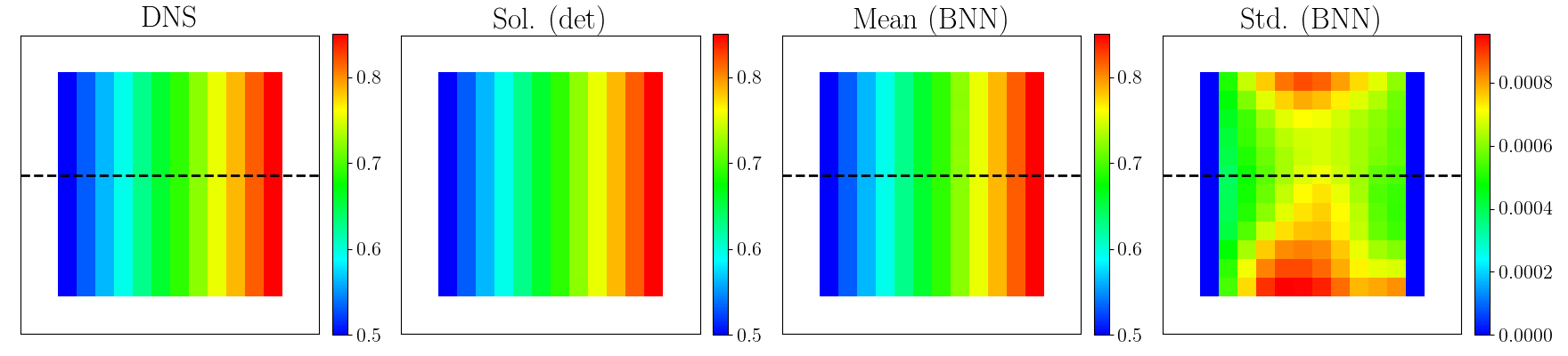}} 
  {\includegraphics[height=0.08\linewidth]{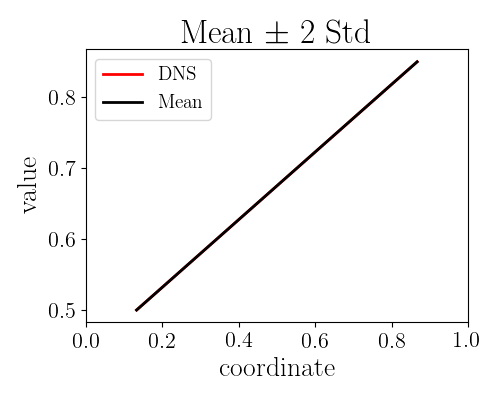}} 
  {\includegraphics[height=0.08\linewidth]{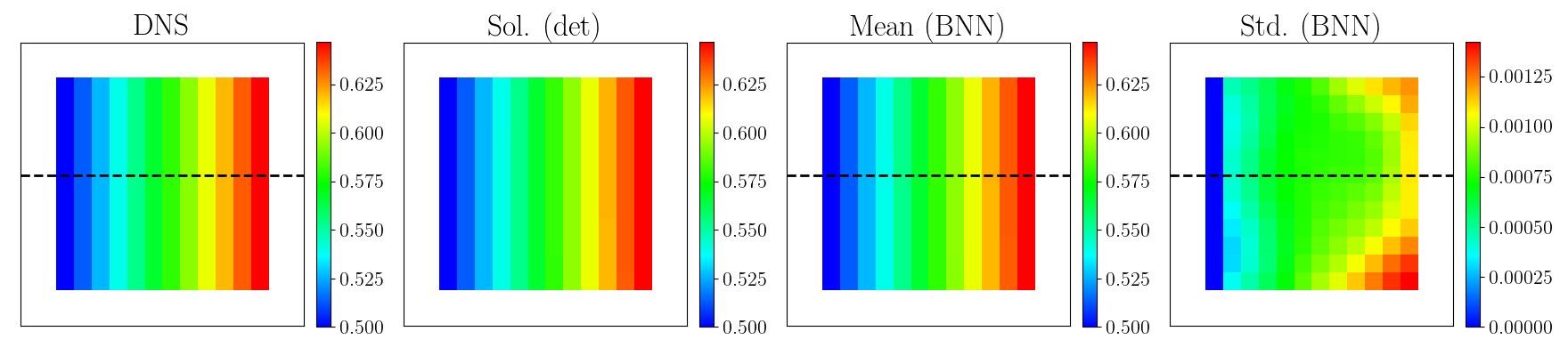}} 
  {\includegraphics[height=0.08\linewidth]{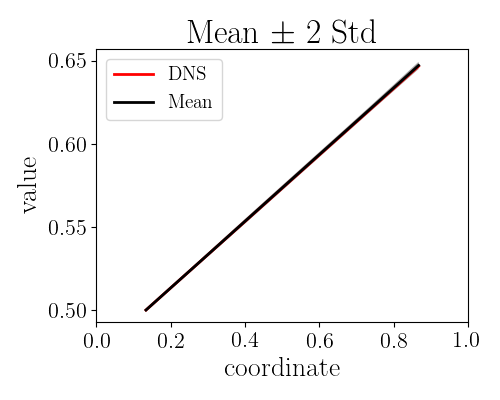}} \\ 
  {\includegraphics[height=0.08\linewidth]{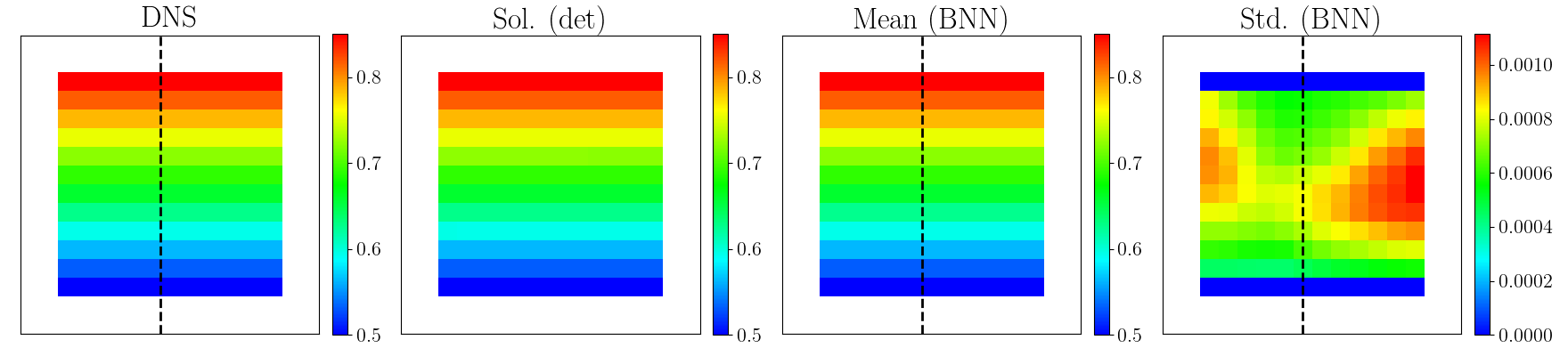}} 
  {\includegraphics[height=0.08\linewidth]{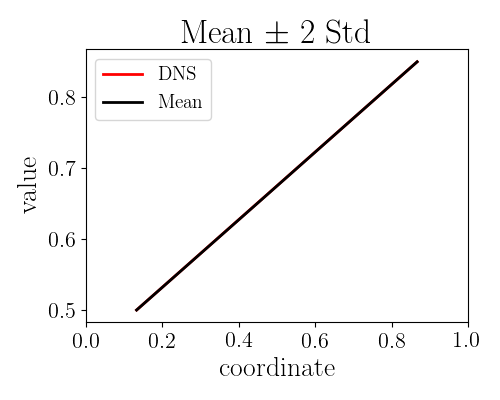}} 
  {\includegraphics[height=0.08\linewidth]{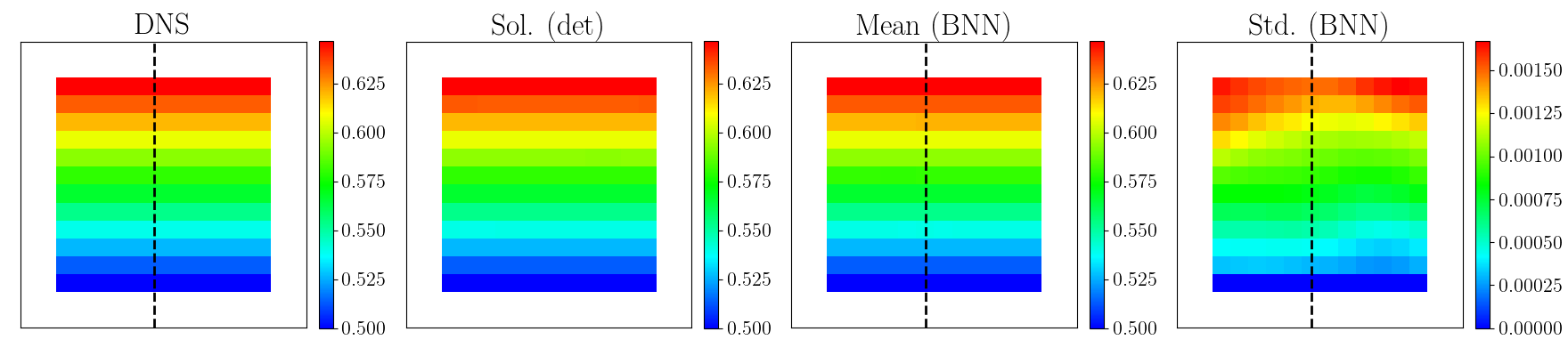}} 
  {\includegraphics[height=0.08\linewidth]{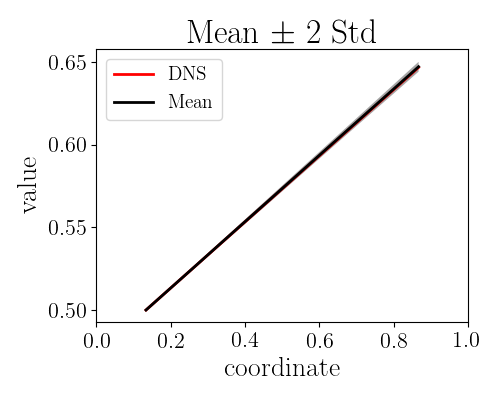}} \\
  {\includegraphics[height=0.08\linewidth]{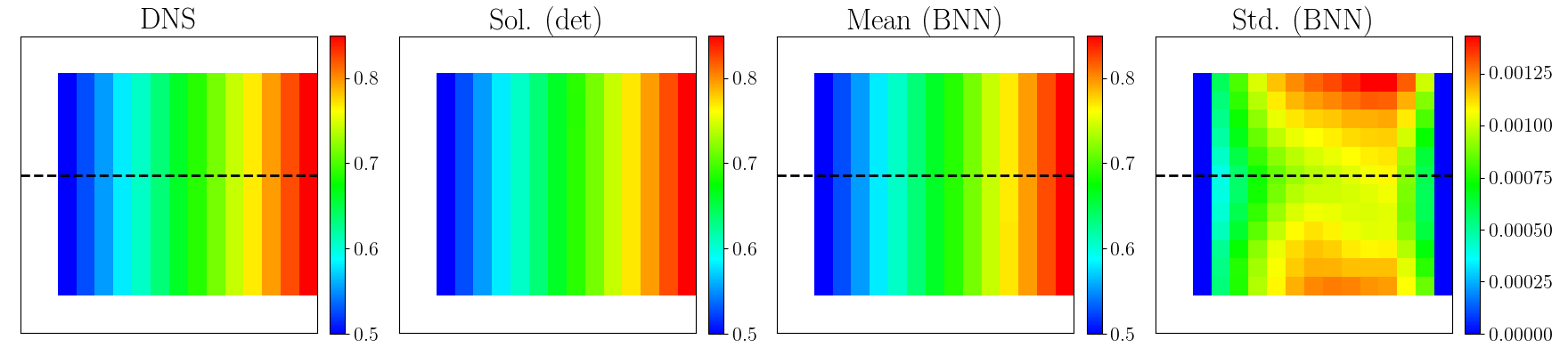}} 
  {\includegraphics[height=0.08\linewidth]{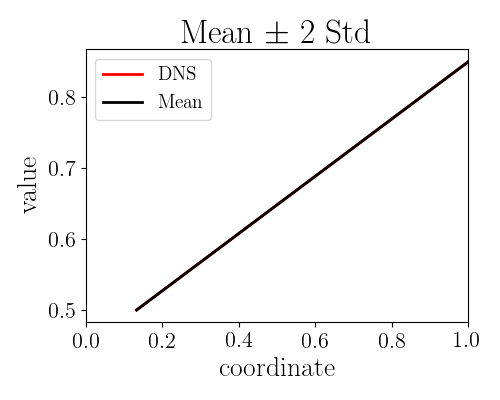}} 
  {\includegraphics[height=0.08\linewidth]{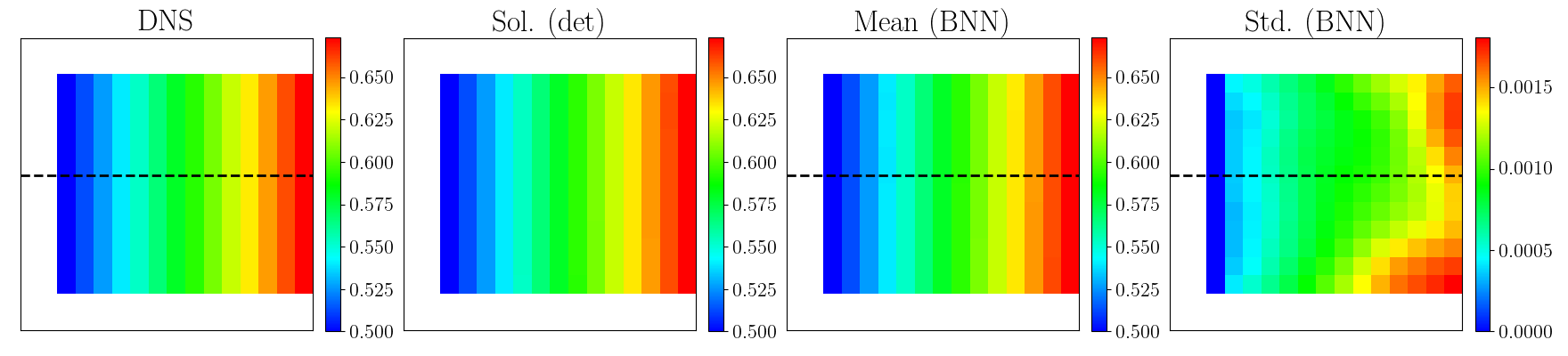}} 
  {\includegraphics[height=0.08\linewidth]{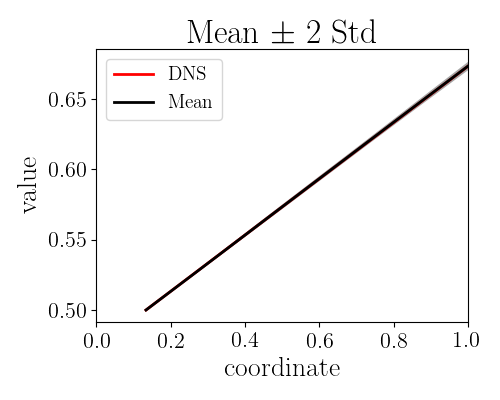}} \\
  {\includegraphics[height=0.08\linewidth]{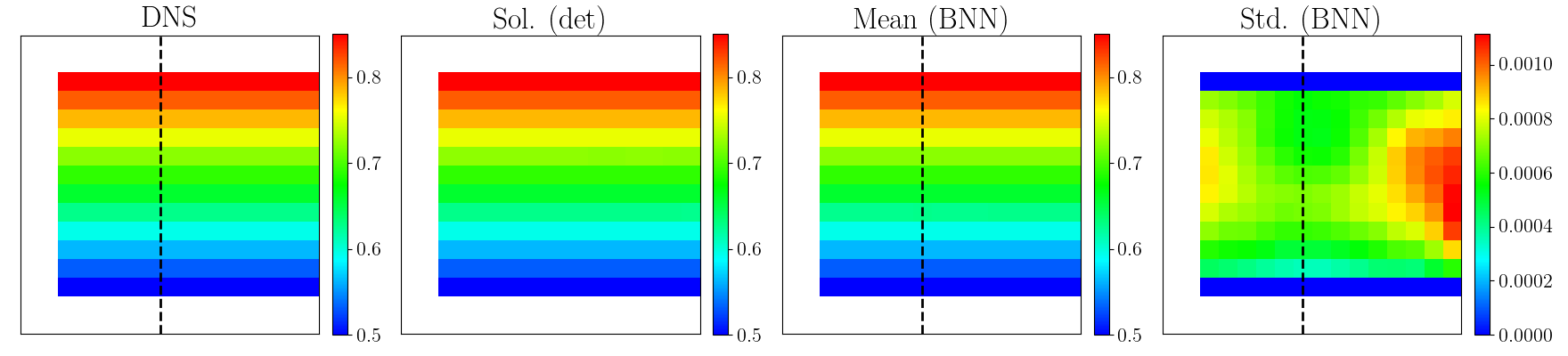}} 
  {\includegraphics[height=0.08\linewidth]{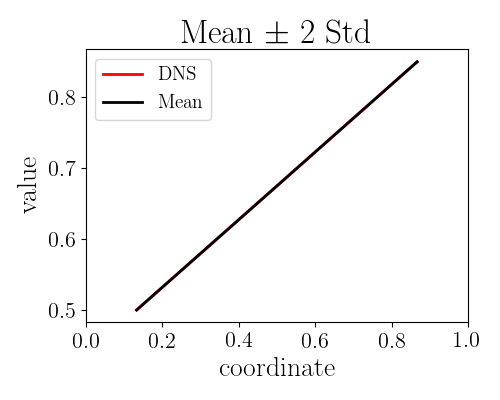}} 
  {\includegraphics[height=0.08\linewidth]{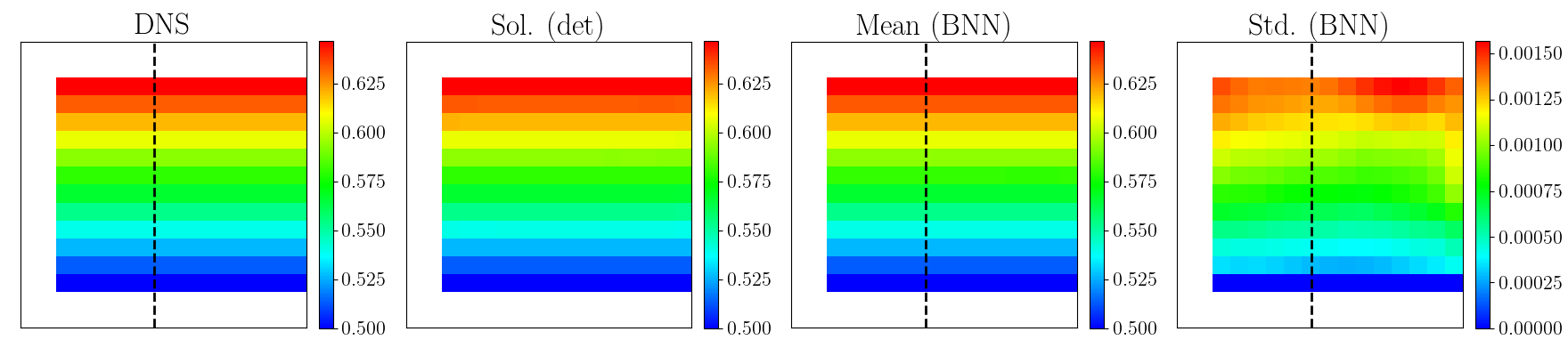}} 
  {\includegraphics[height=0.08\linewidth]{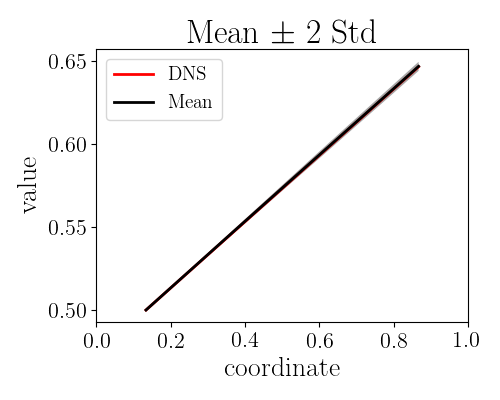}} \\
  \caption{NN results for steady-state diffusion BVPs on rectangle domains. Comparison between the FEM solution (DNS), the deterministic NN solution, the mean and std of BNN solutions over 50 MC samplings, and solution distributions along the dashed lines. }
  \label{fig:diffusion-20bvp-results}
\end{figure}

We also use the proposed PDE constrained NNs to simultaneously solve 20 steady-state diffusion BVPs with a resolution of $16\times16$.
The architectures of both deterministic and probabilistic NNs and other training related NN parameters are summarized in Table \ref{tab:diffusion-20bvp-NNs} and \ref{tab:diffusion-20bvp-NNs-others}, respectively.
The NN hyperparameters are manually tuned to achieve a desired performance.
We follow the training procedures in Algorithm Boxes \ref{algo:cnn-training} and \ref{algo:bnn-training} to first train the deterministic NN with zero initialization, followed by training the BNNs with the optimal parameter initialization process.
The results are shown in Fig. \ref{fig:diffusion-20bvp-results}, which confirm the accuracy of the proposed method an d demonstrate that the proposed method can simultaneously solve multiple BVPs.  
The statistical moments of the BNN predictions are evaluated based on 50 MC samplings.

\subsection{Linear elasticity - small dataset}
\subsubsection{Multiple rectangular domains with different BCs}\label{sec:linear-rectangle}

\begin{figure}[t!]
  \centering
  \includegraphics[width=0.65\linewidth]{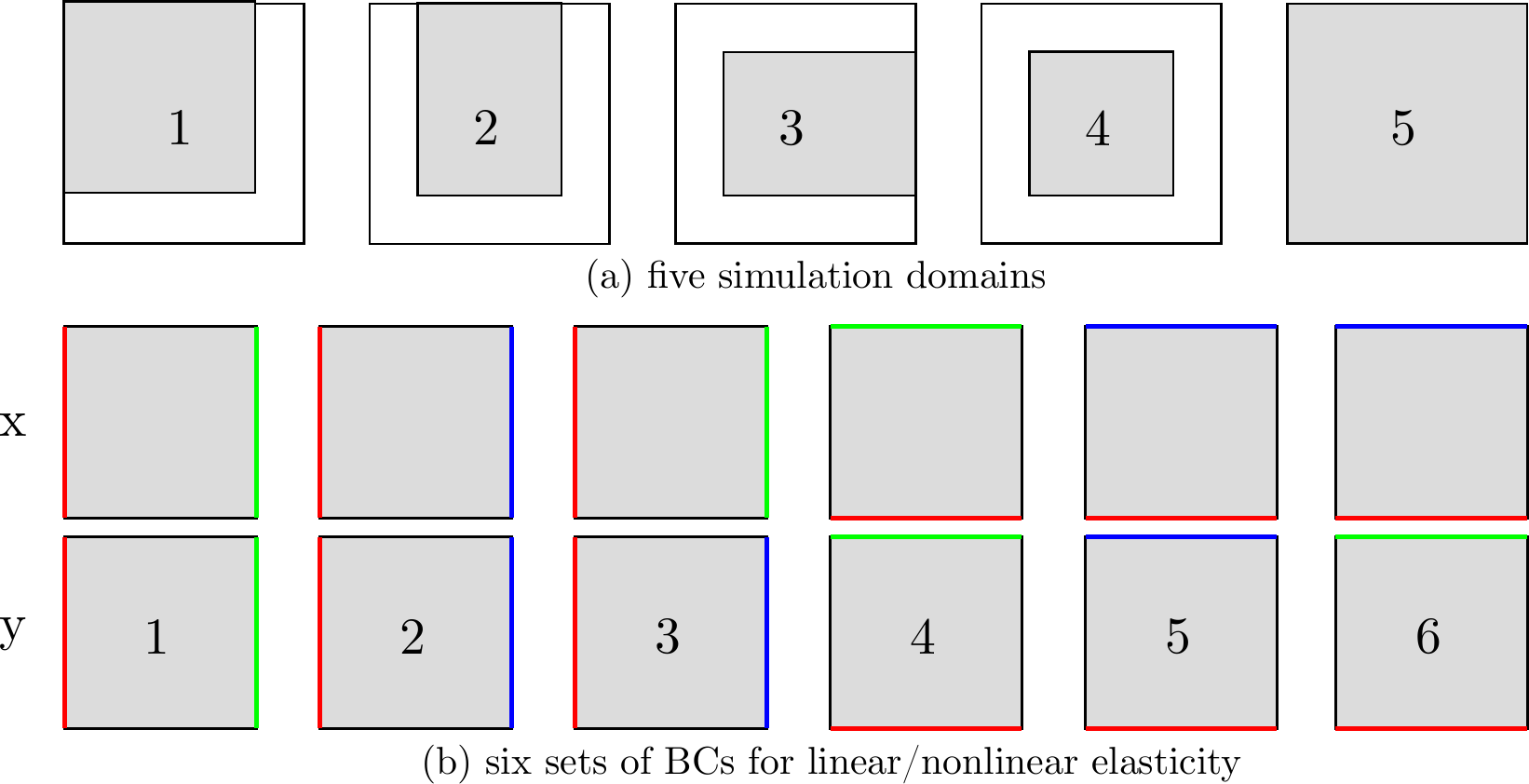}
  \caption{Illustration of the setup of 30 BVPs on different domains for linear/nonlinear elasticity problems. In these drawings, red represents a zero Dirichlet BC. Green represents a non-zero Dirichlet BC. Blue represents a non-zero Neumann BC. No color is assigned to Zero Neumann BCs. (a) Setup of five rectangle simulation domains of different sizes and locations on a fixed background grid with different applied BCs. (b) For linear/nonlinear elasticity, 6 sets of BCs are assigned to each simulation domain, leading to 30 linear/nonlinear elasticity BVPs.}
  \label{fig:bvp-problem-setup}
\end{figure}

\begin{figure}[t!]
    \centering
    \subfloat[three selected BVPs on rectangle domains]{\includegraphics[height=0.20\linewidth]{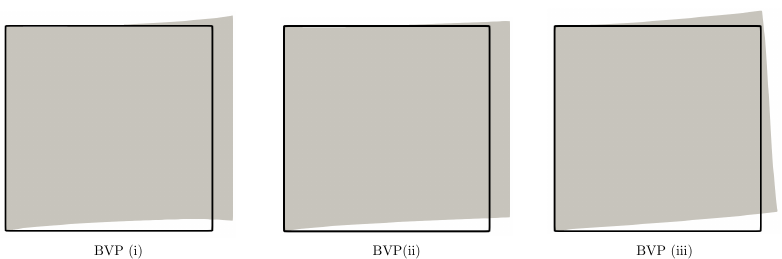}} 
    \caption{Illustration of the deformed shape of selected linear elasticity BVPs. The wireframe and the gray region indicate the undeformed and deformed problem domain, respectively. Three selected BVPs of out the 30 BVPs solved in section \ref{sec:linear-rectangle}. }
    \label{fig:linear-bvp-deformed}
\end{figure}

\begin{figure}[p!]
  \centering
  \subfloat[$u_x$ results for BVP (i)]  {\includegraphics[trim=0.1cm 0.5cm 0.1cm 0.25cm, clip,height=0.16\linewidth]{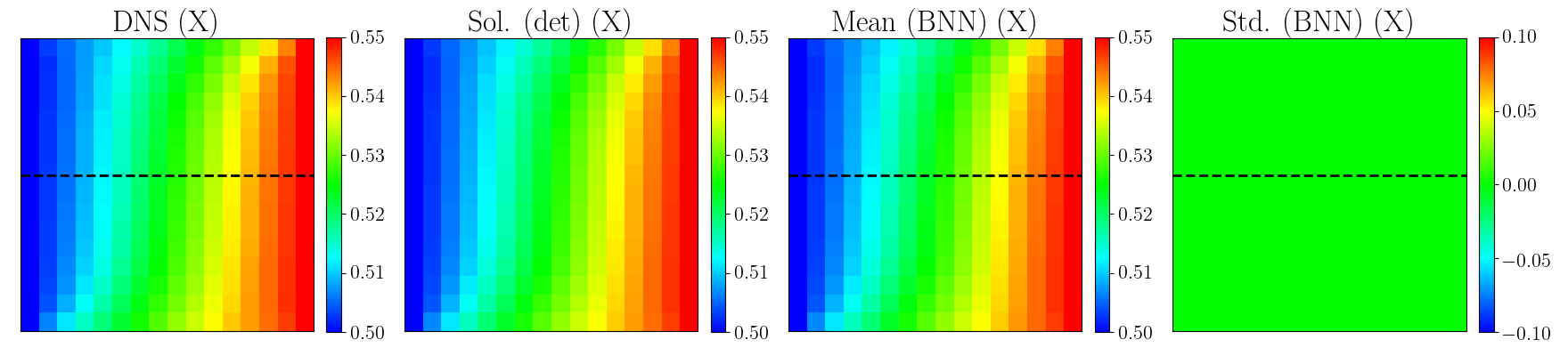}} 
  \subfloat[$u_x$ UQ for BVP (i)]       {\includegraphics[trim=0.1cm 0.5cm 0.1cm 0.25cm, clip,height=0.16\linewidth]{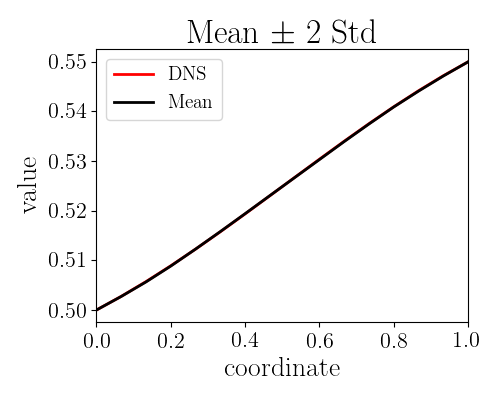}} \\[-2mm]
  \subfloat[$u_y$ results for BVP (i)]  {\includegraphics[trim=0.1cm 0.5cm 0.1cm 0.25cm, clip,height=0.16\linewidth]{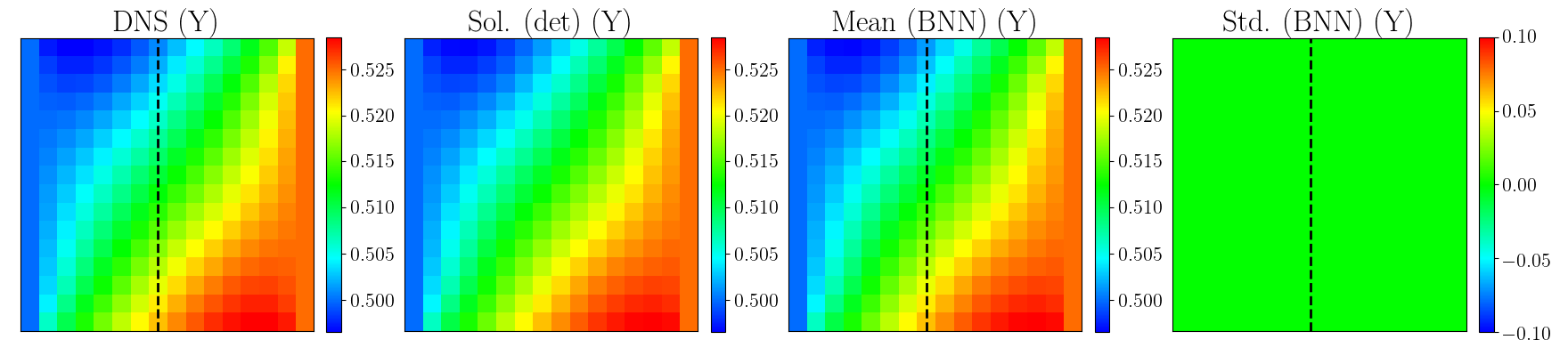}} 
  \subfloat[$u_y$ UQ for BVP (i)]       {\includegraphics[trim=0.1cm 0.5cm 0.1cm 0.25cm, clip,height=0.16\linewidth]{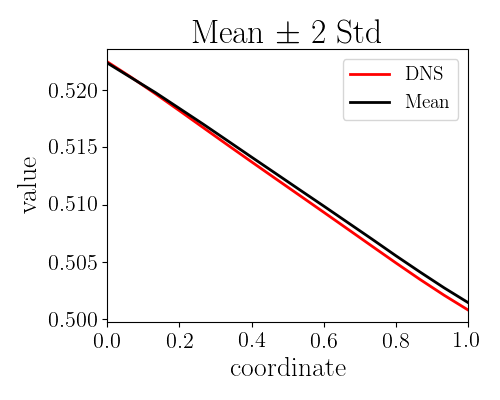}} \\[-2mm]
  \subfloat[$u_x$ results for BVP (ii)] {\includegraphics[trim=0.1cm 0.5cm 0.1cm 0.25cm, clip,height=0.16\linewidth]{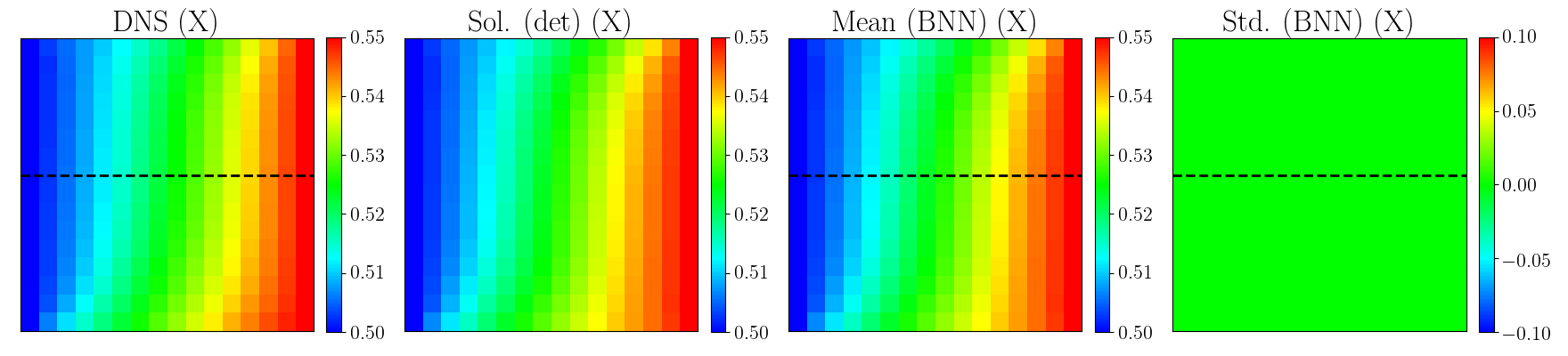}} 
  \subfloat[$u_x$ UQ for BVP (ii)]      {\includegraphics[trim=0.1cm 0.5cm 0.1cm 0.25cm, clip,height=0.16\linewidth]{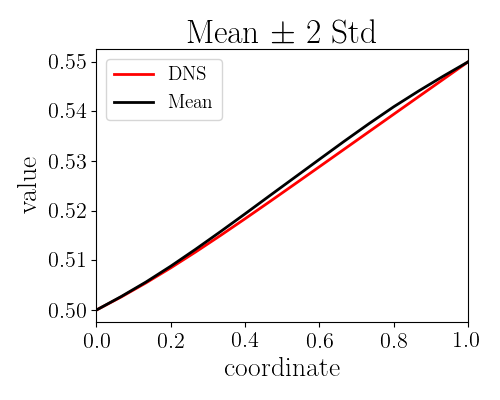}} \\[-2mm]
  \subfloat[$u_y$ results for BVP (ii)] {\includegraphics[trim=0.1cm 0.5cm 0.1cm 0.25cm, clip,height=0.16\linewidth]{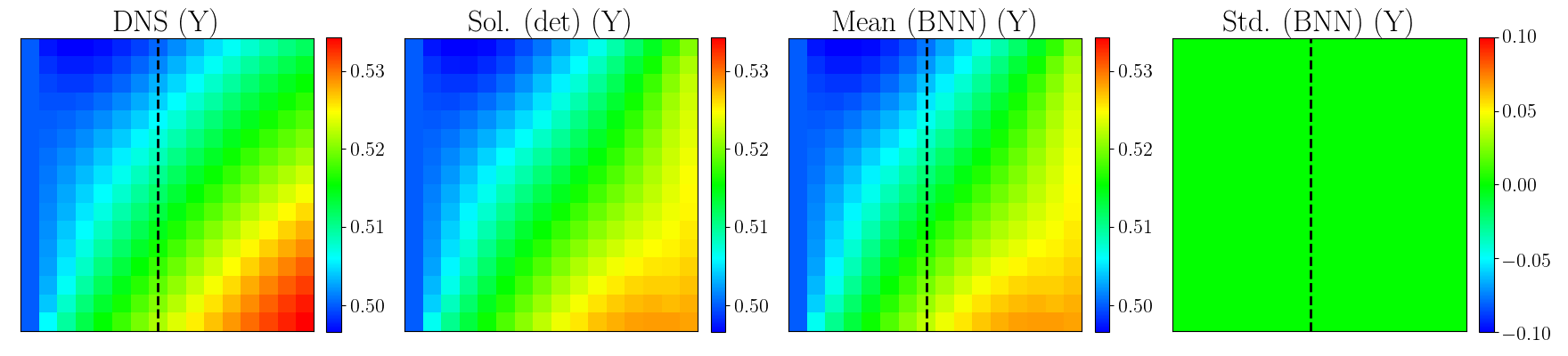}} 
  \subfloat[$u_y$ UQ for BVP (ii)]      {\includegraphics[trim=0.1cm 0.5cm 0.1cm 0.25cm, clip,height=0.16\linewidth]{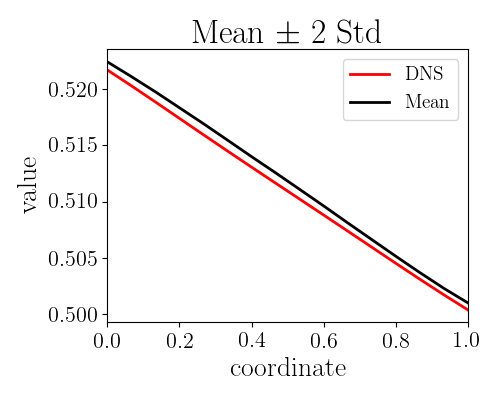}} \\[-2mm]
  \subfloat[$u_x$ results for BVP (iii)]{\includegraphics[trim=0.1cm 0.5cm 0.1cm 0.25cm, clip,height=0.16\linewidth]{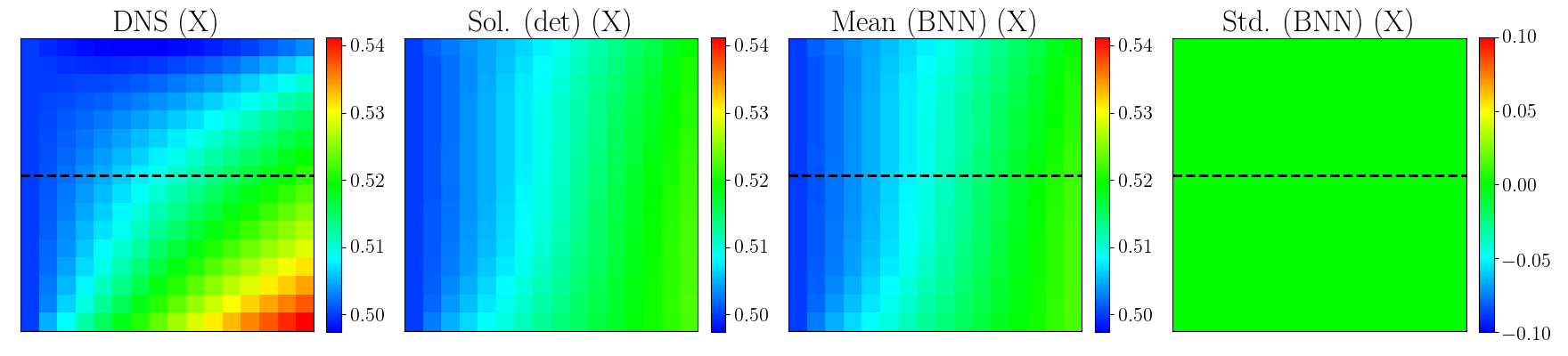}} 
  \subfloat[$u_x$ UQ for BVP (iii)]     {\includegraphics[trim=0.1cm 0.5cm 0.1cm 0.25cm, clip,height=0.16\linewidth]{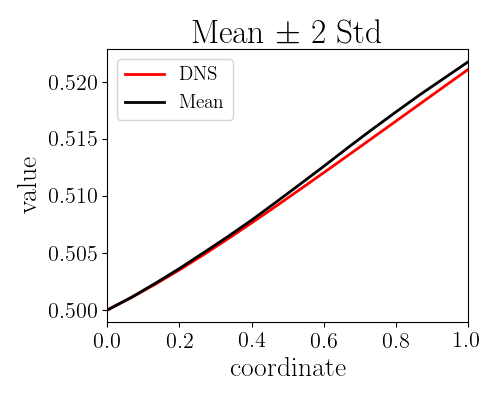}} \\[-2mm]
  \subfloat[$u_y$ results for BVP (iii)]{\includegraphics[trim=0.1cm 0.5cm 0.1cm 0.25cm, clip,height=0.16\linewidth]{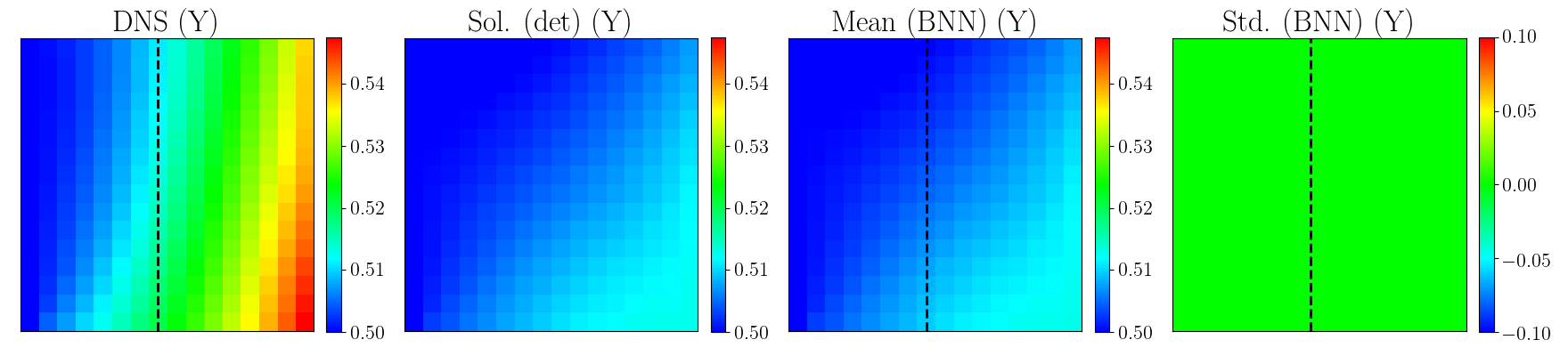}} 
  \subfloat[$u_y$ UQ for BVP (iii)]     {\includegraphics[trim=0.1cm 0.5cm 0.1cm 0.25cm, clip,height=0.16\linewidth]{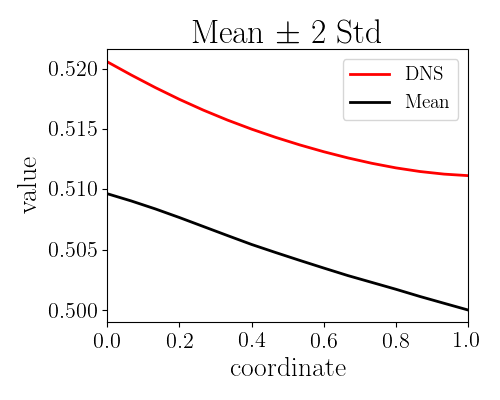}} \\
  \caption{Results of three selected BVPs out of the 30 linear elasticity BVPs with varying domains and different applied BCs simultaneously solved by a single deterministic or probabilistic NN with the proposed method. BVP (i), (ii), (iii) correspond to bc id 1, 2, and 3 for domain id 5, as shown in Fig. \ref{fig:bvp-problem-setup}(a). 
      (a, c, e, g, i, k) Solutions from DNS, deterministic (det) NNs, and BNNs (Mean, Std.) for different BVPs.
    (b, d, f, h, j, l) Quantitative comparison of the solution distribution between DNS and BNNs along the dashed lines.
  }
  \label{fig:linear-30bvp-results}
\end{figure}

\begin{table}
  \centering
  \begin{tabular}{l | l | l }
    \hline
    Description                     & Deterministic                 & Probabilistic         \\ \hline
    Total parameters                & 34,010                        & 67,803                 \\
    Size of $\calD$                 & 30 $\times$ Aug: $2^{9}$      & 30 $\times$ Aug: $2^{9}$      \\
    Epochs                          & 20,000                        & 100                   \\
    Zero initialization epochs      & 100                           & -                     \\
    Optimizer                       & Nadam                         & Nadam                 \\
    Learning Rate                   & 2.5e-4                        & 1e-8                  \\
    Batch Size                      & 128                           & 64                    \\
    $\Sigma_1$                      & -                             & 1e-8                  \\
    Initial value of $\Sigma_2$     & -                             & 1e-8                  \\
    \hline
  \end{tabular}
  \caption{Training related parameters for solving 30 linear elasticity BVPs. Aug: data augmentation. }
  \label{tab:linear-30bvp-NNs-others}
\end{table}

In this section, we use the proposed PDE constrained NNs to simultaneously solve 30 linear elasticity BVPs, as shown in Fig. \ref{fig:bvp-problem-setup}(a), with a resolution of $16\times16$. 
The deformed problem domains from DNSs for three representative BVPs are shown in Fig. \ref{fig:linear-bvp-deformed}(a).
Similar architectures of both deterministic and probabilistic NNs as summarized in Table \ref{tab:diffusion-20bvp-NNs} are used, except that.  the last layer has two filters, representing $u_x$ and $u_y$, instead of one for the steady-state diffusion problem.
The other training related NN parameters are summarized in Table \ref{tab:linear-30bvp-NNs-others}. 
We follow the procedures described in Section \ref{sec:NN-training} to train both types of NNs. 
The NN results of three selected BVPs, as shown in Fig. \ref{fig:linear-bvp-deformed}(a),  are presented in Fig. \ref{fig:linear-30bvp-results}
The statistical moments of the BNN predictions are evaluated based on 50 MC samplings.
In Fig. \ref{fig:linear-30bvp-results}, BVP (i), (ii), and (iii) correspond to bc id 1 (non-zero Dirichlet loading), bc id 2 (non-zero Neumann loading), bc id 3 (mixed loading) applied to domain id 5, respectively. 
The comparison of solutions between DNSs, the deterministic NN, and the BNN for these threeBVPs is shown qualitatively in Fig. \ref{fig:linear-30bvp-results}(a,c,e,g,i,k), with quantitative comparison of the solution distribution along the dashed lines between DNSs and the BNN given in Fig. \ref{fig:linear-30bvp-results}(b,d,f,h,j,l).
Such a comparison shows that the proposed method has successfully solved most of the BVPs with desired accuracy.
The results from NNs in Fig. \ref{fig:linear-30bvp-results}(i,k,l) are slightly worse than the DNSs. 
This happens mainly because the deformation for linear elasticity is small. 
The scaled results have a narrow range of $[0.5, 0.55]$, which is challenging for NNs to learn to distinguish, particularly for purely non-zero traction loads.
For the mathematically more complex nonlinear elasticity BVPs, in which the deformation is large, are solved by the NNs more successfully, as shown in Fig. \ref{fig:nonlinear-30bvp-results}.

\subsubsection{L-shape domain with solution interpolation}\label{sec:linear-lshape}
NN architectures for the L-shaped domain simulation are summarized in Table \ref{tab:linear-lshape-NNs} and \ref{tab:linear-lshape-NNs-others}.
We follow the procedures described in Section \ref{sec:NN-training} to train both types of NNs with an output resolution of $32\times32$. 

\begin{table}
  \centering
  \begin{tabular}{l | l | l | l}
    \hline
    Deterministic         & Probabilistic         & Size         & Layer arguments \\ \hline
    Input                 & Input                 & -            & - \\
    LayerFillRandomNumber & LayerFillRandomNumber & -            & - \\
    Conv2D                & Convolution2DFlipout  & filters = 8  & kernel (5,5), padding: same, ReLU \\
    MaxPooling2D          & MaxPooling2D          & -            & kernel (2,2), padding: same\\
    Conv2D                & Convolution2DFlipout  & filters = 8 & kernel (5,5), padding: same, ReLU \\
    MaxPooling2D          & MaxPooling2D          & -            & kernel (2,2), padding: same\\
    Conv2D                & Convolution2DFlipout  & filters = 16 & kernel (5,5), padding: same, ReLU \\
    MaxPooling2D          & MaxPooling2D          & -            & kernel (2,2), padding: same\\
    Flatten               & Flatten               & -            & - \\
    Dense                 & DenseFlipout          & units = 32   & ReLU \\
    Dense                 & DenseFlipout          & units = 128   & ReLU \\
    Reshape               & Reshape               & -            & $[4,4,8]$ \\
    Conv2D                & Convolution2DFlipout  & filters = 16 & kernel (5,5), padding: same, ReLU \\
    UpSampling2D          & UpSampling2D          & -            & size (2,2) \\
    Conv2D                & Convolution2DFlipout  & filters = 16 & kernel (5,5), padding: same, ReLU \\
    UpSampling2D          & UpSampling2D          & -            & size (2,2) \\
    Conv2D                & Convolution2DFlipout  & filters = 16 & kernel (5,5), padding: same, ReLU \\
    UpSampling2D          & UpSampling2D          & -            & size (2,2) \\
    Conv2D                & Convolution2DFlipout  & filters = 16 & kernel (5,5), padding: same, ReLU \\
    Conv2D                & Convolution2DFlipout  & filters = 2  & kernel (5,5), padding: same, ReLU \\
    \hline
  \end{tabular}
  \caption{Details of both deterministic and probabilistic NNs for solving linear elasticity L-shape BVPs.}
  \label{tab:linear-lshape-NNs}
\end{table}

\begin{table}
  \centering
  \begin{tabular}{l | l | l }
    \hline
    Description                     & Deterministic                 & Probabilistic         \\ \hline
    Total parameters                & 41,346                        & 82,435                 \\
    Size of $\calD$                 & 5 $\times$ Aug: $2^{10}$      & 5 $\times$ Aug: $2^{9}$      \\
    Epochs                          & 10,000                        & 100                 \\
    Zero initialization epochs      & 100                           & -                     \\
    Optimizer                       & Adam                          & Nadam                 \\
    Learning Rate                   & 2.5e-4                        & 1e-8                  \\
    Batch Size                      & 256                           & 64                    \\
    $\Sigma_1$                      & -                             & 1e-8                  \\
    Initial value of $\Sigma_2$     & -                             & 1e-8                  \\
    \hline
  \end{tabular}
  \caption{Training related parameters for solving linear elasticity on L-shaped BVPs. Aug: data augmentation.}
  \label{tab:linear-lshape-NNs-others}
\end{table}

\subsection{Nonlinear elasticity - small dataset}
Even with the zero-initialization process, the NN outputs at early stages of training  could violate the physics, e.g. with a negative or zero determinant of the deformation gradient $J$. To ensure that the residual can be evaluated and to prevent residuals from these ``bad'' pixels values from contributing to the final loss, we regularize the loss by omitting the residual contribution with $J<0.1$ and $J>5.0$. As the training continues towards a later stage, the NN predicted solutions gradually fulfill the governing PDEs, and the regularization on $J$ will cease to function.

\subsubsection{Multiple rectangular domains with different BCs}

\begin{figure}[t!]
    \centering
    \includegraphics[width=0.75\linewidth]{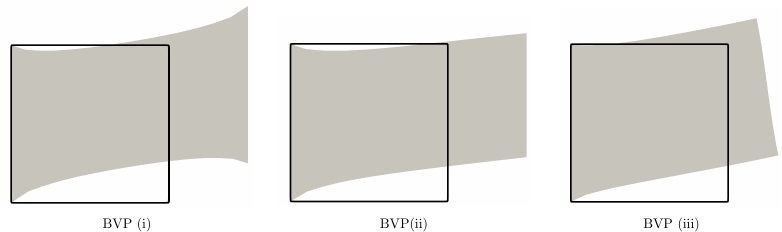}
    \caption{Illustration of the deformed shape of the three selected nonlinear elasticity BVPs. The wireframe and the gray region indicate the undeformed and deformed problem domains, respectively.}
    \label{fig:nonlinear-30bvp-deformed}
\end{figure}

\begin{figure}[h!]
  \centering
  \subfloat[$u_x$ results for BVP (i)]  {\includegraphics[trim=0.1cm 0.5cm 0.1cm 0.25cm, clip,height=0.15\linewidth]{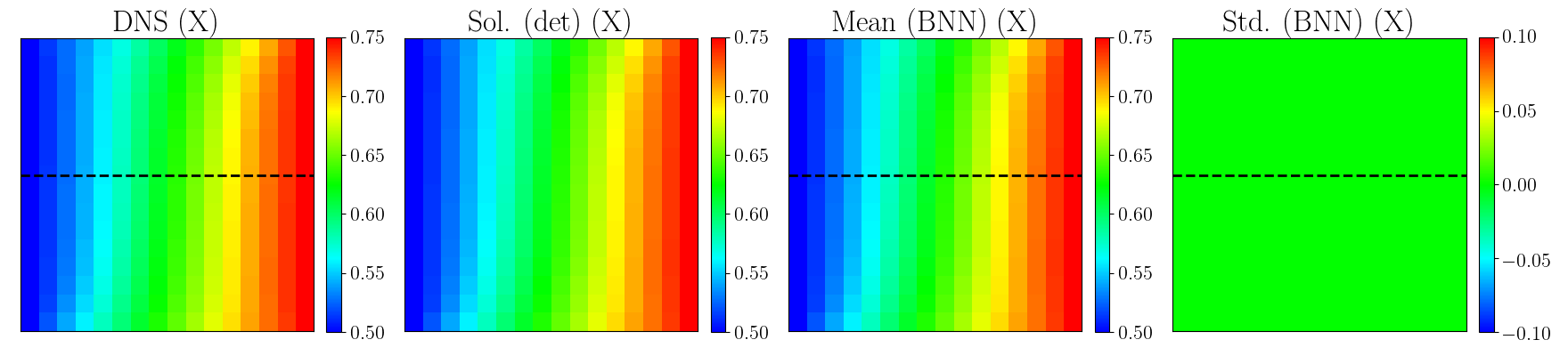}} 
  \subfloat[$u_x$ UQ for BVP (i)]       {\includegraphics[trim=0.1cm 0.5cm 0.1cm 0.25cm, clip,height=0.15\linewidth]{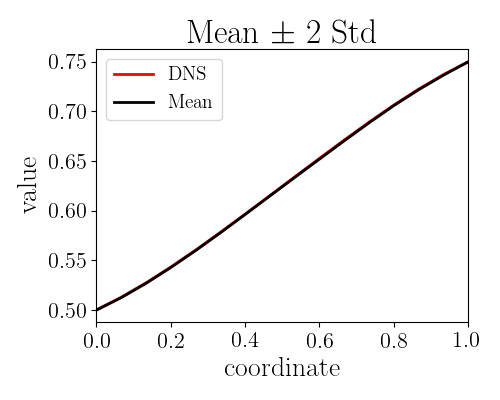}} \\[-2mm]
  \subfloat[$u_y$ results for BVP (i)]  {\includegraphics[trim=0.1cm 0.5cm 0.1cm 0.25cm, clip,height=0.15\linewidth]{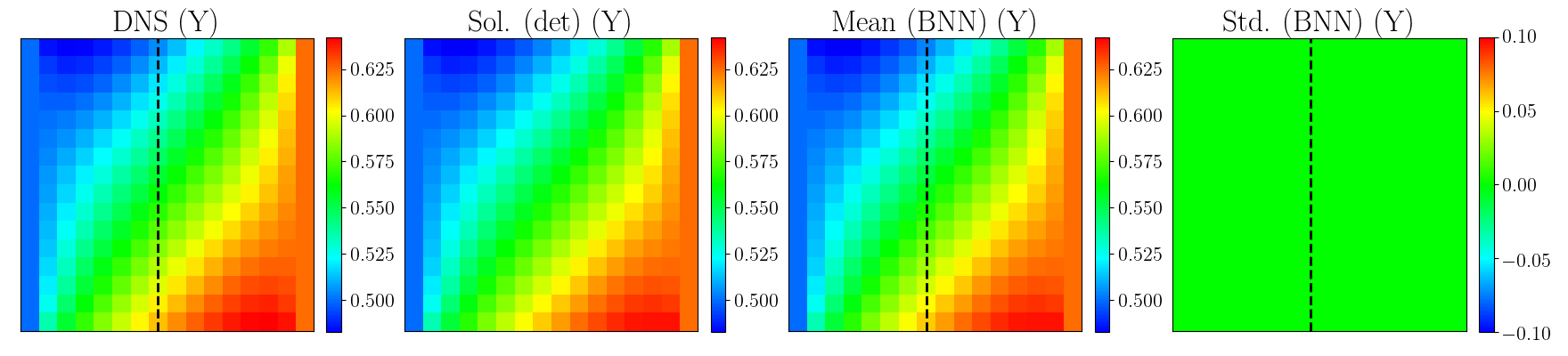}} 
  \subfloat[$u_y$ UQ for BVP (i)]       {\includegraphics[trim=0.1cm 0.5cm 0.1cm 0.25cm, clip,height=0.15\linewidth]{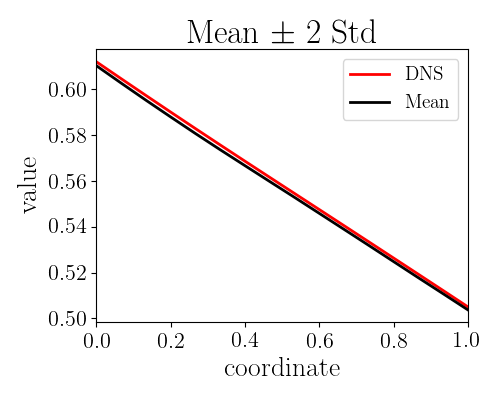}} \\[-2mm]
  \subfloat[$u_x$ results for BVP (ii)] {\includegraphics[trim=0.1cm 0.5cm 0.1cm 0.25cm, clip,height=0.15\linewidth]{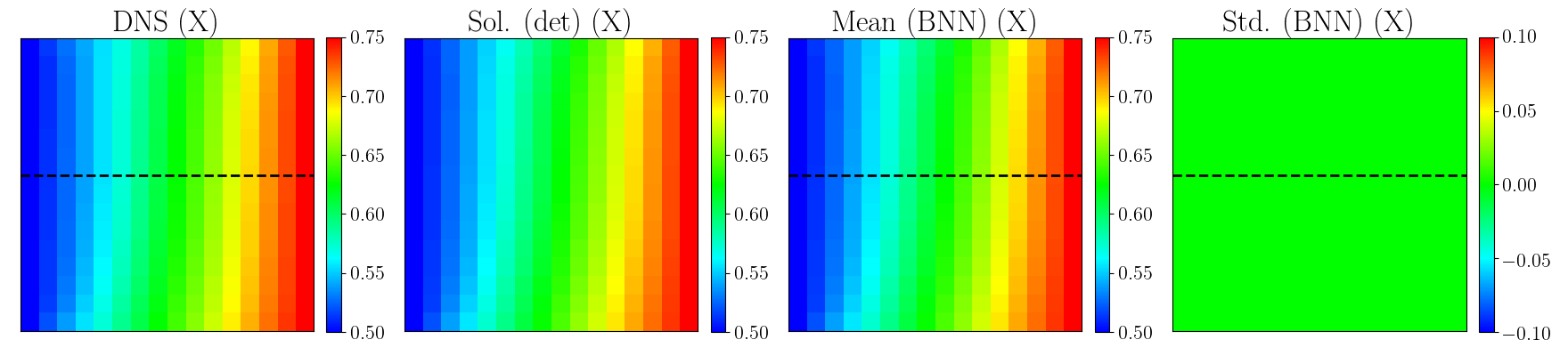}} 
  \subfloat[$u_x$ UQ for BVP (ii)]      {\includegraphics[trim=0.1cm 0.5cm 0.1cm 0.25cm, clip,height=0.15\linewidth]{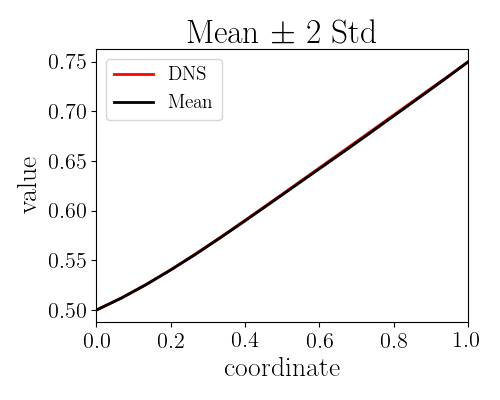}} \\[-2mm]
  \subfloat[$u_y$ results for BVP (ii)] {\includegraphics[trim=0.1cm 0.5cm 0.1cm 0.25cm, clip,height=0.15\linewidth]{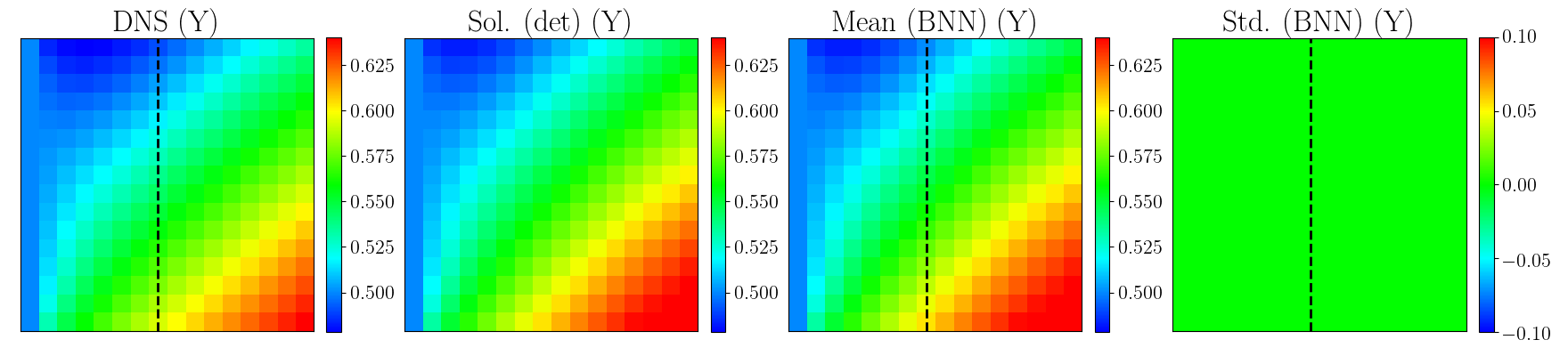}} 
  \subfloat[$u_y$ UQ for BVP (ii)]      {\includegraphics[trim=0.1cm 0.5cm 0.1cm 0.25cm, clip,height=0.15\linewidth]{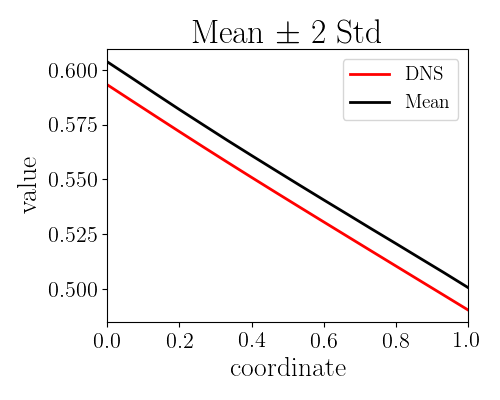}} \\[-2mm]
  \subfloat[$u_x$ results for BVP (iii)]{\includegraphics[trim=0.1cm 0.5cm 0.1cm 0.25cm, clip,height=0.15\linewidth]{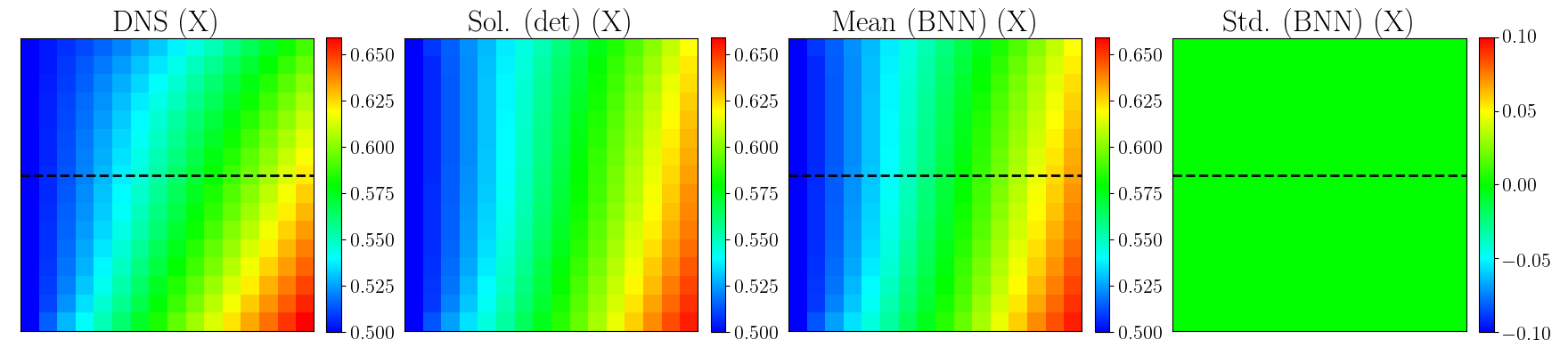}} 
  \subfloat[$u_x$ UQ for BVP (iii)]     {\includegraphics[trim=0.1cm 0.5cm 0.1cm 0.25cm, clip,height=0.15\linewidth]{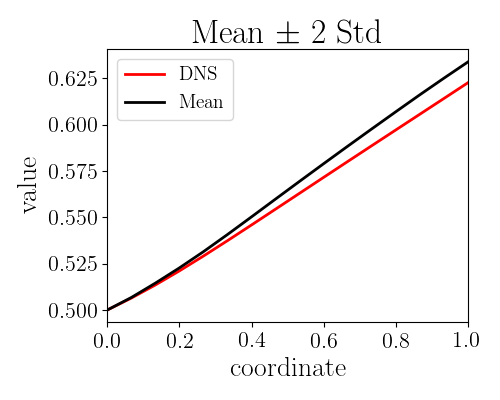}} \\[-2mm]
  \subfloat[$u_y$ results for BVP (iii)]{\includegraphics[trim=0.1cm 0.5cm 0.1cm 0.25cm, clip,height=0.15\linewidth]{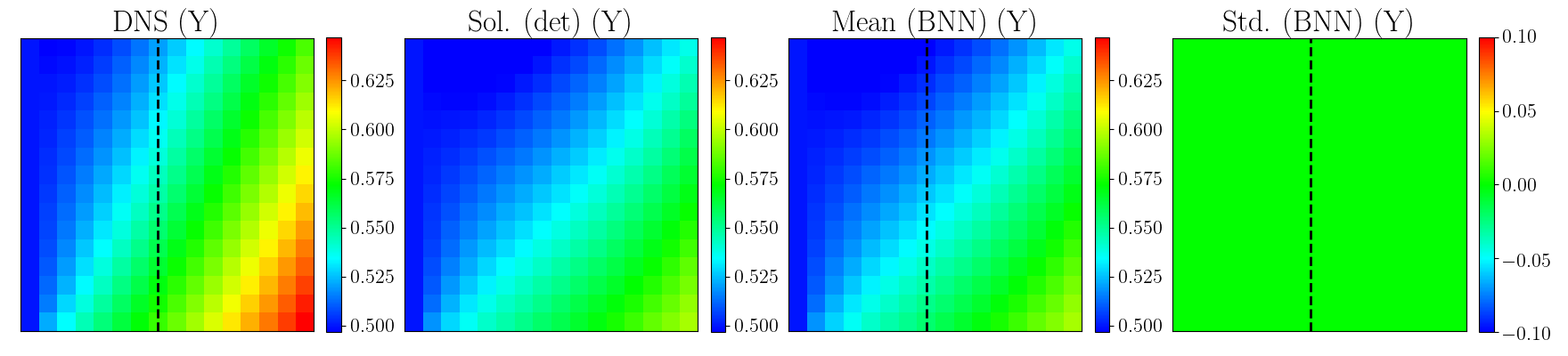}} 
  \subfloat[$u_y$ UQ for BVP (iii)]     {\includegraphics[trim=0.1cm 0.5cm 0.1cm 0.25cm, clip,height=0.15\linewidth]{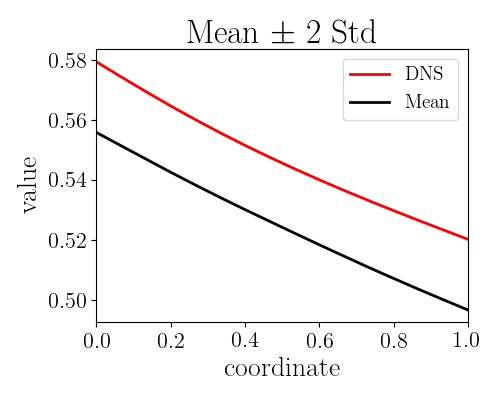}} \\
  \caption{Results of three selected BVPs out of the 30 nonlinear elasticity BVPs with varying domains and different applied BCs simultaneously solved by a single deterministic or probabilistic NN with the proposed method. BVP (i), (ii), (iii) correspond to bc id 1, 2, and 3 for domain id 5, as shown in Fig. \ref{fig:bvp-problem-setup}(a). 
      (a, c, e, g, i, k) Solutions from DNS, deterministic (det) NNs, and BNNs (Mean, Std.) for different BVPs.
  (b, d, f, h, j, l) Quantitative comparison of the solution distribution between DNS and BNNs along the dashed lines. 
}
  \label{fig:nonlinear-30bvp-results}
\end{figure}

In this section, we use the proposed PDE constrained NNs to simultaneously solve 30 nonlinear elasticity BVPs, as show in Fig. \ref{fig:bvp-problem-setup}(a), with a resolution of $16\times16$. 
The deformed problem domains from DNSs for three representative setups are shown in Fig. \ref{fig:nonlinear-30bvp-deformed}.
The architectures of both deterministic and probabilistic NNs and the training related NN parameters used in this section are identical to those used in section \ref{sec:linear-rectangle} for solving linear elasticity BVPs. 
We follow the procedures described in Section \ref{sec:NN-training} to train both types of NNs. 
The NN results of three selected BVPs, as shown in Fig. \ref{fig:nonlinear-30bvp-deformed},  are presented in Fig. \ref{fig:nonlinear-30bvp-results}.
The statistical moments of the BNN predictions are evaluated based on 50 MC samplings.
In Fig. \ref{fig:nonlinear-30bvp-results}, BVP (i), (ii), and (iii) correspond to bc id 1 (non-zero Dirichlet loading), bc id 2 (non-zero Neumann loading), bc id 3 (mixed loading) applied to domain id 5. 
The comparison of solutions between DNSs, the deterministic NN, and the BNN for these three BVPs is shown qualitatively in Fig. \ref{fig:nonlinear-30bvp-results}(a,c,e,g,i,k), with quantitative comparison of the solution distribution along the dashed lines between DNSs and the BNN given in Fig. \ref{fig:nonlinear-30bvp-results}(b,d,f,h,j,l). We draw attention to the improved accuracy of the BNN for BVP (iii) seen in Fig. \ref{fig:nonlinear-30bvp-results}(i,j,k,l) as a consequence of larger deformation (strain) in comparison to the same BVP with inifinitesimal strain in Fig. \ref{fig:linear-30bvp-results}(i,j,k,l). The BNN is able to better learn nonlinear than linear elasticity because of the stronger expression of physics in the nonlinear solution.
The comparison demonstrates that the proposed method has successfully solved multiple BVPs with desired accuracy.

\subsubsection{Rectangular domain with solution interpolation and extrapolation}\label{sec:nonlinear-rectangle-interpolation}

\begin{figure}[h!]
  \centering
  {\includegraphics[height=0.08\linewidth]{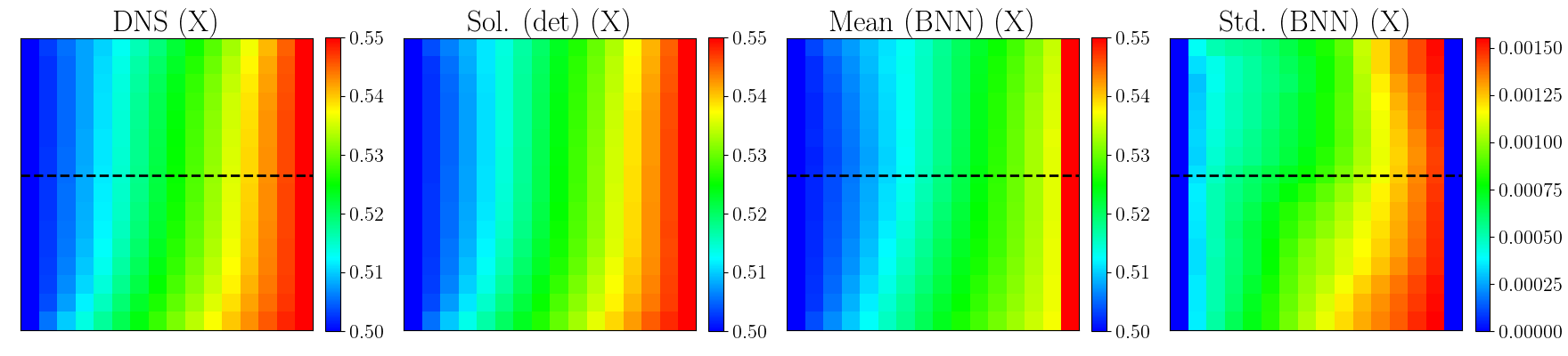}}
  {\includegraphics[height=0.08\linewidth]{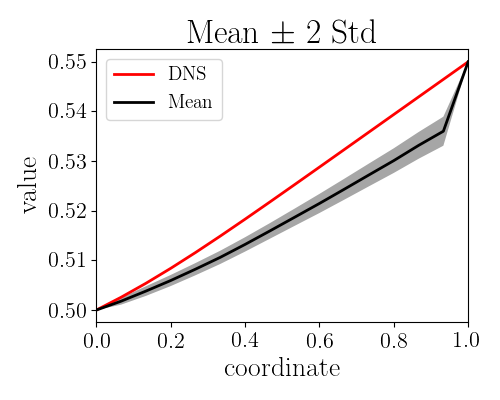}} 
  {\includegraphics[height=0.08\linewidth]{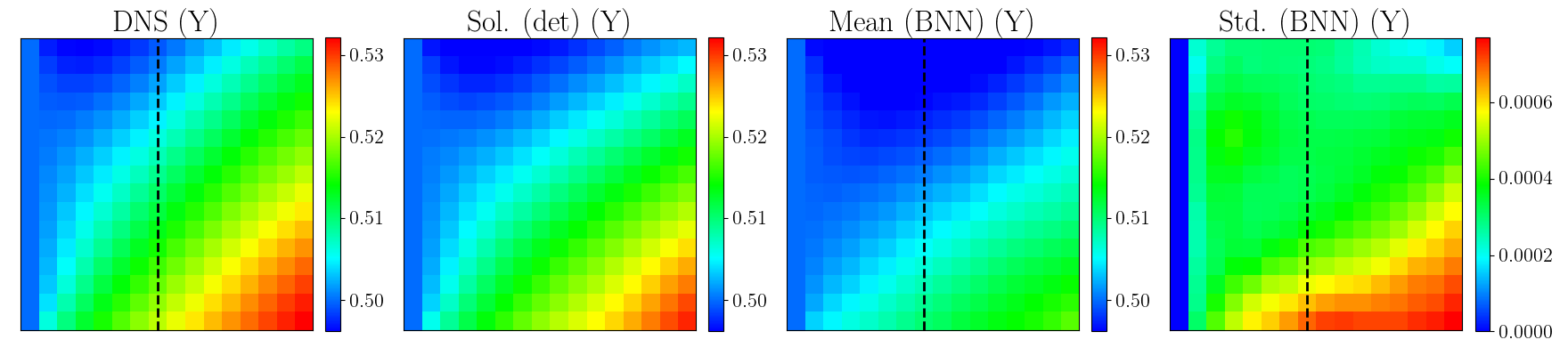}}
  {\includegraphics[height=0.08\linewidth]{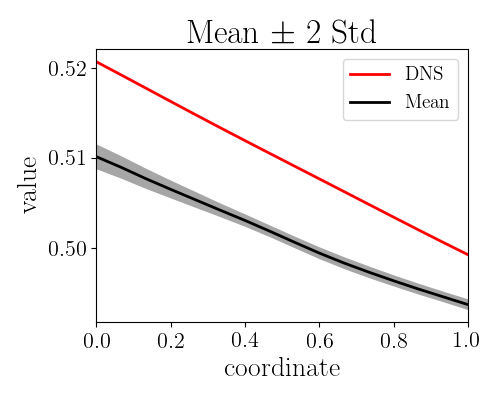}} \\
  {\includegraphics[height=0.08\linewidth]{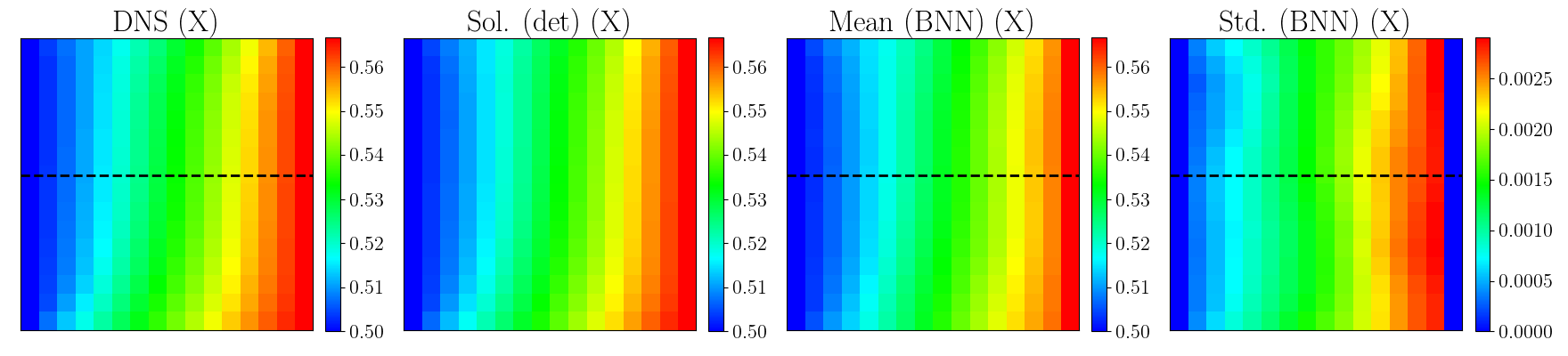}}
  {\includegraphics[height=0.08\linewidth]{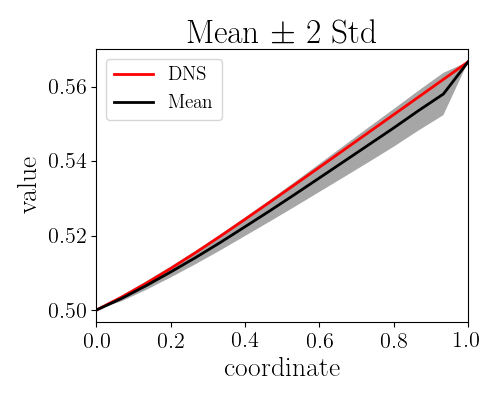}} 
  {\includegraphics[height=0.08\linewidth]{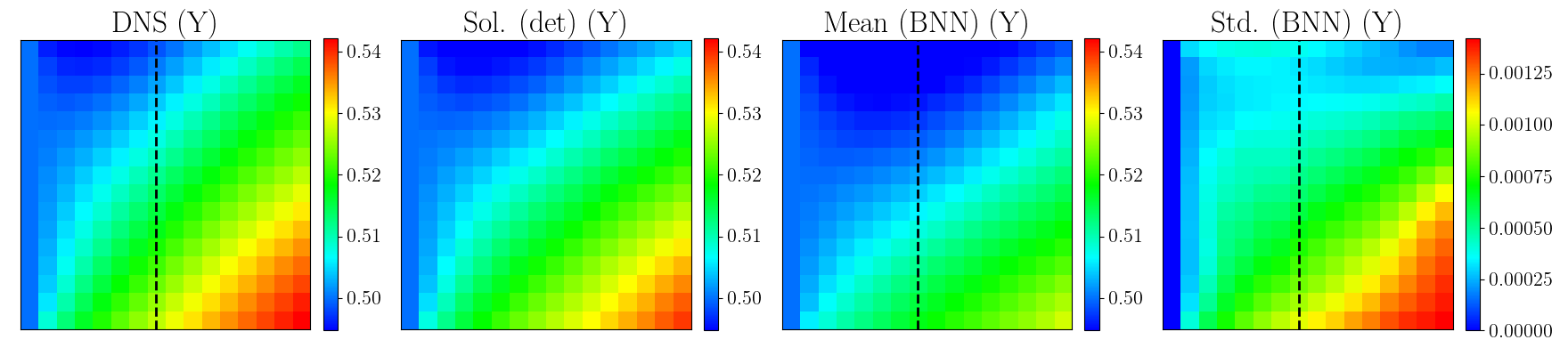}}
  {\includegraphics[height=0.08\linewidth]{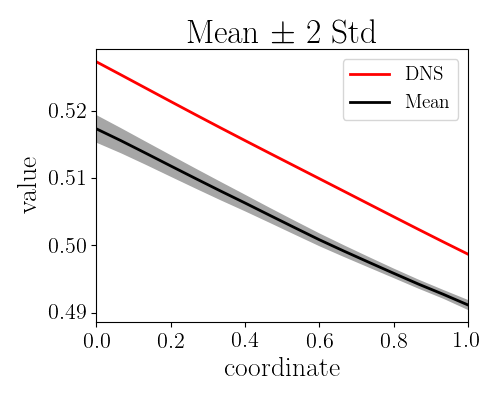}} \\
  {\includegraphics[height=0.08\linewidth]{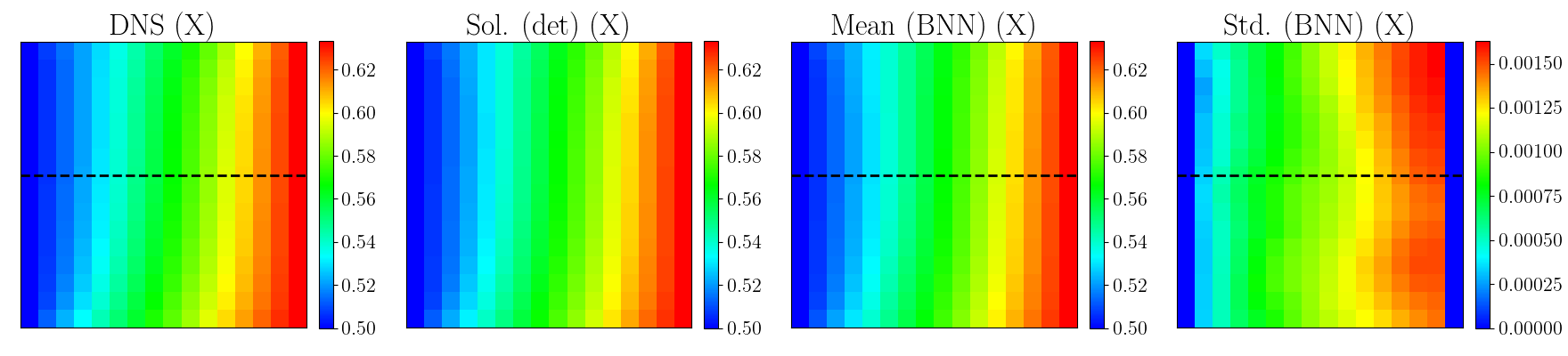}}
  {\includegraphics[height=0.08\linewidth]{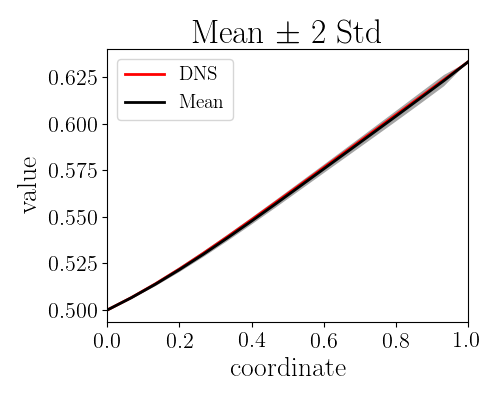}} 
  {\includegraphics[height=0.08\linewidth]{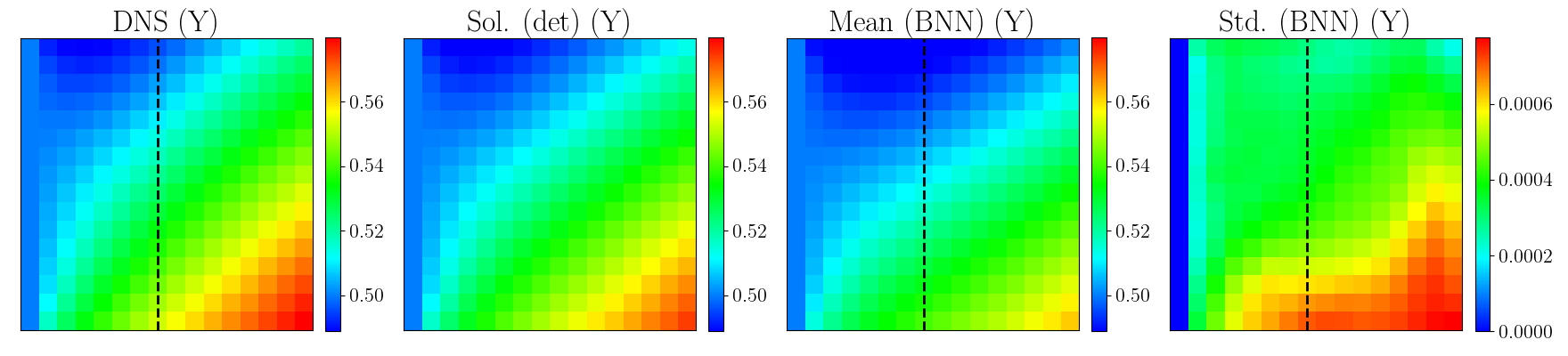}}
  {\includegraphics[height=0.08\linewidth]{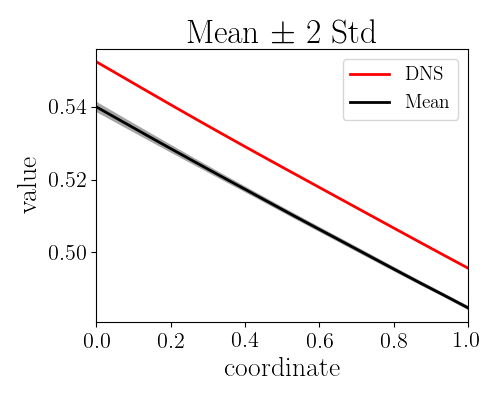}} \\
  \caption{Additional interpolated NN prediction results for case (i). }
  \label{fig:nonlinear-inter-additional}
\end{figure}

\begin{figure}[h!]
  \centering
  {\includegraphics[height=0.08\linewidth]{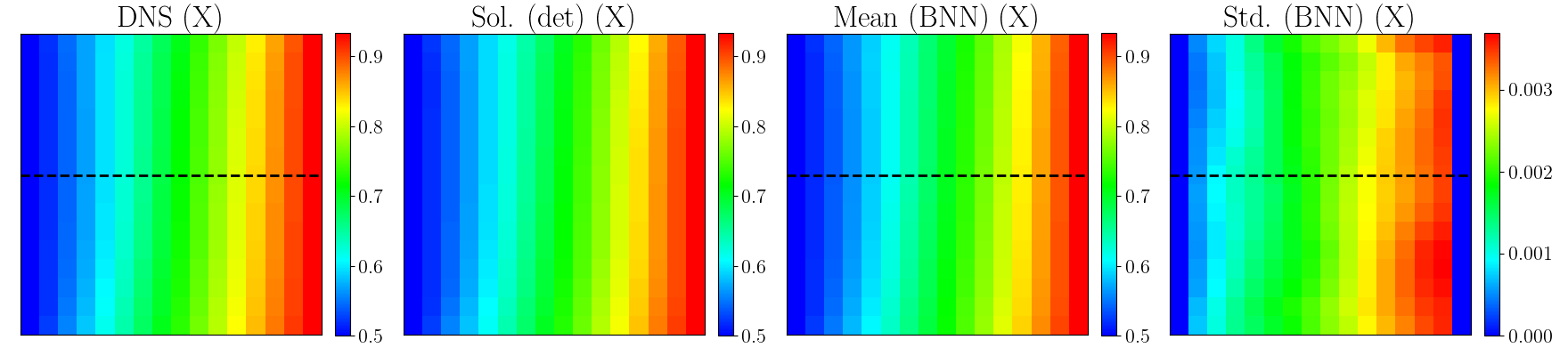}}
  {\includegraphics[height=0.08\linewidth]{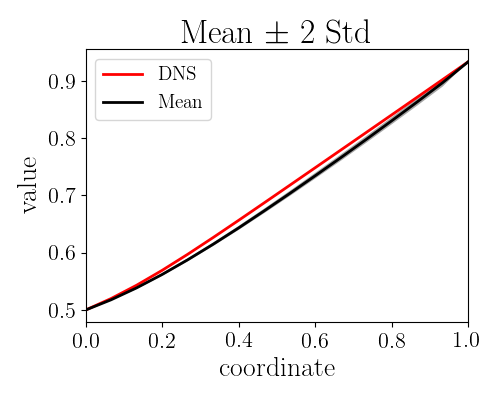}} 
  {\includegraphics[height=0.08\linewidth]{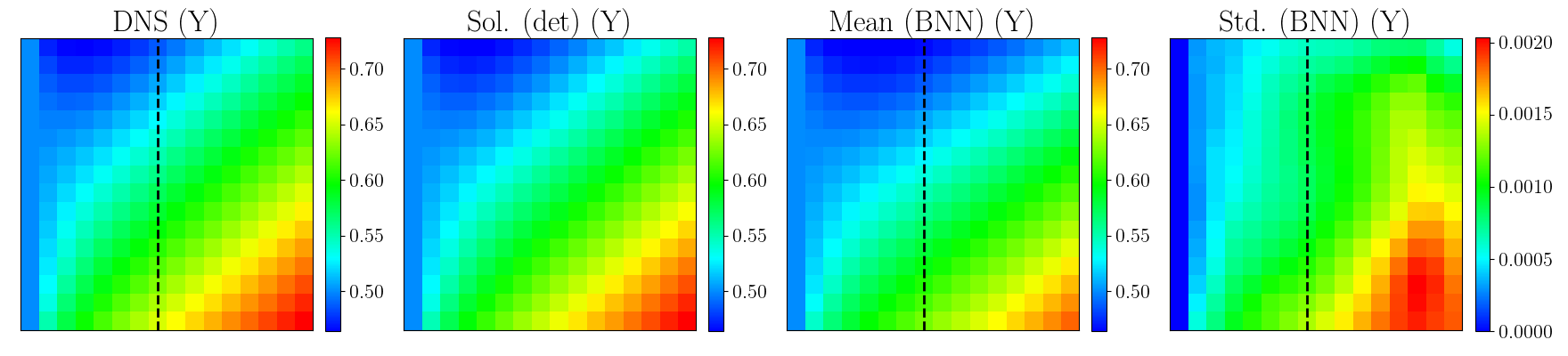}}
  {\includegraphics[height=0.08\linewidth]{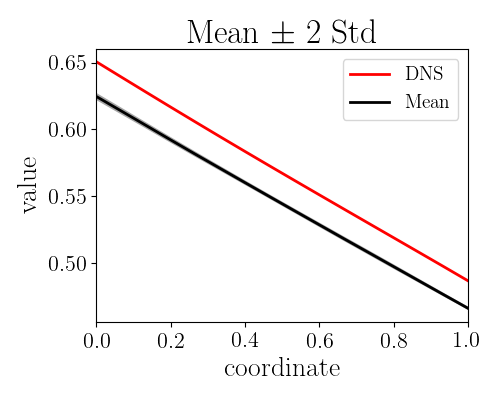}} \\
  {\includegraphics[height=0.08\linewidth]{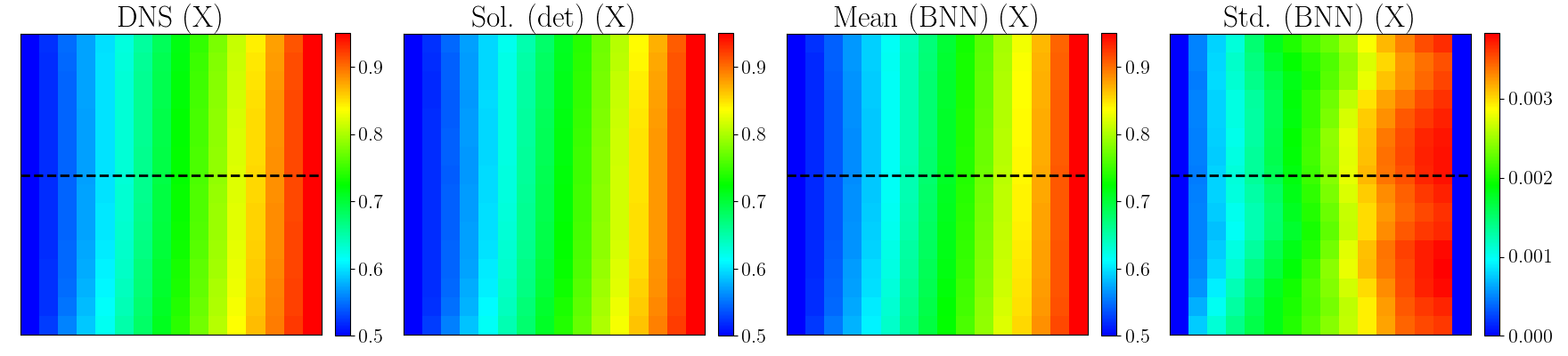}}
  {\includegraphics[height=0.08\linewidth]{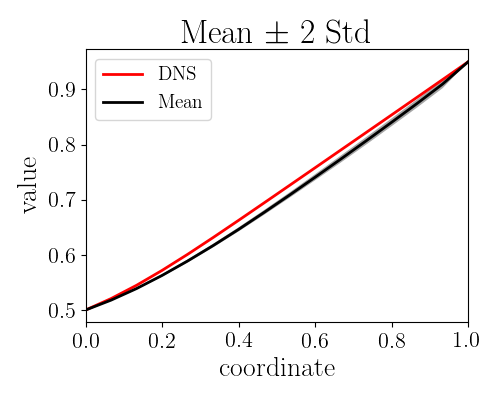}} 
  {\includegraphics[height=0.08\linewidth]{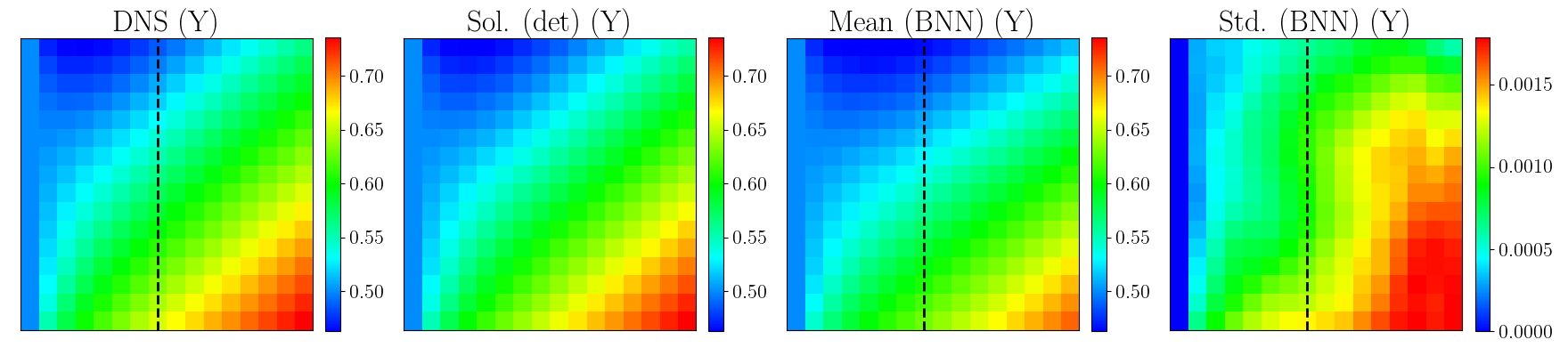}}
  {\includegraphics[height=0.08\linewidth]{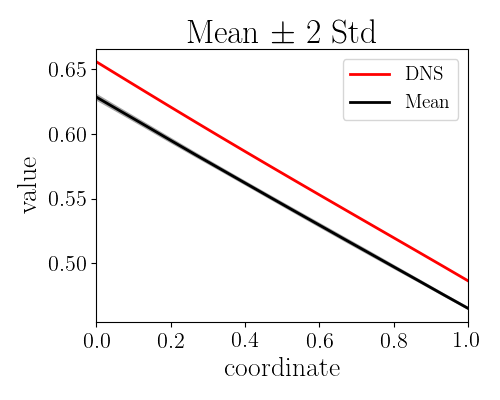}} \\
  {\includegraphics[height=0.08\linewidth]{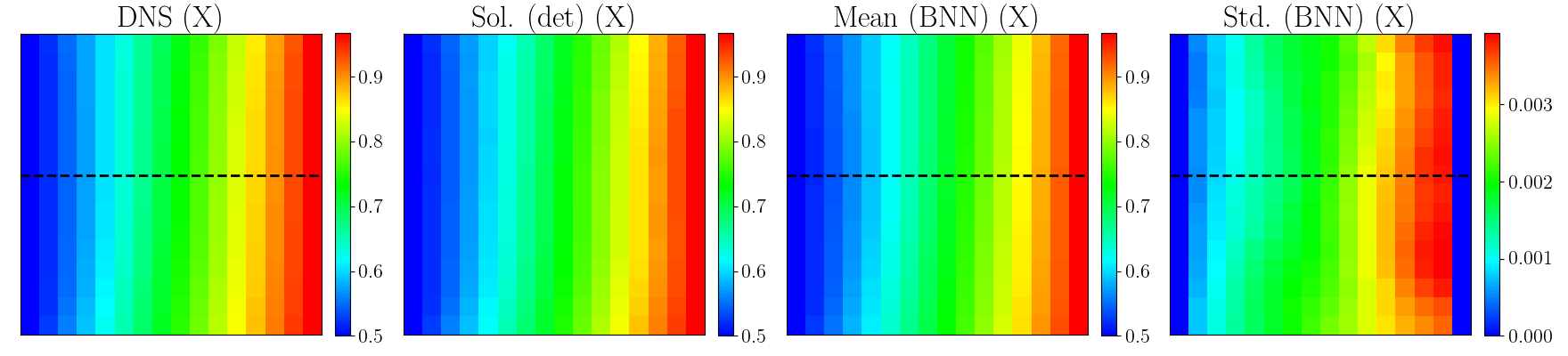}}
  {\includegraphics[height=0.08\linewidth]{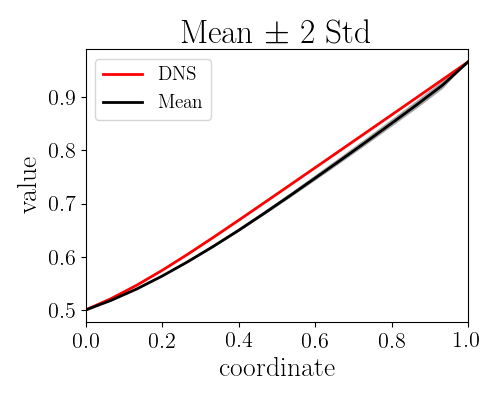}} 
  {\includegraphics[height=0.08\linewidth]{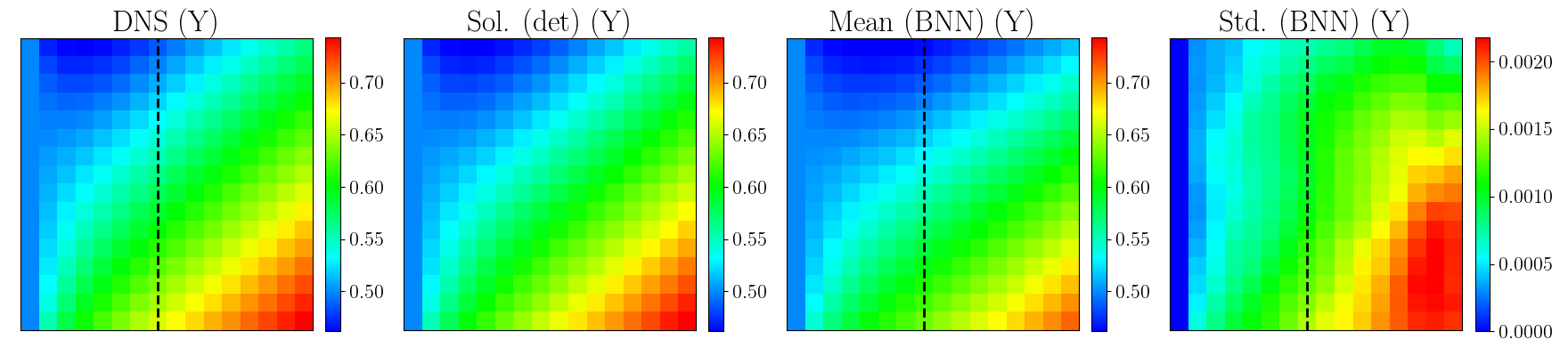}}
  {\includegraphics[height=0.08\linewidth]{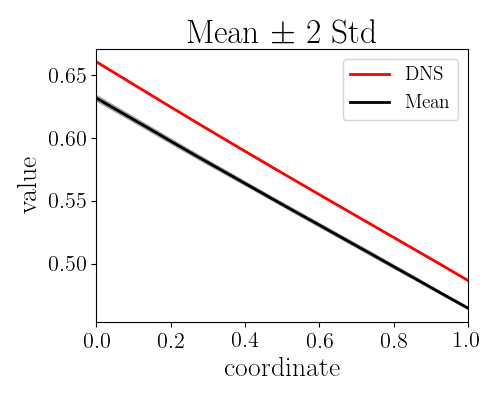}} \\
  \caption{Additional extrapolated NN prediction results for case (i). }
  \label{fig:nonlinear-extra-additional}
\end{figure}

In this section, we explore the interpolating and extrapolating capacity of the proposed framework for the BVP setup shown in Fig. \ref{fig:nonlinear-30bvp-deformed}, case (i). Both the DNS and NN solution have resolutions of $16\times16$.
The architectures of both deterministic and probabilistic NNs and the training related NN parameters used in this section are identical to those used in section \ref{sec:linear-rectangle} for solving linear elasticity BVPs. 
Additional interpolated and extrapolated NN prediction results for case (i) are given in Fig. \ref{fig:nonlinear-inter-additional} and \ref{fig:nonlinear-extra-additional}, respectively.

\subsection{Steady state diffusion - large dataset}
In this section, three large training datasets (D1, D2, D3) and two small testing datasets (T1, T2) with a mesh resolution of $64\times64$ are prepared for the steady-state diffusion problem to test the performance of the proposed method.
The script to synthetically generate the BVPs in these datasets are provided in the source code on GitHub. 

\subsubsection{D1/T1: BVPs on regular pentagon domains}
\begin{figure}[h!]
  \centering
  \includegraphics[width=0.7\linewidth]{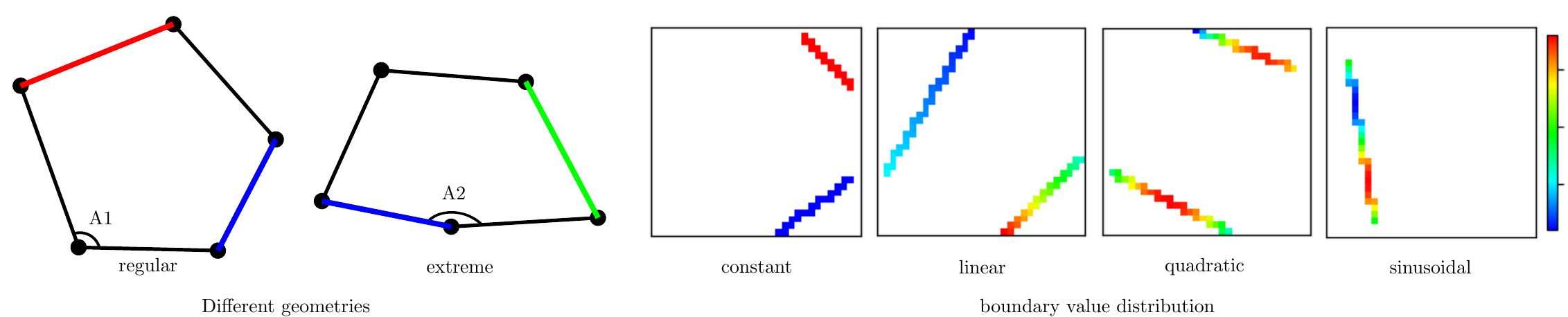}
  \caption{Illustration of dataset D1/T1. D1 contains regular pentagons, whose inner angles are all smaller than 160 degree. T1 contains both regular and irregular pentagons. The latter has one innger angle greater than 160 degree. The boundary values could have a constant/linear/quadratic/sinusoidal distribution.}
  \label{fig:large-D1}
\end{figure}
\begin{figure}[h!]
  \centering
  \includegraphics[height=45mm]{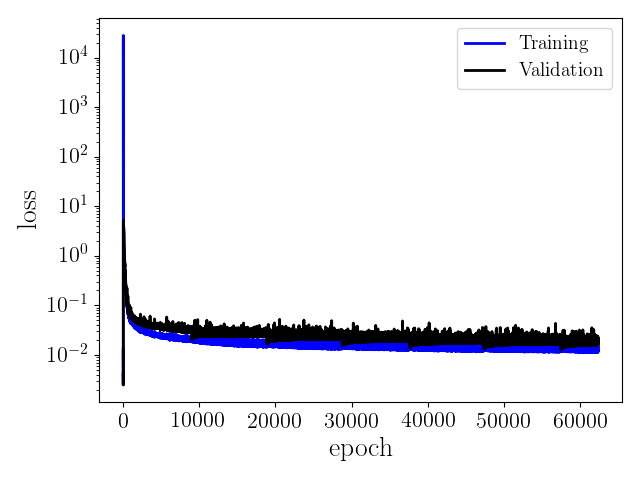}
  \includegraphics[height=45mm]{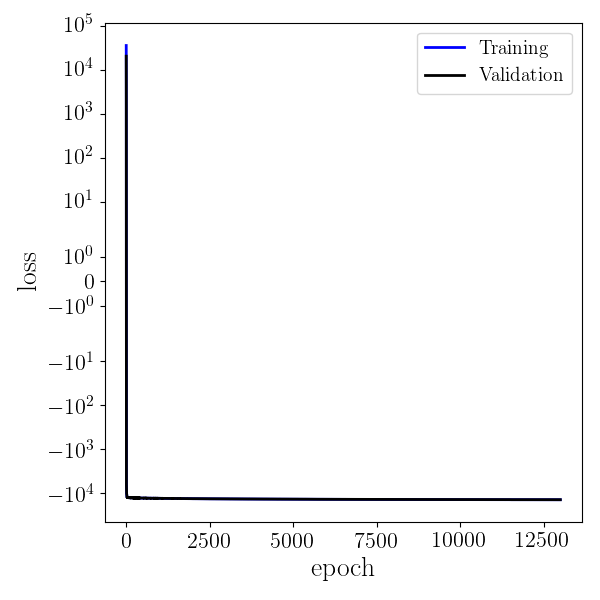}
  \includegraphics[height=45mm]{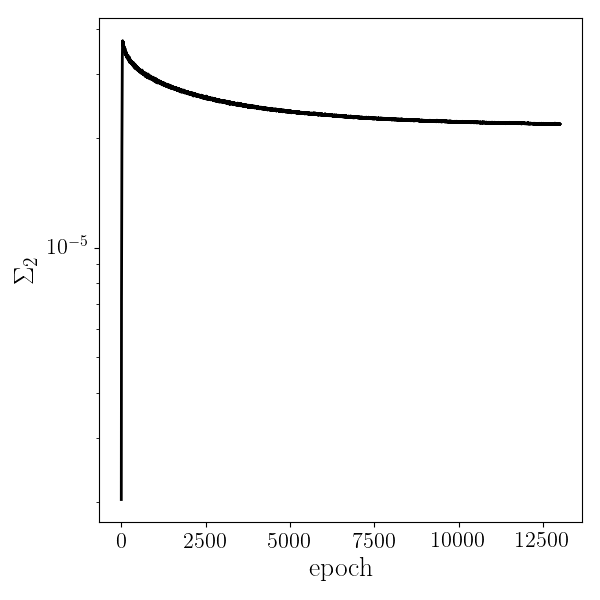}
  \caption{The training loss of the residual constrained NNs to solve all BVPs in D1.}
  \label{fig:D1-loss}
\end{figure}

D1 contains 176K unique BVPs, which covers 1100 regular pentagons, whose inner angles are all smaller than 160 degree, as shown in Fig. \ref{fig:large-D1}. The shape of the pentagons are randomly generated. We apply boundary conditions to two non-adjacent edges, which are randomly sampled from all the possible locations. The boundary values are randomly generated, which could have a constant/linear/quadratic/sinusoidal distribution along the edges, as shown in Fig. \ref{fig:large-D1}. The testing dataset T1 contains 1000 BVPs, which covers both regular shapes and extreme shapes. The latter has one inner angle greater than 160 degree. Same strategies to generate the BCs for D1 are used to generate BCs for T1.
The training loss is shown in Fig. \ref{fig:D1-loss}.
NN architectures and training details are summarized in Table \ref{tab:D1-NNs} and \ref{tab:D1-NNs-others}.

\begin{table}
  \centering
  \begin{tabular}{l | l | l | l}
    \hline
    Deterministic         & Probabilistic         & Size         & Layer arguments \\ \hline
    Input                 & Input                 & -            & - \\
    LayerFillRandomNumber & LayerFillRandomNumber & -            & - \\
    Conv2D                & Convolution2DFlipout  & filters = 8  & kernel (5,5), padding: same, ReLU \\
    MaxPooling2D          & MaxPooling2D          & -            & kernel (2,2), padding: same\\
    Conv2D                & Convolution2DFlipout  & filters = 8 & kernel (5,5), padding: same, ReLU \\
    MaxPooling2D          & MaxPooling2D          & -            & kernel (2,2), padding: same\\
    Conv2D                & Convolution2DFlipout  & filters = 16 & kernel (5,5), padding: same, ReLU \\
    MaxPooling2D          & MaxPooling2D          & -            & kernel (2,2), padding: same\\
    Flatten               & Flatten               & -            & - \\
    Dense                 & DenseFlipout          & units = 32   & ReLU \\
    Dense                 & DenseFlipout          & units = 128   & ReLU \\
    Reshape               & Reshape               & -            & $[4,4,8]$ \\
    Conv2D                & Convolution2DFlipout  & filters = 16 & kernel (5,5), padding: same, ReLU \\
    UpSampling2D          & UpSampling2D          & -            & size (2,2) \\
    Conv2D                & Convolution2DFlipout  & filters = 16 & kernel (5,5), padding: same, ReLU \\
    UpSampling2D          & UpSampling2D          & -            & size (2,2) \\
    Conv2D                & Convolution2DFlipout  & filters = 16 & kernel (5,5), padding: same, ReLU \\
    UpSampling2D          & UpSampling2D          & -            & size (2,2) \\
    Conv2D                & Convolution2DFlipout  & filters = 16 & kernel (5,5), padding: same, ReLU \\
    Conv2D                & Convolution2DFlipout  & filters = 2  & kernel (5,5), padding: same, ReLU \\
    \hline
  \end{tabular}
  \caption{Details of both deterministic and probabilistic NNs for solving D1.}
  \label{tab:D1-NNs}
\end{table}

\begin{table}
  \centering
  \begin{tabular}{l | l | l }
    \hline
    Description                     & Deterministic                 & Probabilistic         \\ \hline
    Total parameters                & 41,346                        & 82,435                 \\
    Size of $\calD$                 & 5 $\times$ Aug: $2^{10}$      & 5 $\times$ Aug: $2^{9}$      \\
    Epochs                          & 10,000                        & 100                 \\
    Zero initialization epochs      & 100                           & -                     \\
    Optimizer                       & Adam                          & Nadam                 \\
    Learning Rate                   & 2.5e-4                        & 1e-8                  \\
    Batch Size                      & 256                           & 64                    \\
    $\Sigma_1$                      & -                             & 1e-8                  \\
    Initial value of $\Sigma_2$     & -                             & 1e-8                  \\
    \hline
  \end{tabular}
  \caption{Training related parameters for solving D1.}
  \label{tab:D1-NNs-others}
\end{table}

Selected NN predicted results are shown in Fig. \ref{fig:large-diffusion}(a,b), with the $L_2$ error of all the training and testing results given in Fig. \ref{fig:large-diffusion}(d).
We can observe that, if the boundary locations and the extremes of boundary values in the training dataset could be more uniformly distributed, the training and prediction are in generally good. 
Also,  the predictions are more accurate if extremes of the testing dataset are not near the extremes of training datasets. 
This suggests that if we know a targeting range to use the NNs to make predictions, we can then design the NN training dataset to have a wider range, so the accuracy of NN predicted solution can be further improved. In another word, make interpolated prediction. 
This study also shows that a single NN could simultaneously solve hundreds of thousands BVPs with reasonable accuracy and make prediction for unseen domains.

\begin{figure}[h!]
  \centering
  \includegraphics[width=1.2\linewidth]{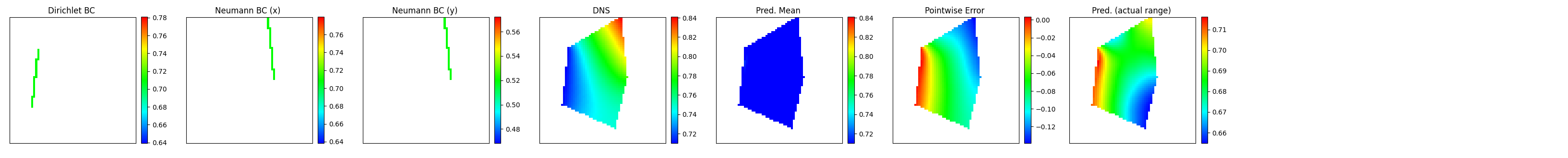}
  \includegraphics[width=1.2\linewidth]{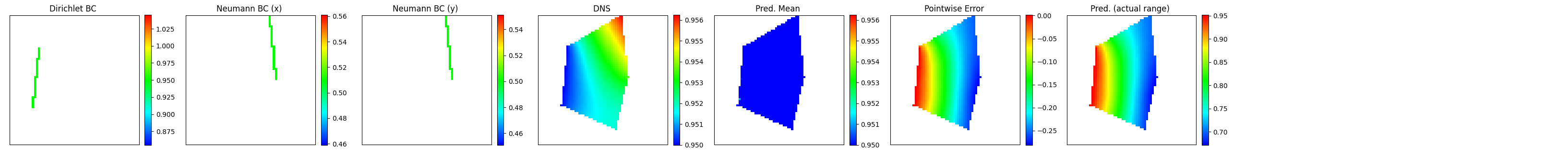}
  \includegraphics[width=1.2\linewidth]{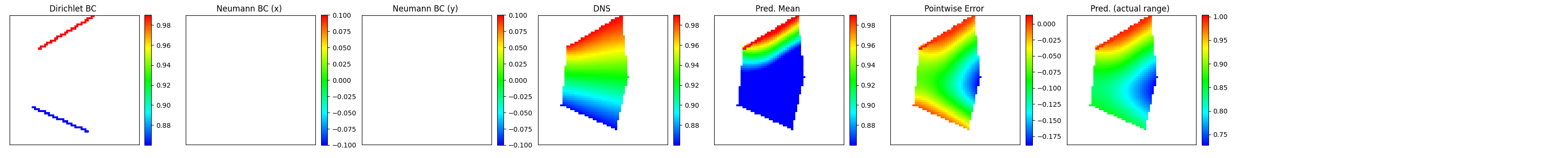}
  \caption{Poor results from NN predictions for T1 from the NN trained over D1.}
  \label{fig:D1-poor}
\end{figure}

\begin{figure}[h!]
  \centering
  \includegraphics[width=0.24\linewidth]{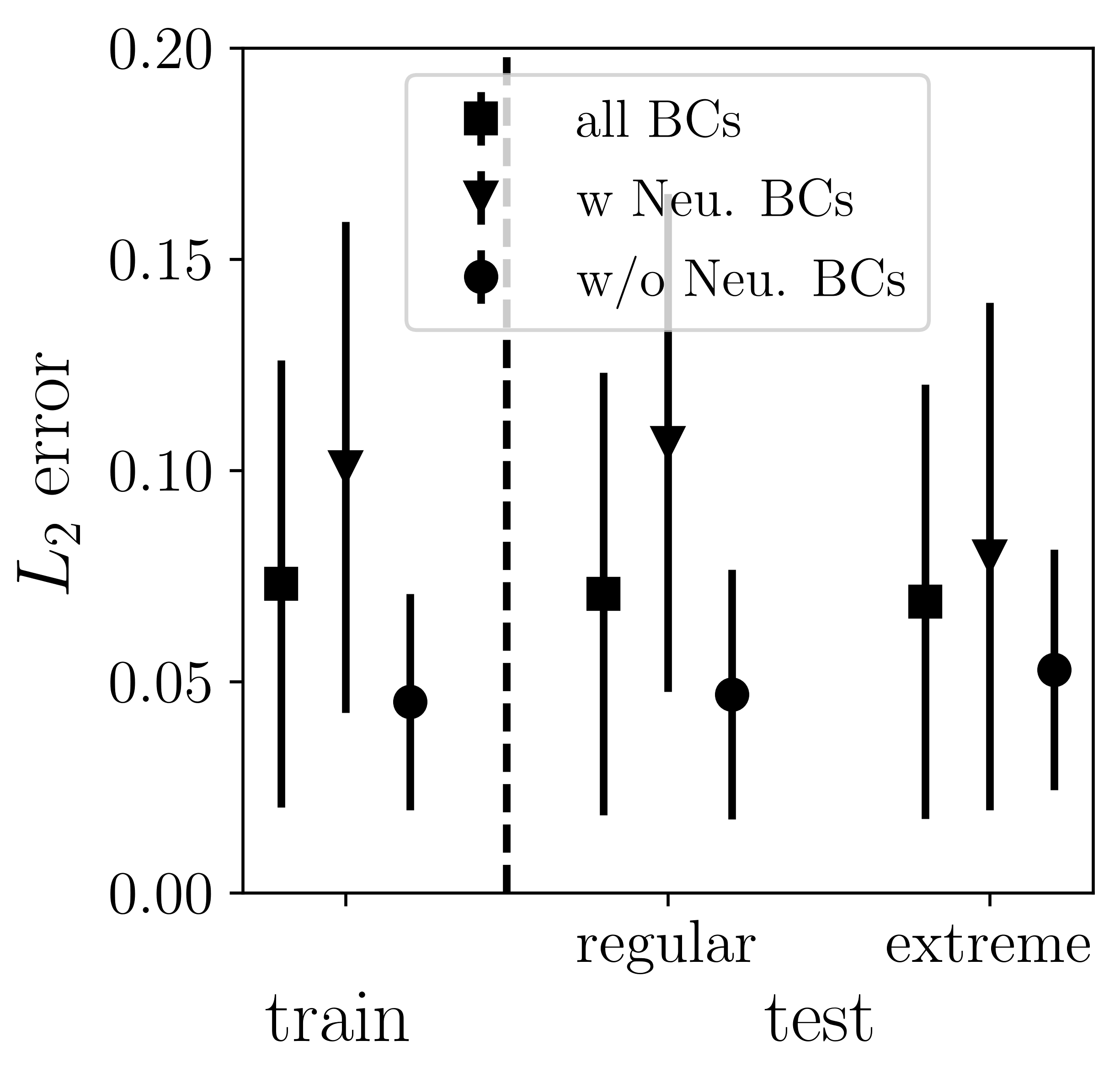}
  \includegraphics[width=0.24\linewidth]{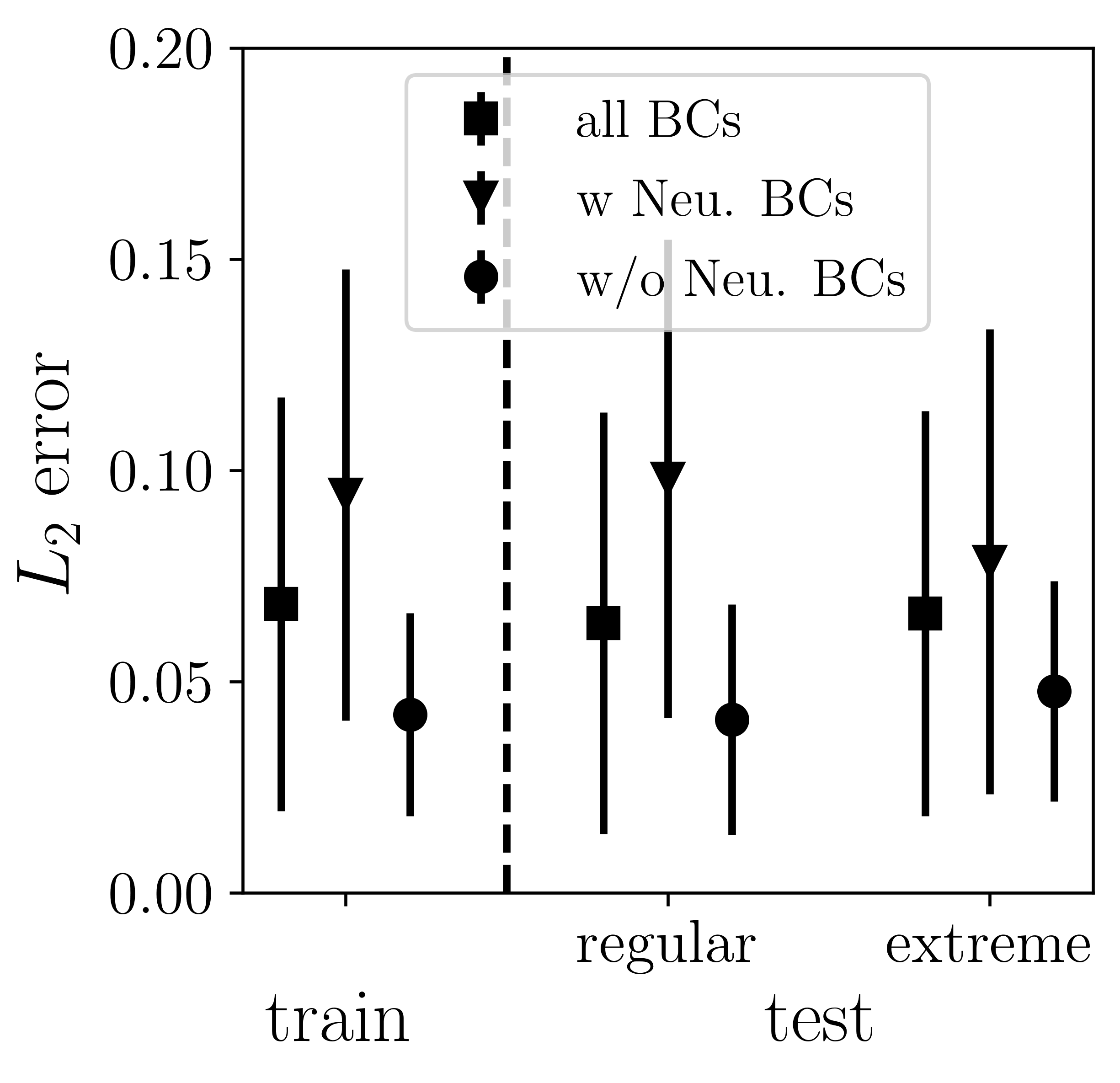}
  \includegraphics[width=0.24\linewidth]{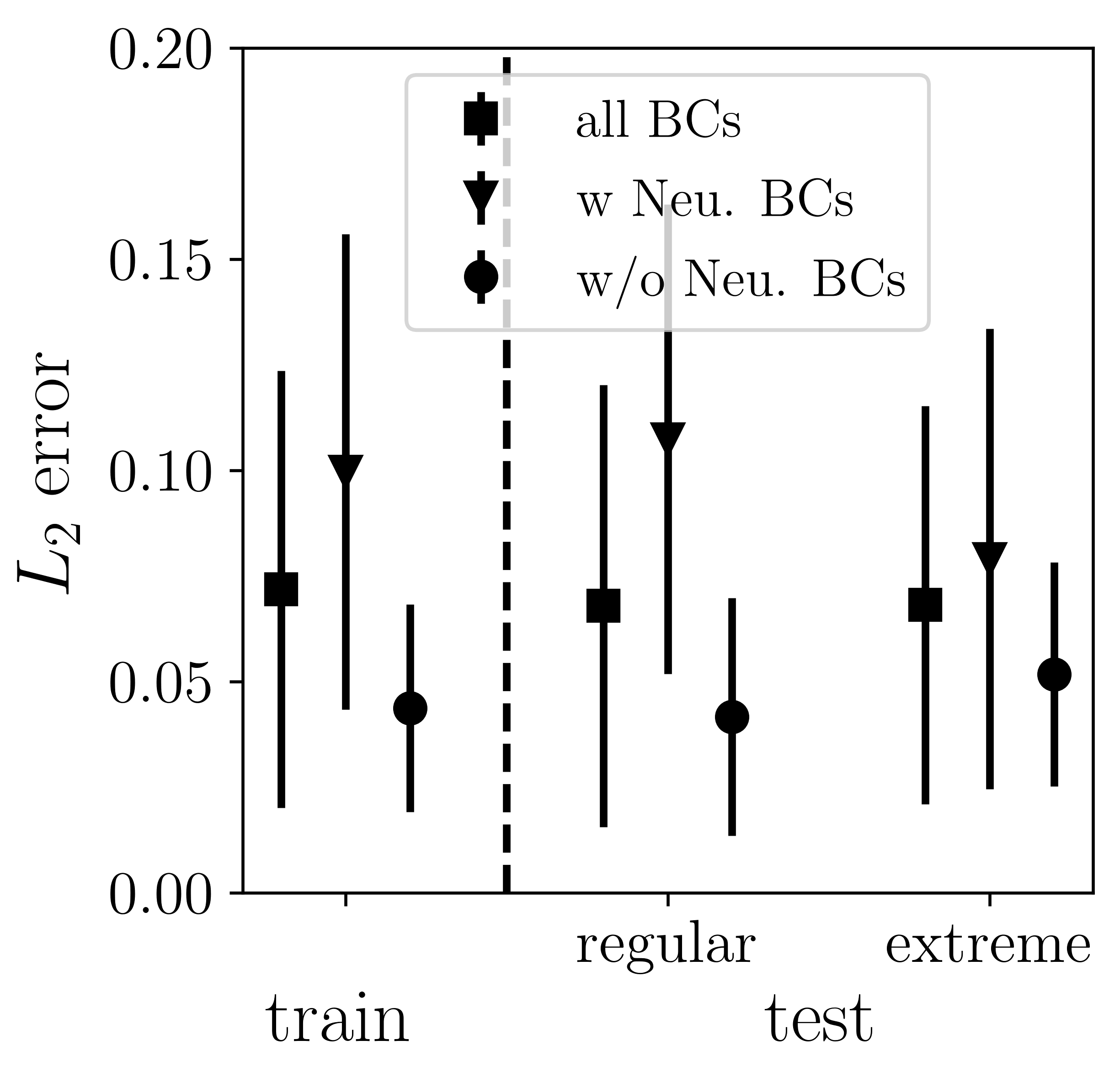}
  \includegraphics[width=0.24\linewidth]{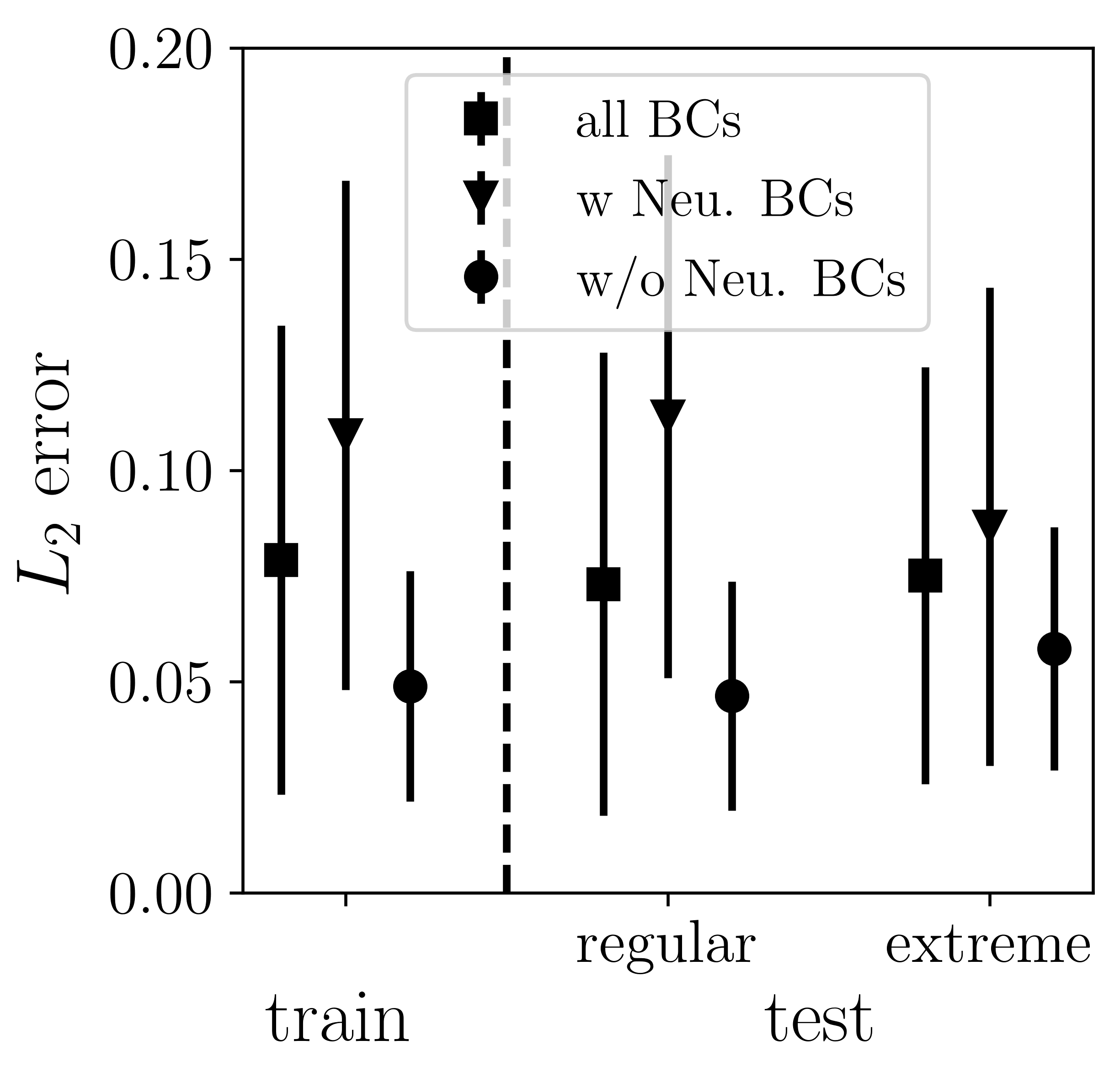}
  \caption{Statistics of NN results for different training to solve D1 to confirm the repeatability of the proposed method. }
  \label{fig:D1-NN-stats}
\end{figure}

\subsubsection{D2/T2: BVPs on quadrilateral/pentagon/hexagon}
\begin{figure}[h!]
  \centering
  \includegraphics[width=1.0\linewidth]{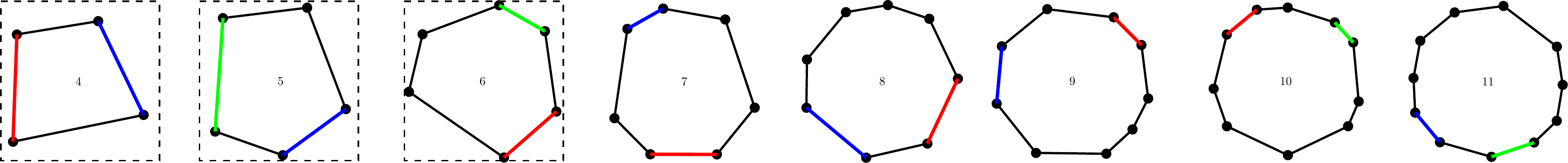}
  \caption{Illustration of dataset D2/T2. D2 contains BVPs that covers quadrilateral, pentagons, and hexagon. T2 contains polygons with the total number of edges ranging from four to eleven. Similar as for D1/T1, the boundary values could have a constant/linear/quadratic/sinusoidal distribution.}
  \label{fig:large-D2}
\end{figure}

\begin{figure}[h!]
  \centering
  \includegraphics[height=45mm]{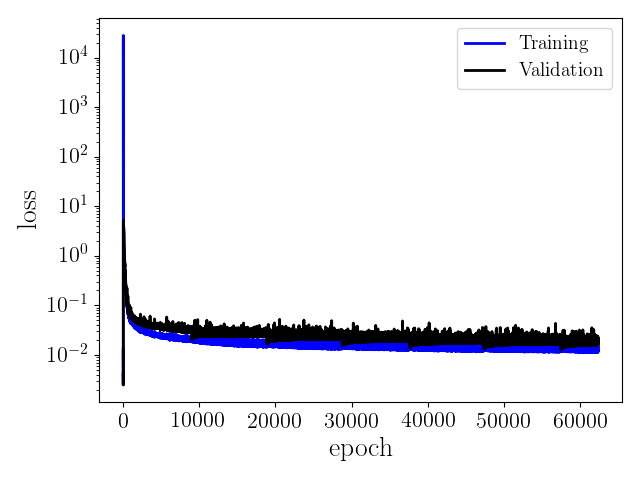}
  \includegraphics[height=45mm]{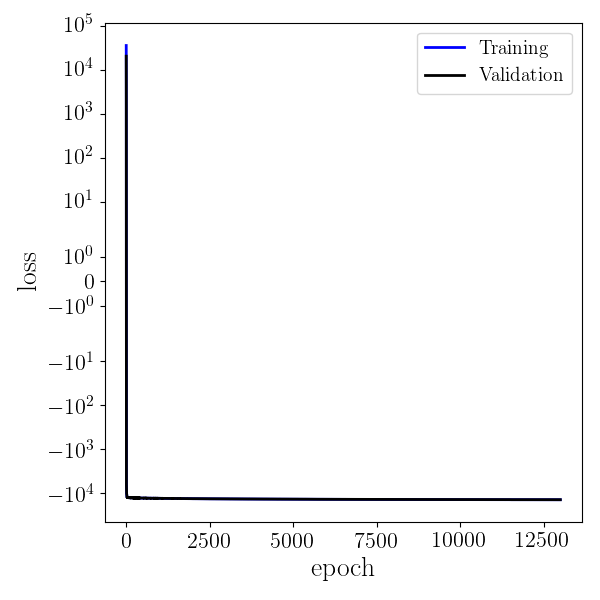}
  \includegraphics[height=45mm]{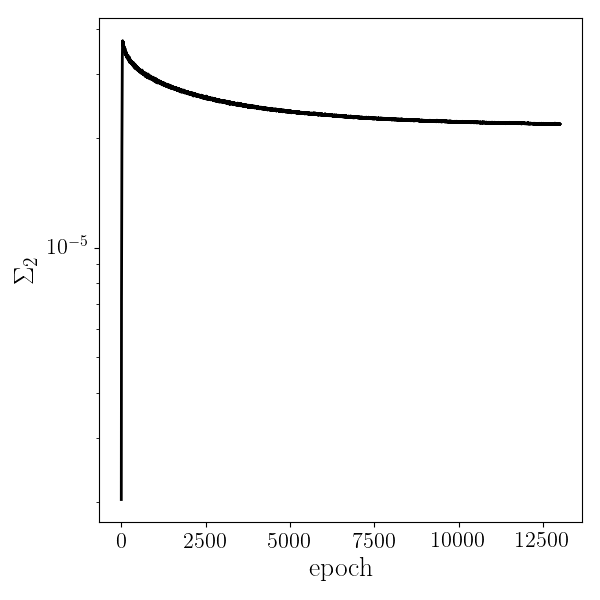}
  \caption{The training loss of the residual constrained NNs to solve all BVPs in D2. }
  \label{fig:D2-loss}
\end{figure}

D2 contains 192K unique BVPs, which consists of 32K BVPs on quadrilateral, 64K BVPs on pentagons, and 96K BVPs on hexagon.
These BVPs covers 400 quadrilaterals, 400 pentagons, and 400 hexagons. 
The shape of these polygons are randomly generated. As for D1/T1, we apply boundary conditions to two non-adjacent edges, which are randomly sampled from all the possible locations. The boundary values are randomly generated, which could have a constant/linear/quadratic/sinusoidal distribution along the edges. 
The testing dataset T2 contains 320 BVPs, which covers polygons with the total number of edges range from four to eleven. 
The training loss is shown in Fig. \ref{fig:D2-loss}.
NN architectures and training details are summarized in Table \ref{tab:D2-NNs} and \ref{tab:D2-NNs-others}. Selected NN predicted results are shown in Fig. \ref{fig:large-diffusion}(c), with the $L_2$ error of all the training and testing results given in Fig. \ref{fig:large-diffusion}(e).

\begin{table}
  \centering
  \begin{tabular}{l | l | l | l}
    \hline
    Deterministic         & Probabilistic         & Size         & Layer arguments \\ \hline
    Input                 & Input                 & -            & - \\
    LayerFillRandomNumber & LayerFillRandomNumber & -            & - \\
    Conv2D                & Convolution2DFlipout  & filters = 8  & kernel (5,5), padding: same, ReLU \\
    MaxPooling2D          & MaxPooling2D          & -            & kernel (2,2), padding: same\\
    Conv2D                & Convolution2DFlipout  & filters = 8 & kernel (5,5), padding: same, ReLU \\
    MaxPooling2D          & MaxPooling2D          & -            & kernel (2,2), padding: same\\
    Conv2D                & Convolution2DFlipout  & filters = 16 & kernel (5,5), padding: same, ReLU \\
    MaxPooling2D          & MaxPooling2D          & -            & kernel (2,2), padding: same\\
    Flatten               & Flatten               & -            & - \\
    Dense                 & DenseFlipout          & units = 32   & ReLU \\
    Dense                 & DenseFlipout          & units = 128   & ReLU \\
    Reshape               & Reshape               & -            & $[4,4,8]$ \\
    Conv2D                & Convolution2DFlipout  & filters = 16 & kernel (5,5), padding: same, ReLU \\
    UpSampling2D          & UpSampling2D          & -            & size (2,2) \\
    Conv2D                & Convolution2DFlipout  & filters = 16 & kernel (5,5), padding: same, ReLU \\
    UpSampling2D          & UpSampling2D          & -            & size (2,2) \\
    Conv2D                & Convolution2DFlipout  & filters = 16 & kernel (5,5), padding: same, ReLU \\
    UpSampling2D          & UpSampling2D          & -            & size (2,2) \\
    Conv2D                & Convolution2DFlipout  & filters = 16 & kernel (5,5), padding: same, ReLU \\
    Conv2D                & Convolution2DFlipout  & filters = 2  & kernel (5,5), padding: same, ReLU \\
    \hline
  \end{tabular}
  \caption{Details of both deterministic and probabilistic NNs for solving D2.}
  \label{tab:D2-NNs}
\end{table}

\begin{table}
  \centering
  \begin{tabular}{l | l | l }
    \hline
    Description                     & Deterministic                 & Probabilistic         \\ \hline
    Total parameters                & 41,346                        & 82,435                 \\
    Size of $\calD$                 & 5 $\times$ Aug: $2^{10}$      & 5 $\times$ Aug: $2^{9}$      \\
    Epochs                          & 10,000                        & 100                 \\
    Zero initialization epochs      & 100                           & -                     \\
    Optimizer                       & Adam                          & Nadam                 \\
    Learning Rate                   & 2.5e-4                        & 1e-8                  \\
    Batch Size                      & 256                           & 64                    \\
    $\Sigma_1$                      & -                             & 1e-8                  \\
    Initial value of $\Sigma_2$     & -                             & 1e-8                  \\
    \hline
  \end{tabular}
  \caption{Training related parameters for solving D2.}
  \label{tab:D2-NNs-others}
\end{table}

\begin{figure}[h!]
  \centering
  \includegraphics[width=0.9\linewidth]{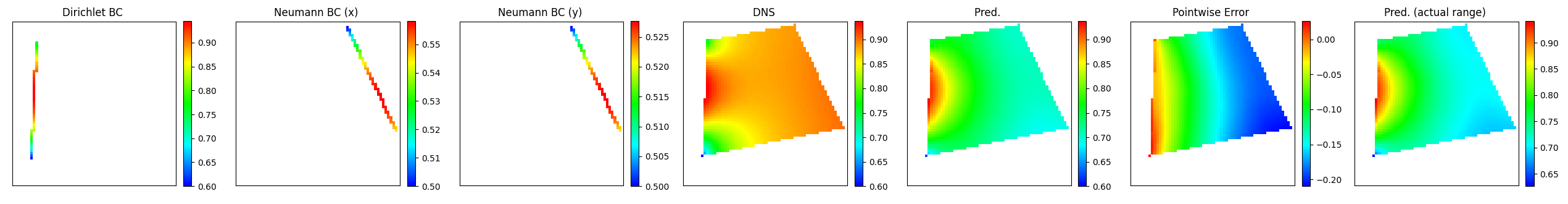}
  \includegraphics[width=0.9\linewidth]{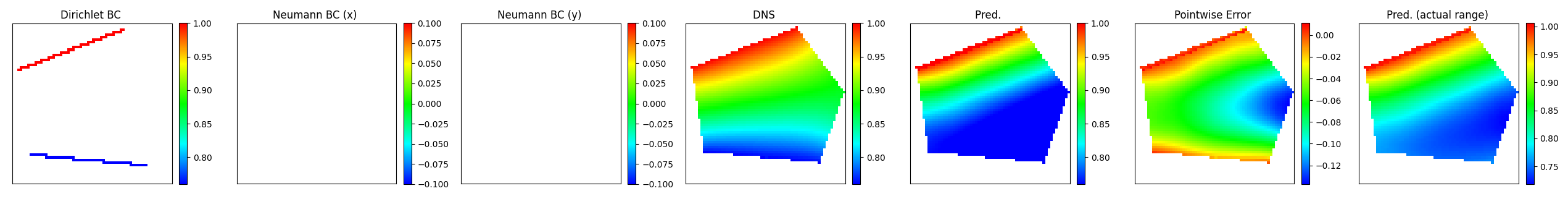}
  \includegraphics[width=0.9\linewidth]{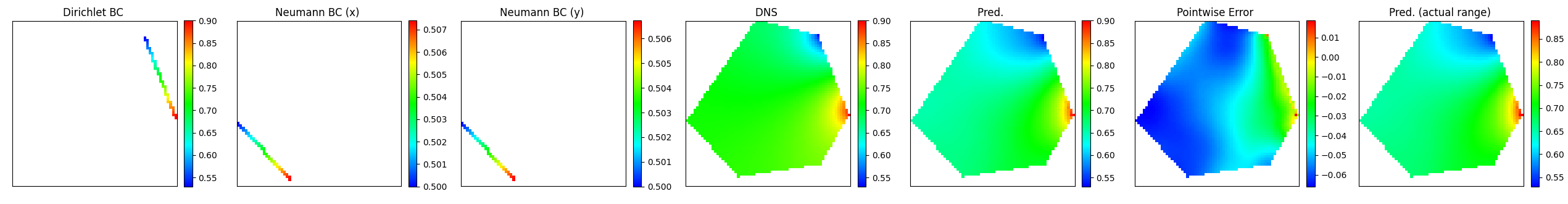}
  \includegraphics[width=0.9\linewidth]{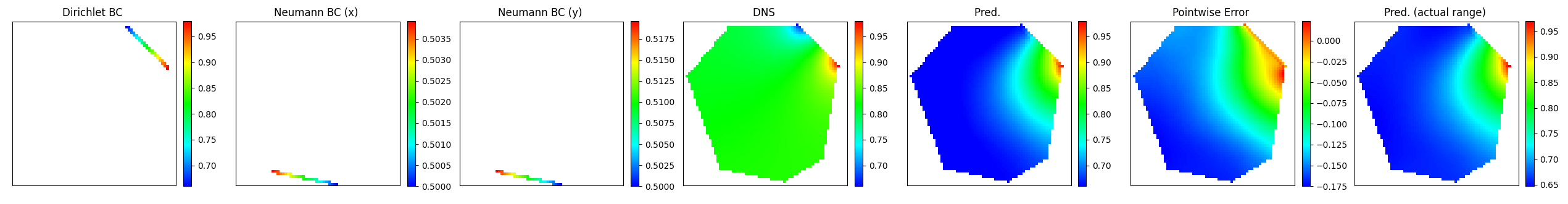}
  \includegraphics[width=0.9\linewidth]{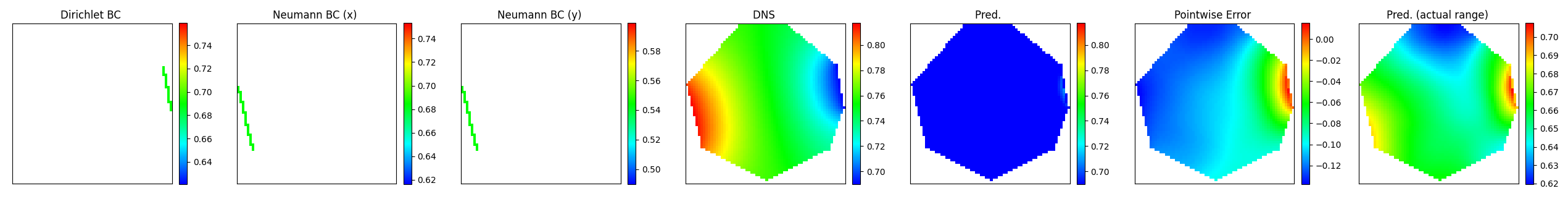}
  \caption{Poor results from NN predictions for T2 from the NN trained over D2. }
  \label{fig:D2-poor}
\end{figure}

\begin{figure}[h!]
  \centering
  \includegraphics[width=0.48\linewidth]{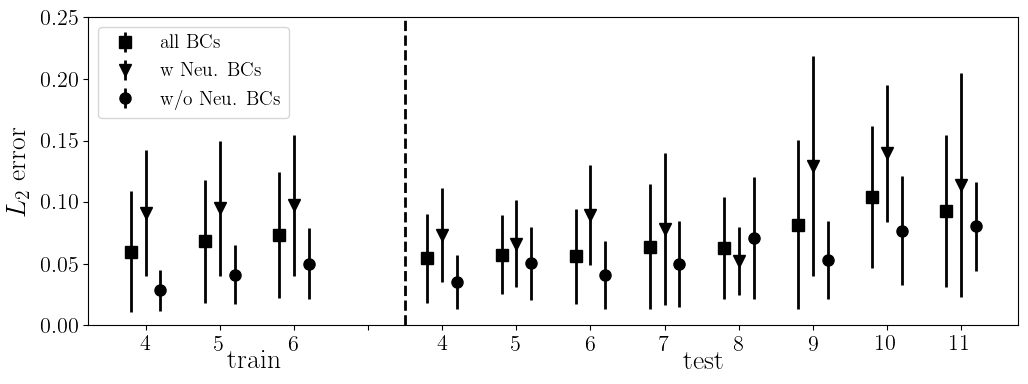}
  \includegraphics[width=0.48\linewidth]{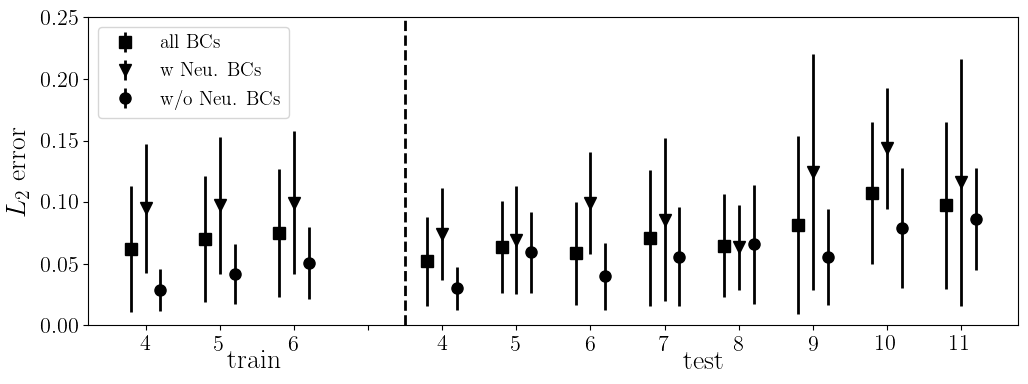} \\
  \includegraphics[width=0.48\linewidth]{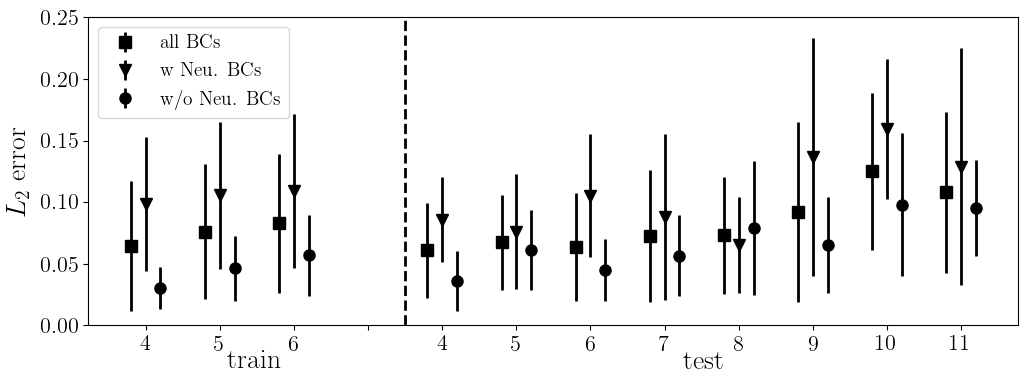}
  \includegraphics[width=0.48\linewidth]{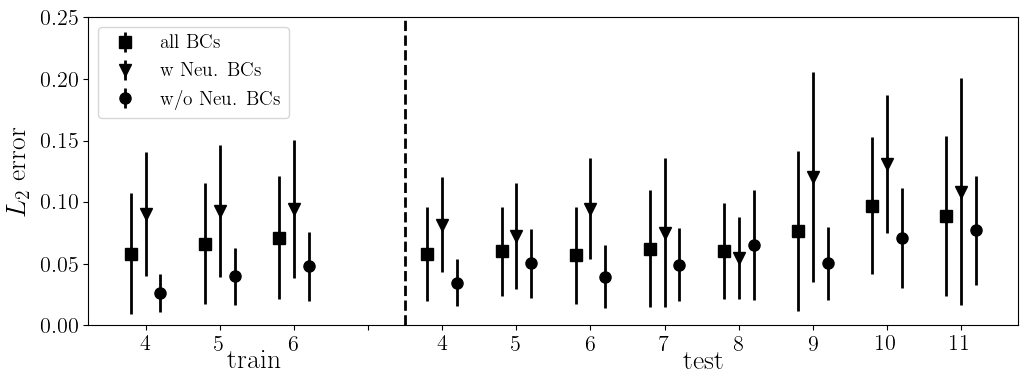}\\
  \includegraphics[width=0.48\linewidth]{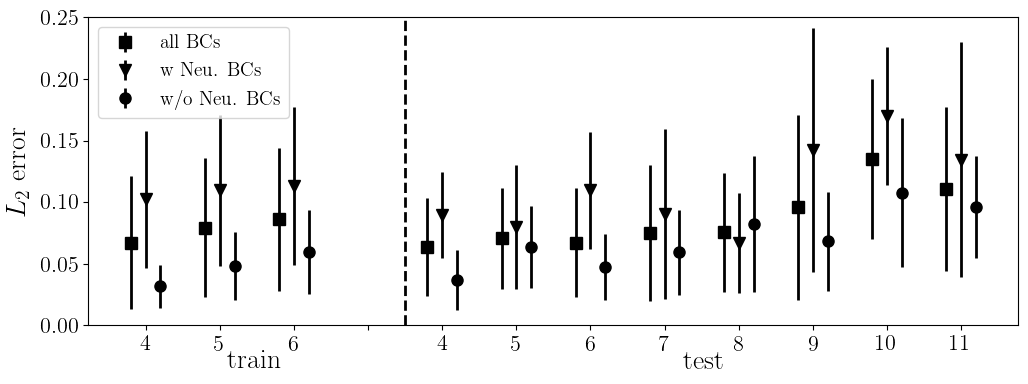}
  \includegraphics[width=0.48\linewidth]{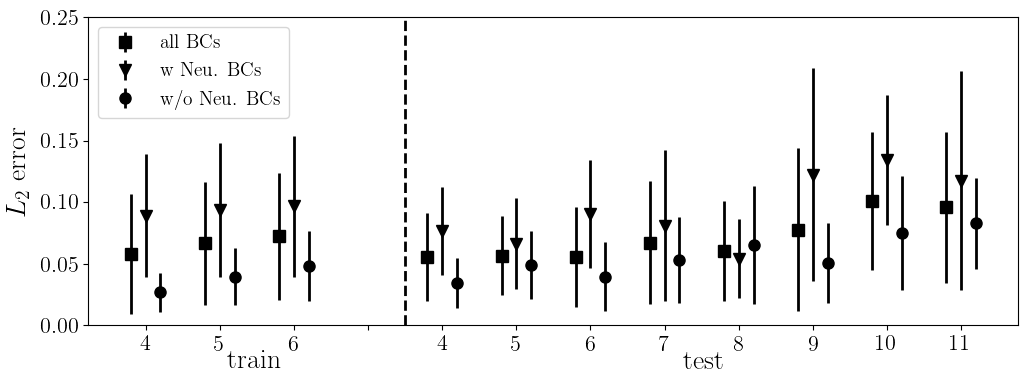}\\
  \includegraphics[width=0.48\linewidth]{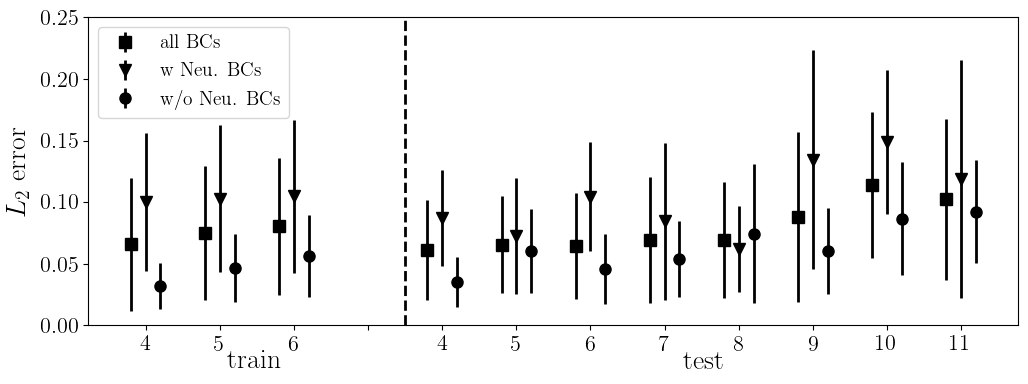}
  \caption{Statistics of NN results for different training cases to solve D2 to confirm the repeatability of the proposed method.}
  \label{fig:D2-NN-stats}
\end{figure}

\subsubsection{D3/T2: BVPs on quadrilateral/pentagon/hexagon/nonagon}
\begin{figure}[h!]
  \centering
  \includegraphics[width=1.0\linewidth]{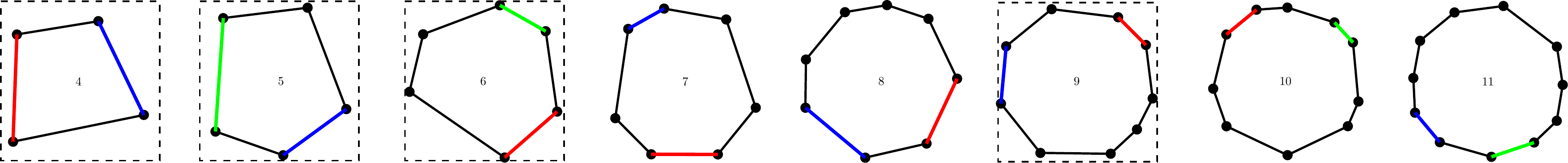}
  \caption{Illustration of dataset D3/T2. D3 contains BVPs that covers quadrilateral, pentagons, and hexagon. T2 contains polygons with the total number of edges ranging from four to eleven. Similar as other large datasets, the boundary values could have a constant/linear/quadratic/sinusoidal distribution.}
  \label{fig:large-D3}
\end{figure}

\begin{figure}[h!]
  \centering
  \includegraphics[height=45mm]{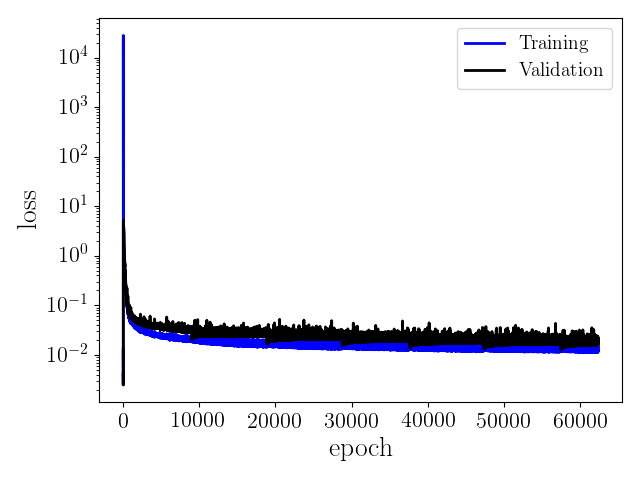}
  \includegraphics[height=45mm]{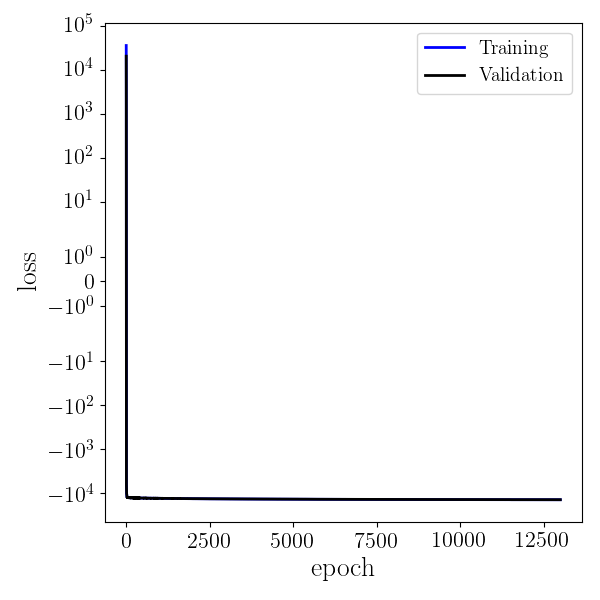}
  \includegraphics[height=45mm]{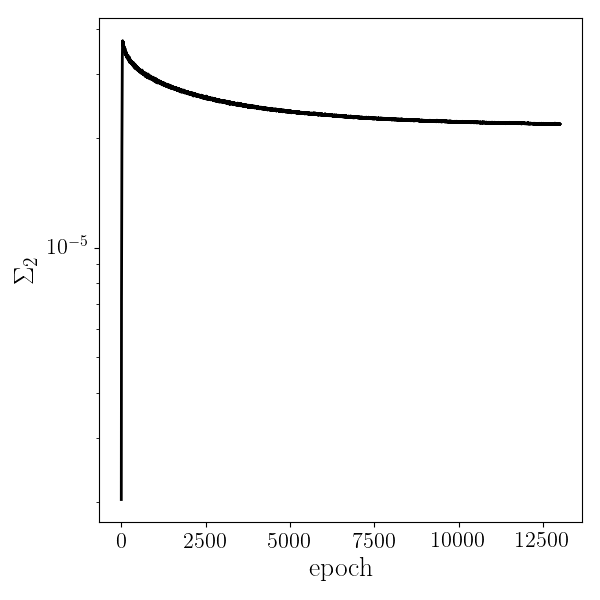}
  \caption{The training loss of the residual constrained NNs to solve all BVPs in D3. }
  \label{fig:D3-loss}
\end{figure}

\begin{figure}[h!]
  \centering
  \includegraphics[width=1.0\linewidth]{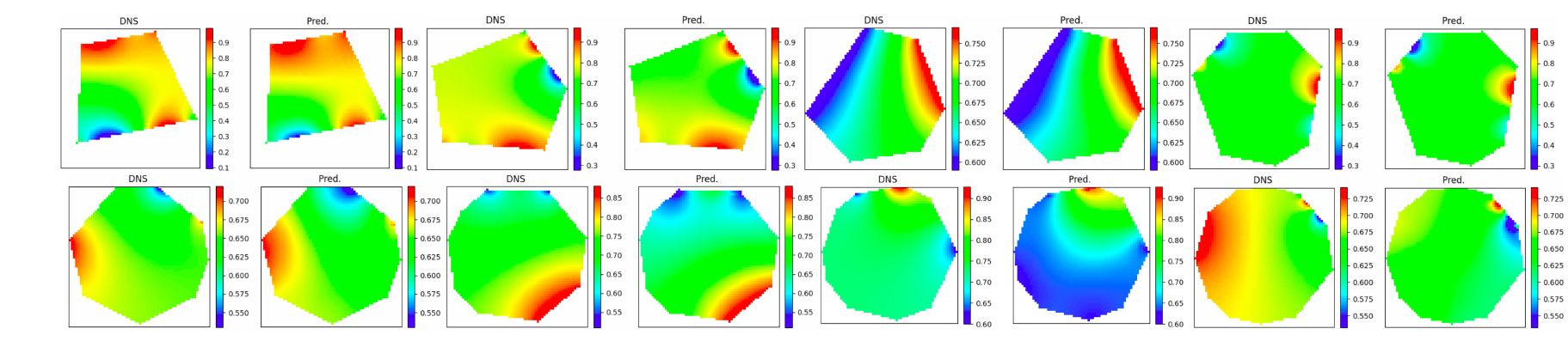}
  \caption{Selected good NN predictions for T2 (newly randomly generated polygons with different total number of edges) from the NN trained over D3. }
  \label{fig:D3-results}
\end{figure}

D3 contains 288K unique BVPs, which consists of 32K BVPs on quadrilateral, 64K BVPs on pentagons, 96K BVPs on hexagon and 96K BVPs on nonagon.
These BVPs covers 400 quadrilaterals, 400 pentagons, 400 hexagons, and 400 nonagons. 
The shape of these polygons are randomly generated. As for other large datasets, we apply boundary conditions to two non-adjacent edges, which are randomly sampled from all the possible locations. The boundary values are randomly generated, which could have a constant/linear/quadratic/sinusoidal distribution along the edges. 
The same testing dataset T2 are used for this problem. 
The training loss is shown in Fig. \ref{fig:D3-loss}.
NN architectures and training details are summarized in Table \ref{tab:D3-NNs} and \ref{tab:D3-NNs-others}. Selected NN predicted results are shown in Fig. \ref{fig:D3-results}(a), with the $L_2$ error of all the training and testing results given in Fig. \ref{fig:D3-results}(b).

\begin{table}
  \centering
  \begin{tabular}{l | l | l | l}
    \hline
    Deterministic         & Probabilistic         & Size         & Layer arguments \\ \hline
    Input                 & Input                 & -            & - \\
    LayerFillRandomNumber & LayerFillRandomNumber & -            & - \\
    Conv2D                & Convolution2DFlipout  & filters = 8  & kernel (5,5), padding: same, ReLU \\
    MaxPooling2D          & MaxPooling2D          & -            & kernel (2,2), padding: same\\
    Conv2D                & Convolution2DFlipout  & filters = 8 & kernel (5,5), padding: same, ReLU \\
    MaxPooling2D          & MaxPooling2D          & -            & kernel (2,2), padding: same\\
    Conv2D                & Convolution2DFlipout  & filters = 16 & kernel (5,5), padding: same, ReLU \\
    MaxPooling2D          & MaxPooling2D          & -            & kernel (2,2), padding: same\\
    Flatten               & Flatten               & -            & - \\
    Dense                 & DenseFlipout          & units = 32   & ReLU \\
    Dense                 & DenseFlipout          & units = 128   & ReLU \\
    Reshape               & Reshape               & -            & $[4,4,8]$ \\
    Conv2D                & Convolution2DFlipout  & filters = 16 & kernel (5,5), padding: same, ReLU \\
    UpSampling2D          & UpSampling2D          & -            & size (2,2) \\
    Conv2D                & Convolution2DFlipout  & filters = 16 & kernel (5,5), padding: same, ReLU \\
    UpSampling2D          & UpSampling2D          & -            & size (2,2) \\
    Conv2D                & Convolution2DFlipout  & filters = 16 & kernel (5,5), padding: same, ReLU \\
    UpSampling2D          & UpSampling2D          & -            & size (2,2) \\
    Conv2D                & Convolution2DFlipout  & filters = 16 & kernel (5,5), padding: same, ReLU \\
    Conv2D                & Convolution2DFlipout  & filters = 2  & kernel (5,5), padding: same, ReLU \\
    \hline
  \end{tabular}
  \caption{Details of both deterministic and probabilistic NNs for solving D3.}
  \label{tab:D3-NNs}
\end{table}

\begin{table}
  \centering
  \begin{tabular}{l | l | l }
    \hline
    Description                     & Deterministic                 & Probabilistic         \\ \hline
    Total parameters                & 41,346                        & 82,435                 \\
    Size of $\calD$                 & 5 $\times$ Aug: $2^{10}$      & 5 $\times$ Aug: $2^{9}$      \\
    Epochs                          & 10,000                        & 100                 \\
    Zero initialization epochs      & 100                           & -                     \\
    Optimizer                       & Adam                          & Nadam                 \\
    Learning Rate                   & 2.5e-4                        & 1e-8                  \\
    Batch Size                      & 256                           & 64                    \\
    $\Sigma_1$                      & -                             & 1e-8                  \\
    Initial value of $\Sigma_2$     & -                             & 1e-8                  \\
    \hline
  \end{tabular}
  \caption{Training related parameters for solving D3.}
  \label{tab:D3-NNs-others}
\end{table}

\begin{figure}[h!]
  \centering
  \includegraphics[width=0.9\linewidth]{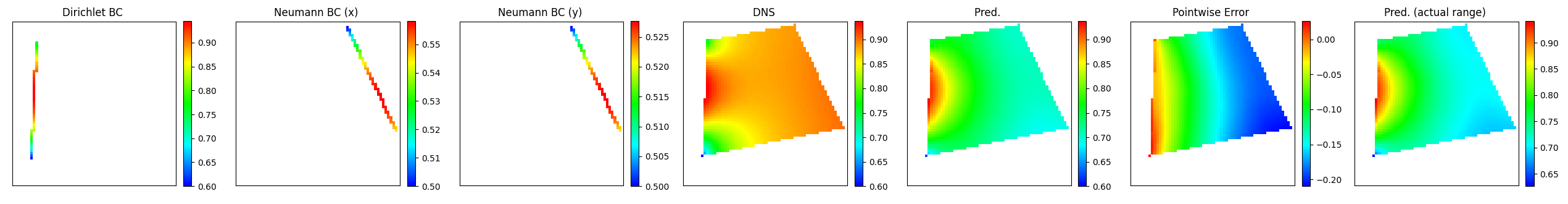}
  \includegraphics[width=0.9\linewidth]{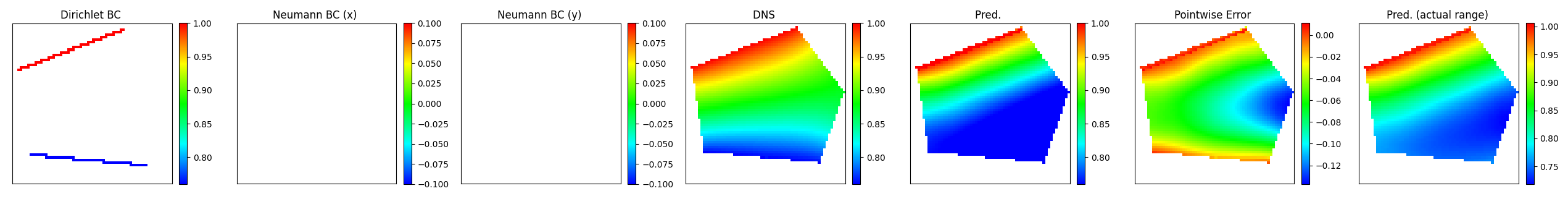}
  \includegraphics[width=0.9\linewidth]{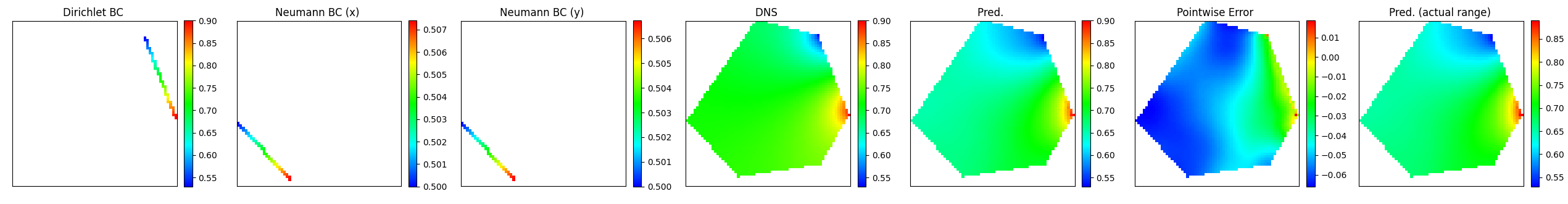}
  \includegraphics[width=0.9\linewidth]{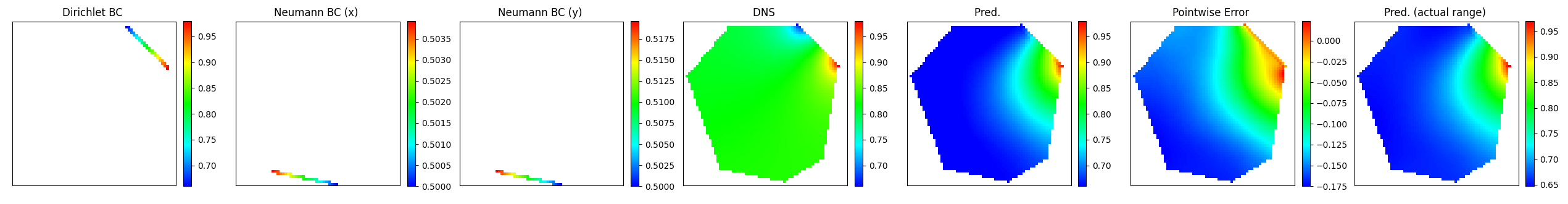}
  \includegraphics[width=0.9\linewidth]{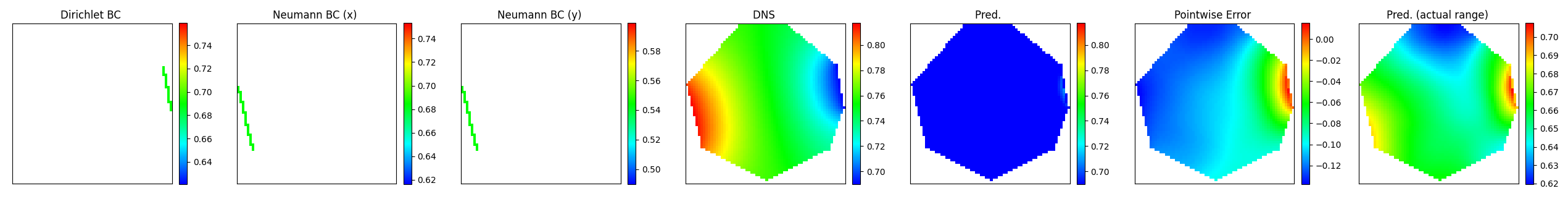}
  \caption{Poor results from NN predictions for T2 from the NN trained over D3. }
  \label{fig:D3-poor}
\end{figure}

\begin{figure}[h!]
  \centering
  \includegraphics[width=0.48\linewidth]{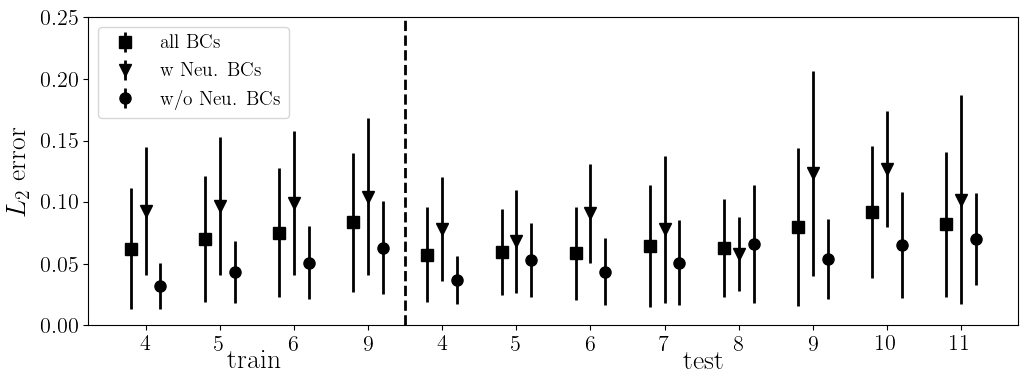}
  \includegraphics[width=0.48\linewidth]{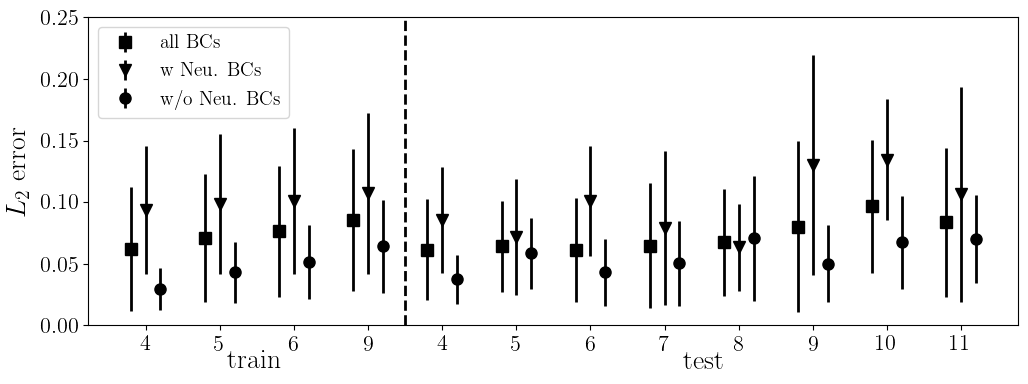} \\
  \includegraphics[width=0.48\linewidth]{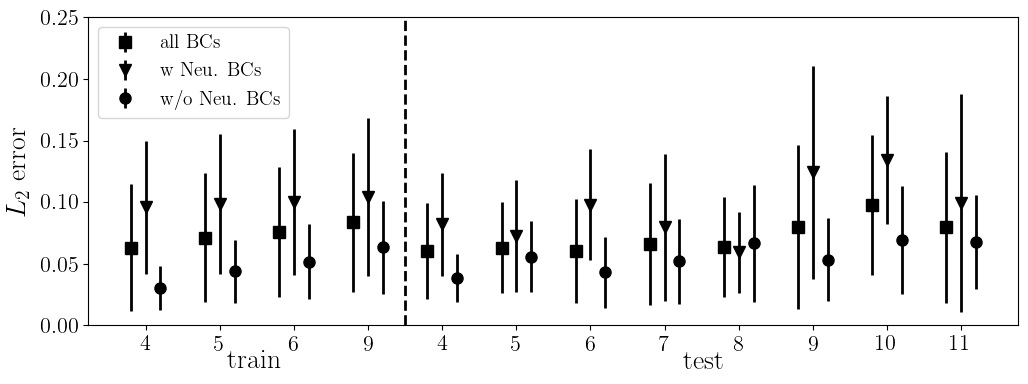}
  \caption{Statistics of NN results for different training cases to solve D3 to confirm the repeatability of the proposed method. }
  \label{fig:D3-NN-stats}
\end{figure}

\bibliographystyle{unsrt} 

\bibliography{lib.bib}